\documentclass[12pt]{article}

\usepackage{amstext    }
\usepackage{amsthm    }
\usepackage{a4}
\usepackage[mathscr]{eucal} 
\usepackage{mathrsfs}

\usepackage{amsmath}
\usepackage{amssymb}
\usepackage{amscd}

\numberwithin{equation}{section}

\newtheorem{theorem}{Theorem}[subsection]
\newtheorem{definition}[theorem]{Definition}
\newtheorem{proposition}[theorem]{Proposition}
\newtheorem{corollary}[theorem]{Corollary}
\newtheorem{lemma}[theorem]{Lemma}
\newtheorem{remark}[theorem]{Remark}
\newtheorem{remarks}[theorem]{Remarks}
\newtheorem{example}[theorem]{Example}

\newcommand{\cali}[1]{\mathscr{#1}}

\newcommand{\Aut}{{\rm Aut}}
\newcommand{\volume}{{\rm volume}}

\newcommand{\supp}{{\rm supp}}
\newcommand{\const}{{\rm const}}
\newcommand{\dist}{{\rm dist}}

\newcommand{\loc}{{loc}}
\newcommand{\ddc}{dd^c}
\newcommand{\dc}{d^c}
\newcommand{\dbar}{\overline\partial}
\newcommand{\ddbar}{\partial\overline\partial}
\newcommand{\DSH}{{\rm DSH}}

\newcommand{\id}{{\rm id}}
\newcommand{\capacity}{{\rm cap}}

\renewcommand{\Re}{{\rm Re}}

\newcommand{\Cc}{\cali{C}}

\newcommand{\Fc}{\cali{F}}

\newcommand{\Hc}{\cali{H}}
\newcommand{\Lc}{\cali{L}}
\newcommand{\Mc}{\cali{M}}

\newcommand{\Pc}{\cali{P}}
\newcommand{\Rc}{\cali{R}}
\newcommand{\Uc}{\cali{U}}
\newcommand{\Vc}{\cali{V}}

\newcommand{\FS}{{\rm FS}}

\newcommand{\C}{\mathbb{C}}

\newcommand{\R}{\mathbb{R}}
\renewcommand{\P}{\mathbb{P}}

\newcommand{\K}{{\cal K}}

\title{Super-potentials of positive closed currents, intersection theory
and dynamics}
\author{Tien-Cuong Dinh and Nessim Sibony}

\begin{document}
\maketitle

\begin{abstract}
We introduce a notion of super-potential for positive closed
currents of bidegree $(p,p)$ on projective spaces. 
This gives a calculus on positive closed currents of
arbitrary bidegree. We define in particular the intersection of
such currents and the pull-back operator by meromorphic maps. One of the
main tools is the introduction of structural discs in the space of
positive closed currents which gives a ``geometry'' on that space.
We apply the theory of super-potentials
to construct Green currents for rational maps and 
to study  equidistribution problems for holomorphic endomorphisms and for polynomial
automorphisms.
\end{abstract}




\tableofcontents

\section{Introduction} \label{section_introduction}

Let $(X,\omega)$ be a compact K{\"a}hler manifold. It is in general
quite difficult to develop a calculus on cycles of codimension $\geq
2$. An important approach has been introduced by Gillet-Soul{\'e}
\cite{GilletSoule} who constructed appropriate
potentials with tame singularities for cycles of
arbitrary codimension. See also Bost-Gillet-Soul\'e
\cite{BostGilletSoule}, Berndtsson \cite{Berndtsson} and
Henkin-Polyakov \cite{HenkinPolyakov} for the resolution of
$\ddbar$ and $\dbar$-equations in the projective space.
 
On the other hand, the calculus on positive closed currents of bidegree
$(1,1)$ using potentials is very useful and quite
well-developped. Demailly's paper \cite{Demailly2} and book
\cite{Demailly3} contain a clear exposition of this subject. It has
many applications in complex geometry and to holomorphic dynamics, see
the surveys \cite{Fornaess, Sibony} for background. The recent papers
\cite{DinhSibony2, DinhSibony3} by the authors give other applications.

Our main goal in this article is to develop a calculus on positive
closed currents of bidegree $(p,p)$. For simplicity, we restrict  here
to the case of the projective space $\P^k$. We first explain the
familiar situation of currents of bidegree $(1,1)$. The reader will find in
Paragraph \ref{section_geometry} some basic notions and  properties of
positive closed currents and of pluri-subharmonic functions.

Denote by $\omega$ the standard Fubini-Study form on $\P^k$ normalized by
$\int_{\P^k}\omega^k=1$. 
Let $S$ be a positive closed $(1,1)$-current on $\P^k$. 
We assume that the mass $\|S\|:=\langle
S,\omega^{k-1}\rangle$ is 1, that is, $S$ is cohomologous to $\omega$.
A {\it quasi-potential} of $S$ is a quasi-plurisubharmonic function
$u$ such that
$$S-\omega=\ddc u.$$
Recall that $\dc:={i\over {2\pi}}(\overline\partial-\partial)$. 
This function $u$ is unique when we normalize
it by  $\int_{\P^k} u\omega^k=0$. The correspondence $S\leftrightarrow u$ is very
useful. Indeed, $u$ has a value at every point if we allow the value
$-\infty$. This makes it possible to consider the pull-back of $S$ by 
dominant meromorphic maps \cite{Meo} or to consider the wedge-product (intersection) $S\wedge
S':=\omega\wedge S'+\ddc (uS')$ when $u$ is integrable with respect to the trace measure of a
positive closed current $S'$. 

From our point of view, the formalism in this case is as follows. Let
$\delta_x$ denote the Dirac mass at $x$. We consider a
$(k-1,k-1)$-current $v$, non uniquely determined, such that $\langle
v,\omega\rangle =0$ and $\ddc v=\delta_x-\omega^k$. We then have
formally
\begin{eqnarray*}
u(x) & = & \langle u,\delta_x\rangle = \langle
u,\delta_x-\omega^k\rangle =\langle u,\ddc v\rangle\\
& = & \langle \ddc u,v\rangle 
 =  \langle S-\omega,v\rangle = \langle S,v\rangle.
\end{eqnarray*}   
So, $\langle S,v\rangle$ is in particular independent of the choice of
$v$. Moreover, we can extend the action of $u$ to $\Cc_k$ the convex set
of probability measures. If $\ddc U_\nu=\nu-\omega^k$ with
$\nu\in\Cc_k$ and $\langle U_\nu,\omega\rangle=0$, we get
$$\langle u,\nu\rangle = \langle S, U_\nu\rangle,$$
where the value $-\infty$ is allowed. We prefer to consider that the
quasi-potential is acting on $\Cc_k$. Define
$$\Uc_S(\nu):=\langle u,\nu\rangle = \langle S, U_\nu\rangle.$$
This is somehow irrelevant in
this case since Dirac masses are the extremal points of $\Cc_k$ and 
$\Uc_S$ is simply the affine extension of $u$ to $\Cc_k$.

Let $\Cc_p$ denote the convex compact set of positive closed currents
$S$ of bidegree $(p,p)$ on $\P^k$ and of mass 1, i.e. $\|S\|:=\langle
S,\omega^{k-p}\rangle =1$. Let $U_S$ denote a solution to the equations
$$\ddc U_S = S-\omega^p, \qquad \langle U_S,\omega^{k-p+1}\rangle =0.$$
We introduce $\Uc_S$ as a function on
$\Cc_{k-p+1}$ that we will call {\it the super-potential of mean $0$
  of $S$}. Suppose $R$ is in $\Cc_{k-p+1}$ and let $\ddc
U_R=R-\omega^{k-p+1}$ with $\langle U_R,\omega^p\rangle=0$. Then, formally
\begin{eqnarray*}
\Uc_S(R) & := & \langle U_S,R\rangle = \langle
U_S,R-\omega^{k-p+1}\rangle  =  \langle U_S,\ddc U_R\rangle \\
&=&\langle \ddc U_S, U_R\rangle = \langle S-\omega^p,U_R\rangle = \langle S, U_R\rangle.
\end{eqnarray*}
The function $\Uc_S$ determines $S$. 
We will show that it is defined everywhere if the value $-\infty$ is allowed.

To develop the calculus, we have to consider $\Cc_p$ and $\Cc_{k-p+1}$
as infinite dimensional spaces with special families of currents
that we parametrize by the unit disc $\Delta$ in $\C$. We call these families
{\it special structural discs of currents}. 
When $\Uc_S$ is restricted to such discs we get
quasi-subharmonic functions. More precisely, if $x\mapsto
R_x$ is a special structural disc of currents parametrized by
$x\in\Delta$, then 
$$\ddc_x \Uc_S(R_x)\geq -\alpha$$
where $\alpha$ is a smooth $(1,1)$-form independent of $S$. The above definition of $\Uc_S(R)$ is
valid for $S$ or $R$ smooth. In general, we have
$$\Uc_S(R)=\lim_{x\rightarrow 0}\Uc_S(R_x)$$
for some special discs with $R_0=R$.

In Paragraph \ref{section_geometry}, we introduce a geometry on the space $\Cc_p$, in
particular, the structural varieties and their curvature forms $\alpha$.
In Paragraph \ref{section_sp}, we establish the basic properties of super-potentials,
in particular, convergence theorems which make the theory useful.
The main point is to extend the definition of the super-potential
$\Uc_S$ from smooth forms in $\Cc_{k-p+1}$ to arbitrary currents in
$\Cc_{k-p+1}$.
We introduce (Definition \ref{def_cv_h}) the notion of {\it Hartogs'
  convergence} (or {\it H-convergence} for short)  for currents, which is technically useful.
Paragraph \ref{section_intersection} deals with a theory of intersection of currents. We give
good conditions for the intersection of currents
of arbitrary bidegrees. 
Two currents $R_1\in\Cc_{p_1}$ and $R_2\in\Cc_{p_2}$ are {\it
  wedgeable} if and only if
a super-potential of $R_1$ is finite at $R_2\wedge
\omega^{k-p_1-p_2+1}$. The calculus on differential forms can be
extended to wedgeable currents: commutativity, associativity,
convergence and continuity of wedge-product for the H-convergence. 
If $R_2$ is of bidegree $(1,1)$, the condition means that the
quasi-potentials of $R_2$ are integrable with respect to the trace
measure of $R_1$.
As a special case, we obtain the usual
intersection of algebraic cycles. The question of developing such a
theory was raised by Demailly in \cite{Demailly2}. We give,  in the last paragraph, a
satisfactory approach to the problem of pulling back  a current in
$\Cc_p$ by meromorphic maps. Also in this
paragraph, we apply the theory of super-potentials to 
complex dynamics in higher dimension.
The main applications are
the following results.

As a first application, we construct Green currents of
bidegree $(p,p)$ for a large class of meromorphic maps on $\P^k$. This
requires a good calculus using the pull-back operation. The following
result holds for holomorphic maps and for Zariski generic meromorphic
maps which are not holomorphic. 

\begin{theorem}
Let $f$ be an algebraically $p$-stable meromorphic map on $\P^k$ with
dynamical degrees $d_s$, $1\leq s\leq k$. Assume that $d_{p-1}<d_p$ and that the union
of the infinite fibers is of dimension $\leq k-p$. 
Then,
$d_p^{-n}(f^n)^*(\omega^p)$ converge to an $f^*$-invariant current $T$ which is
is extremal among $f^*$-invariant currents in $\Cc_p$.
\end{theorem}

Note that the convergence result holds also for regular polynomial
automorphisms. 
The current $T$ is called {\it the Green current of bidegree $(p,p)$}
of $f$. The convergence is still valid if we replace $\omega^p$ by a
current with bounded super-potentials. The case $p=1$
was considered by the second author in \cite{Sibony}.

Let $\Mc_d(\P^k)$  denote the space 
of dominant meromorphic self-maps of algebraic degree $d\geq 2$ on
$\P^k$.  Such a map can be lifted to a homogeneous polynomial self-map
of $\C^{k+1}$ of degree $d$. The lift is unique up to a multiplicative
constant.   
The space 
$\Mc_d(\P^k)$ has the structure of a Zariski dense open set in $\P^N$ with
$N:=(k+1)(d+k)!/(d!k!) -1$.
The space $\Hc_d(\P^k)$ 
of holomorphic self-maps of algebraic degree $d\geq 2$ on
$\P^k$ is a Zariski open subset of  $\Mc_d(\P^k)$ and
 $\Mc_d(\P^k)\setminus \Hc_d(\P^k)$ is an irreducible hypersurface of
 $\Mc_d(\P^k)$, see \cite{BassanelliBerteloot} and \cite[p.427]{GuelfandKapranov}.

\begin{theorem} There is a Zariski dense open set $\Hc^*_d(\P^k)$ in
  $\Hc_d(\P^k)$ such that if $f$ is in $\Hc_d^*(\P^k)$ and if $S$ is a
  current in $\Cc_p$, then $d^{-pn}(f^n)^*(S)$ converges to the Green
  current of bidegree $(p,p)$ of $f$ uniformly on $S$.
\end{theorem}

A more precise description is known for $p=1$ and $k=2$ in 
\cite{FornaessSibony3, FavreJonsson}, for $p=1$ and $k\geq 2$ in \cite{DinhSibony9}, and
for $p=k$ in \cite{DinhSibony1, DinhSibony9}, see also
\cite{FornaessSibony1, BriendDuval}. 
Applying the previous theorem to the currents of integration on 
subvarieties $H$ gives the equidistribution of $f^{-n}(H)$ in $\P^k$.
Another application is a rigidity theorem for polynomial automorphisms
of $\C^k$ that we consider as birational maps on $\P^k$. 
 
\begin{theorem} \label{th_2}
Let $f$ be a polynomial automorphism of $\C^k$ which is regular in the
sense of \cite{Sibony}. Let $I_+$ denote the indeterminacy
set of $f$ at infinity and $p$ the integer such that $\dim
I_+=k-p-1$. Let $\K_+$ be the set of points
$z\in\C^k$ with bounded orbits. Then, 
the Green $(p,p)$-current associated to $f$ is the
unique positive closed $(p,p)$-current of mass $1$ with support in
$\overline\K_+$.
\end{theorem}

The result was proved by Forn\ae ss and the second author in 
dimension $k=2$ \cite{FornaessSibony1}. 
Note that when $k=2$ and $p=1$, regular automorphisms are the H{\'e}non type
automorphisms of $\C^2$. It is known that dynamically interesting
polynomial automorphisms in $\C^2$ are conjugated to the regular ones \cite{FriedlandMilnor}.
Let $H$ be an
analytic subset of pure dimension $k-p$ which does not intersect the
indeterminacy set $I_-$ of $f^{-1}$. We obtain as a consequence
of Theorem \ref{th_2} that the currents of
integration on $f^{-n}(H)$, properly normalized, converge to the Green
$(p,p)$-current of $f$. The case $k=2$ and $p=1$ of this result was
proved by Bedford-Smillie in \cite{BedfordSmillie}.

\begin{remark}\rm
The super-potential $\Uc_S$ can be extended to a function on weakly
positive closed currents of bidegree $(k-p+1,k-p+1)$. For simplicity,
we consider only (strongly) positive currents. We can also define
super-potentials for weakly positive closed $(p,p)$-currents; they are
functions on (strongly) positive closed currents  of bidegree
$(k-p+1,k-p+1)$.
The super-potentials are introduced on currents of mass 1 but they can
be easily extended by linearity to currents of arbitrary mass. Their
domain of definition can be also extended to positive closed currents of
arbitrary mass. 
\end{remark}

\noindent
{\bf Other notation.}
$\Delta_r$ is the disc of center 0 and of radius $r$ in $\C$, 
$\Delta$ denotes the unit disc, $\Delta^k$ the unit polydisc in $\C^k$
and $\Delta^*:=\Delta\setminus\{0\}$. 
The group of automorphisms of $\P^k$ is 
a complex Lie group of dimension $k^2+2k$ that we
denote by
$\Aut(\P^k)\simeq {\rm PGL}(k+1,\C)$.
We will
work with a fixed holomorphic chart 
and local holomorphic coordinates $y$ of $\Aut(\P^k)$. 
The automorphism with coordinates
$y$ is denoted by $\tau_y$.
Choose
$y$ so that $|y|< 2$ and $y=0$ at the identity $\id\in\Aut(\P^k)$.
In order to simplify the notation, choose a norm $|y|$ of $y$ which
is invariant under the
involution $\tau\mapsto\tau^{-1}$.
Fix a smooth probability measure $\rho$
with compact support in $\{|y|<1\}$. 
Choose $\rho$ radial and decreasing when $|y|$ increases.
 So, the involution $\tau\mapsto \tau^{-1}$
preserves $\rho$.
The {\it mass} of a positive or negative $(p,p)$-current $S$ on $\P^k$
is defined by $\|S\|:=|\langle S,\omega^{k-p}\rangle|$. Throughout the
paper, $S_\theta$, $R_\theta$, $\ldots$ will denote the
regularization of $S$, $R$, $\ldots$ defined in Paragraph
\ref{section_topology} below.

\bigskip
\noindent
{\bf Aknowledgement.} We thank the referee who has read carefully the
first version of this paper. He suggested several clarifications which
permited to improve the exposition.


\section{Geometry of currents on projective spaces} \label{section_geometry}

In this paragraph, 
we introduce some basic facts about the convex set $\Cc_p$ of positive closed
$(p,p)$-currents of mass 1 in $\P^k$.


\subsection{Topology and distances on the spaces of currents} \label{section_topology}

Let $X$ be a complex manifold of dimension $k$. Recall that a
$(p,p)$-form $\Phi$ on $X$ is {\it (strongly) positive} if it is
positive at every point $a\in X$, that is, $\Phi$ is equal at
the point $a$ to a linear combination with positive coefficients of forms of type
$$(i\varphi_1\wedge \overline \varphi_1) \wedge \ldots \wedge
(i\varphi_p\wedge \overline \varphi_p)$$
where $\varphi_i$ are $(1,0)$-forms on $X$. Positive
$(0,0)$-forms are positive functions and positive 
$(k,k)$-forms are products of volume forms with 
positive functions.

A $(p,p)$-form $\Phi$ is {\it weakly positive} if $\Phi\wedge \Psi$ is
a positive form of maximal bidegree for every positive
$(k-p,k-p)$-form $\Psi$. 
A $(p,p)$-current
$T$ on $X$ is {\it positive} (resp. {\it weakly positive}) if $T\wedge
\Psi$ is a positive measure for every weakly positive (resp. positive)
smooth $(k-p,k-p)$-form $\Psi$. 
Positive forms and currents are weakly positive.   
The notions of positivity and of weak positivity
coincide only for bidegrees $(0,0)$, $(1,1)$, $(k-1,k-1)$ and $(k,k)$.  
We also say that $\Phi$ and $T$ are {\it negative} or
{\it weakly negative} if $-\Phi$ and $-T$ are positive or weakly
positive. For real $(p,p)$-currents $T,T'$, we will write $T\geq T'$
and $T'\leq T$ when $T-T'$ is positive.

Assume that $X$ is a compact K{\"a}hler manifold and $\omega_X$ is a
K{\"a}hler form on $X$. If $T$ is a positive or
negative $(p,p)$-current, {\it the mass} of $T$ on a Borel set $K\subset X$ is
the mass of {\it the trace measure} $T\wedge \omega_X^{k-p}$ of $T$ on $K$; that is
$$\|T\|_K:=|\langle T,\omega_X^{k-p}\rangle_K|.$$  
{\it The mass} of $T$ means its mass $\|T\|$ on $K=X$. Assume that $T$ is
positive and closed. Then, $\|T\|$ depends only on the class of $T$ in the 
Hodge cohomology group $H^{p,p}(X,\C)$. We recall the notion of density
of positive closed currents. Let $x$ denote 
local coordinates in a neighbourhood of a point $a\in X$ such that
$x=0$ at $a$ and $\beta:=\ddc |x|^2$ denote the standard Euclidean form. 
Let $B_r$ denote the ball $\{|x|<r\}$. {\it The Lelong number} of $T$ at $a$ is defined by
$$\nu(T,a):=\lim_{r\rightarrow 0}
\frac{\|T\wedge \beta^{k-p}\|_{B_r}}{\pi^{k-p}r^{2k-2p}}.$$
When $r$ decreases to 0, the expression on the right hand side
decreases to $\nu(T, a)$ which does not depend on the choice of
coordinates $x$ \cite{Siu}. The Lelong number compares the mass of the current on $B_r$ with
the  Euclidean volume $\pi^{k-p}r^{2k-2p}/(k-p)!$ of a ball of radius
$r$ in $\C^{k-p}$. 
A theorem of Siu says that
$\{\nu(T,a)\geq c\}$ is an analytic subset of dimension $\leq
k-p$ of $X$ for every $c>0$ \cite{Siu}.

The K{\"a}hler manifolds we consider in this paper are the projective space $\P^k$ and the
product $\P^k\times \P^k$. Let $\pi_1$ and $\pi_2$ be the canonical
projections of $\P^k\times\P^k$ onto its factors.
Let $\omega$ denote the Fubini-Study form
on $\P^k$ normalized so that $\int_{\P^k} \omega^k=1$, and define
$$\widetilde \omega:= \pi_1^*(\omega)+\pi_2^*(\omega)$$
the canonical K{\"a}hler form on $\P^k\times\P^k$.
If $T$ is a positive closed $(p,p)$-current on $\P^k$,  
one proves easily that $\nu(T,a)\leq \|T\|$ for every $a\in\P^k$. 

\begin{example} \label{ex_Lelong}
\rm
Let $V$ be an analytic subset of pure dimension
$k-p$ in $\P^k$. Lelong showed in \cite{Lelong} that the integration on the
regular part of $V$ defines a positive closed $(p,p)$-current
$[V]$. The mass of $[V]$ is equal to the degree of $V$, i.e. the
number of points in the intersection of $V$ with a generic projective
plane $P$ of dimension $p$. By a theorem of Thie, the Lelong number of
$[V]$ at $a$ is the multiplicity of $V$ at $a$, i.e. the multiplicity
at $a$ of $V\cap P$ for $P$ generic passing through $a$. 
This number is also equal to the number of
points, in a small neighbourhood of $a$, of $V\cap P'$
for $P'$ generic close enough to $P$. We deduce from the definition of
the Lelong number that there are constants $c,c'>0$ such that 
$$cr^{2k-2}\leq \volume(V\cap B)\leq c'r^{2k-2}$$
for every ball $B$ with center in $V$ of radius $r\leq 1$.
\end{example}

We will use the {\it weak topology} in $\Cc_p$,
i.e. the topology induced
by the weak topology of currents. Recall that a sequence of
$(p,p)$-currents $(R_n)$ converges weakly to a current $R$ if $\langle R_n,\Phi\rangle\rightarrow \langle
R,\Phi\rangle$ for every smooth $(k-p,k-p)$-form $\Phi$ on $\P^k$. 
Since the currents in $\Cc_p$ are positive, we obtain the same topology on
$\Cc_p$ if we consider real continuous forms $\Phi$ instead of smooth forms.
For this topology, $\Cc_p$ is compact.

We introduce some natural {\it distances} on $\Cc_p$ as
follows. 
For $\alpha\geq 0$ let $[\alpha]$ denote the integer part of $\alpha$. Let
$\Cc_{p,q}^\alpha$ be the space of $(p,q)$-forms whose coefficients admit
derivatives of all orders $\leq [\alpha]$ and these derivatives are $(\alpha-[\alpha])$-H{\"o}lder continuous.  
We use here the sum of $\Cc^\alpha$-norms of the coefficients for a
fixed atlas.
If $R$ and $R'$ are
currents in $\Cc_p$, define 
$$\dist_\alpha(R,R'):=\sup_{\|\Phi\|_{\Cc^\alpha}\leq 1} |\langle R-R',\Phi\rangle|$$
where $\Phi$ is a smooth $(k-p,k-p)$-form on $\P^k$.
Observe that $\Cc_p$ has finite diameter with respect to these
distances since $\langle R,\Phi\rangle$ and $\langle R',\Phi\rangle$
are bounded.

\begin{lemma} \label{lemma_compare_dist}
For every $0<\alpha<\beta<\infty$, there is a constant
$c_{\alpha,\beta}>0$ such that 
$$\dist_\beta\leq\dist_\alpha\leq
c_{\alpha,\beta}[\dist_\beta]^{\alpha/\beta}.$$
In particular, a function on $\Cc_p$ is H{\"o}lder continuous for $\dist_\alpha$
if and only if it is H{\"o}lder continuous for $\dist_\beta$.
\end{lemma}
\proof
The first inequality is clear. Let
$L:\Cc^\infty_{k-p,k-p}\rightarrow \C$ be a continuous linear
form. Assume that there are constants $A$ and $B$ such that 
$|L(\Phi)|\leq A\|\Phi\|_{\Cc^0}$ and $|L(\Phi)|\leq
B\|\Phi\|_{\Cc^\beta}$. The
theory of interpolation between Banach spaces \cite{Triebel}
implies that $|L(\Phi)|\leq c_{\alpha,\beta} A^{1-\alpha/\beta}B^{\alpha/\beta}\|\Phi\|_{\Cc^\alpha}$ with
$c_{\alpha,\beta}$ independent of $A$, $B$ and $L$. Applying this to
$L:=R-R'$ with $R,R'$ as above, gives the second
inequality in the lemma. 
\endproof

When $p=k$,  $\Cc_k$ is the convex of the probability measures on
$\P^k$ and its extremal elements are the Dirac masses. One can
identify the set of extremal elements of $\Cc_k$ with $\P^k$. Let
$\delta_a, \delta_b$ denote the Dirac masses at $a, b$ and
$\|a-b\|$ the distance between $a$, $b$  induced by the Fubini-Study metric.

\begin{lemma} We have
$$\dist_\alpha(\delta_a,\delta_b)\simeq
\|a-b\|^{\min\{\alpha,1\}}.$$
\end{lemma}
\proof
It is enough to consider the case where $a$ and $b$ are close. Let $x=(x_1,\ldots,x_k)$
be local coordinates so that $a$ and $b$ are close to 0. 
Without loss of generality, one can assume $a=0$ and
$b=(t,0,\ldots,0)$. 
It is clear that
$$\dist_\alpha(\delta_a,\delta_b)=\sup_{\|\Phi\|_{\Cc^\alpha}\leq 1}
|\Phi(a)-\Phi(b)|\lesssim \|a-b\|^{\min\{\alpha,1\}}.$$
Using a cut-off function, one construct easily a function $\Phi$ with
bounded $\Cc^\alpha$-norm such that near 0, $\Phi(x)=|\Re(x_1)|^\alpha$ if
$\alpha<1$ and $\Phi(x)=\Re(x_1)$ if $\alpha\geq 1$.
Hence
 $$\dist_\alpha(\delta_a,\delta_b)\gtrsim |\Phi(a)-\Phi(b)|=
 \|a-b\|^{\min\{\alpha,1\}}.$$
This implies the lemma.
\endproof

\begin{proposition}
For $\alpha>0$,  the topology  induced by 
$\dist_\alpha$ coincides with the weak topology on $\Cc_p$. In
particular, $\Cc_p$ is a compact separable metric space.
\end{proposition}
\proof
It is clear that the convergence with respect to $\dist_\alpha$
implies the weak convergence. Conversely, if a sequence converges
weakly in $\Cc_p$, then it converges uniformly on compact sets of test
forms with uniform norm. By Dini's theorem, the set of test forms $\Phi$
with $\|\Phi\|_{\Cc^\alpha} \leq 1$ is relatively compact for the uniform convergence. The
proposition follows.
\endproof

Note that since the convex set $\Cc_p$ is a Polish space, measure theory on $\Cc_p$ is
quite simple.
We show in Lemma \ref{lemma_deform} and Proposition
\ref{prop_regularization} 
below that smooth forms are dense in $\Cc_p$, see
\cite{DinhSibony3} 
for the case of arbitrary compact K{\"a}hler manifolds. 
Here, since $\P^k$ is homogeneous, one can use the group $\Aut(\P^k)$ of
automorphisms of $\P^k$ in order to regularize currents, see also
\cite{deRham, Demailly3}.

Let $h_\theta(y):=\theta y$ denote the multiplication by $\theta\in\C$
and for $|\theta|\leq 1$ define $\rho_\theta:=(h_\theta)_*\rho$, see
Introduction for the notation.
Then, $\rho_0$ is the Dirac mass at the identity $\id\in\Aut(\P^k)$ and $\rho_\theta$ is a
smooth probability measure if $\theta\not=0$. Moreover, for every
$\alpha\geq 0$ there is a
constant $c_\alpha>0$ such that
$$\|\rho_\theta\|_{\Cc^\alpha}\leq c_\alpha
|\theta|^{-2k^2-4k-\alpha}$$
where $2k^2+4k$ is the real dimension of $\Aut(\P^k)$.
Define
for any positive or negative $(p,p)$-current $R$ on $\P^k$ not necessarily closed
$$R_\theta:=\int_{\Aut(\P^k)} (\tau_y)_* R
\ d\rho_\theta(y)=\int_{\Aut(\P^k)}(\tau_{\theta y})_*R
\ d\rho(y)=\int_{\Aut(\P^k)}(\tau_{\theta y})^*R\ d\rho(y).$$
The last equality follows from the fact that $\rho$ is radial and the
involution $\tau\mapsto \tau^{-1}$ preserves the norm of $y$.

Define $R_{\theta y}:=(\tau_{\theta
  y})_*R$. If $R$ is positive and closed, then $R_{\theta y}$ and $R_\theta$
are also positive and closed. Observe that since $\rho$ is radial,
$R_\theta=R_{\theta'}$ when $|\theta|=|\theta'|$.

\begin{lemma} \label{lemma_deform}
When $\theta$ tends to $0$, $R_{\theta y}$ and
  $R_\theta$ converge weakly to $R$. 
If the restriction of $R$ to an open set $W\subset \P^k$ is a form of
class $\Cc^\alpha$, then $R_{\theta y}$ and $R_\theta$ converge to $R$
in $\Cc^\alpha(W')$ for any $W'\Subset W$.
\end{lemma}
\proof
The convergence of $R_{\theta y}$ is deduced  from the fact that
$\tau_{\theta y}$ converge to the identity in the $\Cc^\infty$
topology. This and the definition of $R_\theta$ imply the convergence of $R_\theta$.
\endproof

\begin{proposition} \label{prop_regularization}
If $\theta\not=0$, then $R_\theta$ is a
smooth form which depends continuously on $R$. Moreover, for every $\alpha\geq 0$ there is a constant $c_\alpha$
independent of $R$ such that
$$\|R_\theta\|_{\Cc^\alpha}\leq c_\alpha \|R\|
|\theta|^{-2k^2-4k-\alpha}.$$
If $K$ is a compact set in $\Delta^*$, there is a constant $c_{\alpha,K}>0$
such that if $\theta$ and $\theta'$ are in $K$ then
$$\|R_\theta-R_{\theta'}\|_{\Cc^\alpha}\leq c_{\alpha,K} \|R\| |\theta-\theta'|.$$
\end{proposition}
\proof
We can assume that $R$ is supported at a point $a$, that is,
$R=\delta_a\wedge \Psi$ for some tangent $(k-p,k-p)$-vector $\Psi$ defined at $a$ with
norm $\leq 1$ (here, we use Federer's notation and we consider the
vector $\Psi$ as a form with negative bidegree $(p-k,p-k)$). 
The general case is deduced using a desintegration of
$R$ as currents with support at a point. We have
$$R_\theta=\int_{\Aut(\P^k)}\big(\delta_{\tau_y(a)}\wedge
(\tau_y)_*\Psi\big) d\rho_\theta(y).$$
Hence, $R_\theta$ is smooth and depends continuously on $R$. The estimate on
$\|R_\theta\|_{\Cc^\alpha}$ follows from the estimate on the
$\Cc^\alpha$-norm of $\rho_\theta$. The last estimate in the
proposition follows from the inequality
$\|\rho_\theta-\rho_{\theta'}\|_{\Cc^\alpha} \lesssim
|\theta-\theta'|$ on $K$.
\endproof

\begin{remark}\rm \label{rk_theta_reg}
We call $R_\theta$ {\it the $\theta$-regularization} of $R$. 
In Proposition \ref{prop_regularization} we can 
replace $|\theta|^{-2k^2-4k-\alpha}$ by $|\theta|^{-2k-\alpha}$
but the estimates are more technical. 
\end{remark}

Let $\dist(\tau,\tau')$ denote the
distance between $\tau$ and $\tau'$ for a fixed smooth metric on
$\Aut(\P^k)$. The following simple lemma will be useful in the next paragraphs.

\begin{lemma} \label{lemma_deformation_c1} 
Let $K$ be a compact subset of $\Aut(\P^k)$. Let $W$ and
$W_0$ be  open sets in $\P^k$ such that $\overline W_0\subset \tau(W)$ for
every $\tau\in K$.  If $R$ is of class $\Cc^\alpha$, $\alpha\geq 0$, on $W$, then $\tau_*(R)$ is of class $\Cc^\alpha$
on $W_0$. Moreover, there is a constant $c>0$ such that
for  all $\tau$ and $\tau'$ in $K$
$$\|\tau_*(R)\|_{\Cc^\alpha(W_0)}\leq c \|R\|_{\Cc^\alpha(W)}$$
and
$$ \|\tau_*(R) -\tau'_*(R)\|_{\Cc^\alpha(W_0)}\leq
c\|R\|_{\Cc^\alpha(W)}\dist(\tau,\tau')^{\min(\alpha,1)}.$$
\end{lemma}
\proof
Since $\overline W_0\subset \tau(W)$, it is clear that  $\tau_*(R)$ is of class $\Cc^\alpha$
on $W_0$.
For $\tau$ in $K$, we have 
$\|\tau^{-1}\|_{\Cc^{\alpha+1}}\leq A$ which implies the first estimate.
For the second one, observe that
$$\tau_*(R)-\tau'_*(R)=\tau_*\big[R-\tau^*\tau'_*(R)\big]=\tau_*\big[R-(\tau^{-1}\circ\tau')_*(R)\big].$$ 
This and the inequality 
$$\|\tau^{-1}\circ\tau'-\id\|_{\Cc^{\alpha+1}}\lesssim
\dist(\tau,\tau')$$
imply the estimate.
\endproof


\subsection{Quasi-plurisubharmonic functions and capacity}  \label{section_quasi_psh}

Positive closed currents of bidegree $(1,1)$ admit
quasi-potentials which are quasi-plurisubharmonic functions (quasi-psh
for short). The compactness properties of these functions are fundamental
in the study of positive closed $(1,1)$-currents.
We recall here some facts, see \cite{Demailly3, DinhSibony6}.

A {\it quasi-psh function} is locally the difference of a psh function and a
smooth one, see \cite{Demailly3}. The first important property we
will use is the following that we state only in dimension 1. It is a
direct consequence of \cite[Theorem 4.4.5]{Hormander}.

\begin{lemma} \label{lemma_exp_estimate}
Let $\Fc$ be a compact family in $\Lc^1_\loc(\Delta)$ of subharmonic functions on
$\Delta$. Then, for every compact subset $K\subset \Delta$ there are
constants $c>0$ and $A>0$ such that 
$$\|e^{-A u}\|_{\Lc^1(K)}\leq c\quad \mbox{for every}\quad u\in\Fc.$$ 
\end{lemma}

Recall that 
a function $\varphi:\P^k\rightarrow \R\cup\{-\infty \}$ is
quasi-psh if and only if 
\begin{enumerate}
\item[$\bullet$] $\varphi$ is integrable with respect to the Lebesgue measure
and $\ddc\varphi\geq -c\omega$ for
some constant $c>0$;
\item[$\bullet$] $\varphi$ is strongly upper semi-continuous (strongly
  u.s.c. for short),
that is, for any Borel subset $A\subset \P^k$ of full Lebesgue
measure, we have $\varphi(x)=\limsup_{y\rightarrow x} \varphi(y)$ with $y\in A\setminus\{x\}$.
\end{enumerate}
A set $E\subset\P^k$ is {\it pluripolar}
or {\it complete pluripolar} if there is a quasi-psh function
$\varphi$ such that $E\subset \varphi^{-1}(-\infty)$ or
$E=\varphi^{-1}(-\infty)$ respectively.

If $\varphi$ is as above, then the $(1,1)$-current $T:=\ddc\varphi+c\omega$ is positive closed and of
mass $c$ since it is cohomologous to $c\omega$. We say that $\varphi$ is a {\it quasi-potential} of
$T$; it is defined everywhere on $\P^k$.
There is a continuous 1-1 correspondence between the positive closed $(1,1)$-currents
of mass 1 and the quasi-psh functions $\varphi$ satisfying $\ddc\varphi\geq
-\omega$, normalized by $\int_{\P^k} \varphi\omega^k=0$ or by
$\max_{\P^k}\varphi=0$. The following compactness property is deduced
from the corresponding properties of psh functions.

\begin{proposition}
Let $(\varphi_n)$ be a sequence of quasi-psh functions on $\P^k$ with
$\ddc\varphi_n\geq -\omega$. Assume that $\varphi_n$ is bounded from
above by a constant independent of $n$. Then, either $(\varphi_n)$ converges
uniformly to $-\infty$ or there is a subsequence $(\varphi_{n_i})$
converging, in $\Lc^p$ for $1\leq p<\infty$, to a quasi-psh function $\varphi$ with $\ddc\varphi\geq
-\omega$.
\end{proposition}

The next proposition is a consequence of the classical Hartogs' lemma
for psh functions.

\begin{proposition} \label{prop_hartogs}
Let $\varphi_n$ and $\varphi$ be quasi-psh functions on $\P^k$ with
$\ddc\varphi_n\geq -\omega$ and $\ddc\varphi\geq -\omega$. Assume that
$\varphi_n$ converge in $\Lc^1$ to $\varphi$. Let $\widetilde\varphi$ be a
continuous function on a compact subset $K$ of $\P^k$ such that
$\varphi<\widetilde\varphi$ on $K$. Then,
$\varphi_n<\widetilde\varphi$ on $K$ for $n$ large enough. In particular, we
have $\limsup\varphi_n\leq\varphi$ on $\P^k$. 
\end{proposition}

We recall a compactness property of 
quasi-psh functions and also an approximation result
(see also Proposition \ref{prop_regularization_sqp} below).

\begin{proposition}
Let $(\varphi_n)$ be a decreasing sequence of quasi-psh functions with
$\ddc \varphi_n\geq -\omega$. Then, either $\varphi_n$ converge
uniformly to $-\infty$ or $\varphi_n$ converge pointwise and also in
$\Lc^p$, $1\leq p<\infty$, to a quasi-psh function $\varphi$ with
$\ddc\varphi\geq -\omega$. Moreover, for every quasi-psh function
$\varphi$ with $\ddc\varphi\geq -\omega$, there is a
sequence $(\varphi_n)$ of smooth functions such that
$\ddc\varphi_n\geq -\omega$ which decreases to $\varphi$.
\end{proposition} 

Consider now a hypersurface $V$ of $\P^k$ of degree $m$ and the
positive closed $(1,1)$-current $[V]$ of integration on $V$ which is
of mass $m$. Let $\varphi$ be a quasi-potential of $[V]$, i.e. a
quasi-psh function such that $\ddc\varphi=[V] -m\omega$. 
Let $\delta$ be an integer such that the multiplicity of
  $V$ is $\leq \delta$ at every point.
The following
lemma will be useful in the next paragraphs.

\begin{lemma} \label{lemma_dist_V}
 There is a  constant $A>0$  such that
$$\delta\log \dist(\cdot,V)-A\leq\varphi\leq \log \dist(\cdot,V)+A.$$
\end{lemma}
\proof
Let $x=(x_1,\ldots,x_k)=(x',x_k)$ denote the coordinates of $\C^k$. Let
$\Pi:\C^k\rightarrow \C^{k-1}$ with $\Pi(x):=x'$ be the projection
on the first $k-1$ factors.
We can reduce the problem to the local situation where $V$ is a
hypersurface of the unit 
polydisc $\Delta^k$ such that the projection $\Pi:V\rightarrow \Delta^{k-1}$ defines a
ramified covering of degree $s\leq \delta$. For $x'\in\Delta^{k-1}$, denote by $x_{k,1}$,
$\ldots$, $x_{k,s}$ the last coordinates of points in
$\Pi^{-1}(x')\cap V$. Here, these points are repeated according to
their multiplicity.
So, $V$ is the zero set of the Weierstrass polynomial
$$P(x):=(x_k-x_{k,1})\ldots(x_k-x_{k,s}).$$
This is a holomorphic function on $\Delta^k$.
It follows that $\varphi(x)-\log|P(x)|$ is a smooth function. We only
have to prove that
$$ \dist(x,V)^s\lesssim |P(x)| \lesssim \dist(x,V)$$
locally in $\Delta^k$. 
The first inequality follows from the definition of $P$.
Since the derivatives of $P$ are locally
bounded, it is clear that for every $a$ in a compact set of $V$
$$|P(x)|=|P(x)-P(a)|\lesssim |x-a|.$$
Hence, $|P(x)| \lesssim \dist(x,V)$. 
\endproof   

Let $V_t$ denote the $t$-neighbourhood of $V$, i.e. the open set of
points whose distance to $V$ is smaller than $t$. Recall that an
integrable function $\varphi$ on $\P^k$ is said to be {\it dsh} if it is equal
outside a pluripolar set to a difference of two quasi-psh
functions \cite{DinhSibony6}. We identify two dsh functions if they
are equal out of a pluripolar set. 
The space of dsh functions is endowed with the following
norm
$$\|\varphi\|_\DSH:=\|\varphi\|_{\Lc^1}+\inf\|T^+\|$$ 
where $T^\pm$ are positive closed $(1,1)$-currents such that
$\ddc\varphi=T^+-T^-$. The currents $T^+$ and $T^-$ are cohomologous
and have the same mass.
Note that the notion of dsh function can be easily extended to compact
K{\"a}hler manifolds.
We have the following lemma.

\begin{lemma} \label{lemma_dsh_compose}
Let $\chi:\R\cup\{-\infty\}\rightarrow
\R$ be a convex increasing function such that $\chi'$ is bounded. Then,
for every dsh function $\varphi$,
$\chi(\varphi)$ is dsh and
$$\|\chi(\varphi)\|_\DSH\lesssim 1+\|\varphi\|_\DSH.$$ 
\end{lemma}
\proof
Up to a linear change of coordinate on $\R\cup\{-\infty\}$, we can assume that
$\|\varphi\|_\DSH\leq 1$. 
Since $\chi(x)\lesssim 1+|x|$,
$\|\chi(\varphi)\|_{\Lc^1}$ is bounded.
So, it is enough to prove that $\chi(\varphi)$ is dsh and to bound
$\ddc \chi(\varphi)$. We can write $\varphi=\varphi^+-\varphi^-$ out
of a pluripolar set where $\varphi^\pm$ are quasi-psh with bounded
DSH-norm such that $\ddc\varphi^\pm\geq -\omega$. Since $\varphi^\pm$ can
approximated by decreasing sequences of smooth quasi-psh
functions, it is enough to consider the case where $\varphi^\pm$ and $\varphi$ are
smooth. It remains to bound $\ddc \chi(\varphi)$. We have since
$\chi''$ is positive
$$\ddc
\chi(\varphi)=\chi'(\varphi)\ddc\varphi+\chi''(\varphi)d\varphi\wedge\dc\varphi\geq
\chi'(\varphi)\ddc\varphi\geq - \|\chi'\|_\infty T^-,$$
Because $\chi'$ is bounded,
$\ddc\chi(\varphi)$ can be written as a difference of positive closed
currents with bounded mass.  The lemma follows. 
\endproof

\begin{lemma} \label{lemma_cutoff}
For every $t>0$
  there is a smooth function $\chi_t$, $0\leq \chi_t\leq 1$,  with compact support in
  $V_{A_1t^{1/\delta}}$, equal to $1$ on $V_t$ and such that
  $\|\chi_t\|_\DSH\leq A_1$, where $A_1>0$ is a constant independent of $t$.
\end{lemma}
\proof
We only have to consider the case $t\ll 1$.
We will construct $\chi_t$ using Lemma \ref{lemma_dsh_compose} applied twice 
to the function $\varphi$ in
Lemma \ref{lemma_dist_V}. 
Let $\chi:\R\cup\{-\infty\}\rightarrow [0,+\infty[$ be a smooth function
which is convex increasing. We choose $\chi$ such that  $\chi(x)=0$
on $[-\infty,-1]$ and $\chi(x)=x$ for $x\geq 1$. So, we have 
$\max(x,0)\leq \chi\leq \max(x,0)+1$. Let 
$\varphi$ and $A$ be as in Lemma \ref{lemma_dist_V}. Define
$$\phi_t:= - \chi \big(\varphi-\log t-A-1\big) \quad \mbox{and}\quad
\chi_t:= \chi(\phi_t+1).$$
Then, $\phi_t$ and $\chi_t$ are smooth, and by Lemma
\ref{lemma_dsh_compose}, their DSH-norms are bounded uniformly on $t$.
 We deduce from the
properties of $\chi$ that $\chi_t\geq 0$, $\phi_t\leq 0$ and $\phi_t=0$ on
$V_t$. It follows that $\chi_t=1$ on $V_t$. 
Out of $V_{A_1t^{1/\delta}}$ with $A_1\gg 1$, by  Lemma
\ref{lemma_dist_V},  we have $\varphi-\log
t-A-1\gg 0$, hence
$\phi_t=-\varphi+\log t+A+1$. We deduce that
$\phi_t+1\leq -1$ and $\chi_t=0$ there.
This implies the lemma.
\endproof

We recall a notion of {\it capacity} that we introduced in
\cite{DinhSibony6} which can be extended to any compact K\"ahler manifold,
see also \cite{SibonyWong, Alexander}. 
Let 
$$\Pc:=\big\{\varphi\mbox{ quasi-psh}, \quad \ddc\varphi\geq -\omega,\
\max_{\P^k}\varphi=0\big\}.$$ 
For $E\subset \P^k$, define 
$$\capacity(E):=\inf_{\varphi\in\Pc}\exp\big(\sup_E\varphi\big).$$
We have $\capacity(\P^k)=1$, and $E$ is pluripolar if and only if $\capacity(E)=0$. 

Consider a quasi-potential $\varphi$ of a current
$T\in\Cc_1$, i.e. a quasi-psh function such that
$\ddc\varphi=T-\omega$. Quasi-potentials of $T$ differ by constants.
We can associate to each point $a\in \P^k$ the Dirac mass $\delta_a$
at $a$. Define a function $\Uc$ on the extremal elements of
$\Cc_k$ by 
$$\Uc(\delta_a):=\varphi(a).$$
We can extend this function in a unique way to an affine function on
$\Cc_k$ by setting 
$$\Uc(\nu):=\int_{\P^k}\varphi d\nu\qquad \mbox{for } \nu\in\Cc_k.$$
The upper semi-continuity of $\varphi$ implies that $\Uc$ is also u.s.c. on $\Cc_k$. 
We say that $\Uc$ is {\it a super-potential} of $T$. Super-potentials
of a given current differ by  constants. 

Let 
$$\Pc_1:=\big\{\Uc\mbox{ super-potential of a current } T\in\Cc_1, \quad 
\max_{\Cc_k}\Uc=0\big\}.$$ 
For each set $E$ of probability measures in $\Cc_k$, define 
$$\capacity(E):=\inf_{\Uc\in\Pc_1}\exp\big(\sup_{\nu\in E}\Uc(\nu)\big).$$
It is easy to check that for a single measure $\nu$,
 $\capacity(\nu)>0$ if and only if quasi-psh functions are
$\nu$-integrable, i.e. $\nu$ is PB in the sense of \cite{DinhSibony1, DinhSibony6}.
A definition of super-potentials for currents of any bidegree  will be given in the next
paragraph.

\begin{lemma}  Let $E'\subset \P^k$ be a 
Borel set. Let $E$ be the set of measures $\nu\in\Cc_k$
with $\nu(E')=1$. Then, $\capacity(E')=\capacity(E)$.
\end{lemma}
\proof
Since $\Uc$ is affine and u.s.c., the supremum can be taken on the set of
extremal points. It follows that $\max_{\Cc_k}\Uc=0$ if and only if
$\max_{\P^k}\varphi=0$. Moreover, we have
$\sup_E\Uc=\sup_{E'}\varphi$. 
It is now clear that $\capacity(E')=\capacity(E)$.
\endproof


\subsection{Green quasi-potentials of currents}

Let $R$ be a current in $\Cc_p$ with $p\geq 1$. 
If $U$ is a  $(p-1,p-1)$-current such
that $\ddc U=R-\omega^p$, we say that $U$ is a
{\it quasi-potential} of $R$. 
The integral $\langle U,\omega^{k-p+1}\rangle$ is {\it the mean} of $U$.
Such currents $U$ exist but they are not
unique.  When $p=1$ the quasi-potentials
of $R$ differ by constants, when $p>1$ they differ by $\ddc$-closed currents which
can be singular. 
Moreover, for $p>1$,  $U$ is not always defined at every point of
$\P^k$. This is one of the difficulties in the study of positive closed currents
of higher bidegree. We will use constantly the following result which
gives potentials with good estimates.

\begin{theorem} \label{thm_ddbar}
Let $R$ be a current in $\Cc_p$. 
Then, there is a  negative  quasi-potential $U$ of $R$ depending
linearly on $R$ such that for every $r$ with $1\leq r<k/(k-1)$ and for
$1\leq s < 2k/(2k-1)$  
$$\|U\|_{\Lc^r} \leq c_r\quad \mbox{and} \quad \|dU\|_{\Lc^s}\leq c_s$$
for some positive constants $c_r, c_s$  independent of $R$. Moreover, $U$ depends
continuously on $R$ with respect to the $\Lc^r$ topology on $U$ and the
weak topology on $R$
\end{theorem}

We will construct $U$ using a kernel solving the
$\ddc$-equation for the diagonal of $\P^k\times\P^k$. We need a
negative kernel with tame singularities. 
In the case of arbitrary compact K\"ahler manifolds, this is not
always possible \cite{BostGilletSoule}.
In order to simplify the
notation, consider the following general situation. Let $X$ be a
homogeneous compact K{\"a}hler manifold of dimension $n$ and let $G$ be a
complex Lie group of dimension $N$ acting transitively on $X$. 
The following proposition gives some precisions on a result
in Bost-Gillet-Soul{\'e} \cite[Prop. 6.2.3]{BostGilletSoule}, see also Andersson \cite{Andersson}.

\begin{proposition} \label{prop_bgs}
Let $D$ be a submanifold of pure dimension $n-p$ in $X$ with $p\geq 1$
and
 $\Omega$ be a real closed $(p,p)$-form cohomologous to the
current $[D]$.
Then, there
is a negative $(p-1,p-1)$-form $K$ on $X$ smooth outside $D$ such that 
$\ddc K=[D]-\Omega$ which satisfies the following inequalities near $D$
$$\|K(\cdot)\|_\infty\lesssim -\log
\dist(\cdot,D)\dist(\cdot,D)^{2-2p},
\qquad \|\nabla K(\cdot)\|_\infty\lesssim \dist(\cdot,D)^{1-2p}.$$
Moreover, there is a negative dsh function $\eta$ and a positive
closed $(p-1,p-1)$-form $\Theta$ smooth outside $D$ such that
$K\geq \eta\Theta$, $\|\Theta(\cdot)\|_\infty \lesssim \dist(\cdot,D)^{2-2p}$
and $\eta+\log\dist(\cdot,D)$ is bounded near $D$.
\end{proposition}

Note that $\|\nabla K\|_\infty$ is the sum $\sum |\nabla K_i|$ where
$K_i$ are the coefficients of $K$ for a fixed atlas of $X$. 
We first prove the following lemmas.

\begin{lemma} \label{lemma_dsh_singular}
There is a negative dsh function $\eta$ on $X$ smooth outside $D$ such that
$\eta-\log\dist(\cdot,D)$ is bounded. 
\end{lemma}
\proof
Let $\pi:\widetilde X\rightarrow X$ be the blow-up of $X$ along
$D$. 
Denote by $\widehat D:=\pi^{-1}(D)$ the exceptional divisor. If
$\alpha$ is a real closed $(1,1)$-form on $\widetilde X$ cohomologous to
$[\widehat D]$, there is a negative quasi-psh function $\widetilde
\eta$ such that $\ddc\widetilde\eta=[\widehat D]-\alpha$. It is
clear that $\widetilde \eta$ is smooth outside $\widehat D$ and
$\widetilde\eta-\log\dist(\cdot,\widehat D)$ is bounded. Define
$\eta:=\widetilde\eta\circ\pi^{-1}$. Hence, $\eta-\log\dist(\cdot,D)$ is
bounded. Moreover, by a theorem of Blanchard \cite{Blanchard},
$\widetilde X$ is K{\"a}hler. Hence, $\ddc \widetilde\eta$ can be
written as a difference of positive closed currents. It follows that
$\ddc\eta=\pi_*(\ddc\widetilde \eta)$ is also a difference of positive
closed currents. We deduce that $\eta$ is dsh. 
\endproof

\noindent
{\bf Proof of Proposition \ref{prop_bgs}.} 
Let $\Gamma_D \subset G\times D\times X$ denote the graph of the map
$(g,x)\mapsto g(x)$ from $G\times D$ to $X$. 
Let $\Pi_G$ and $\Pi_X$ 
denote the projections
of $\Gamma_D$ onto $G$ and $X$
respectively. 
Observe that $\Pi_G$ defines a trivial fibration.
The map $\Pi_X$ also defines a fibration which is locally trivial. 
Indeed, we can pass from a fiber to another one using the action
$(g,x,g(x))\mapsto (\tau(g), x, \tau(g(x))$ on $G\times D\times X$, of an element $\tau$ of
$G$. So, $\Pi_X$ is a submersion. 
The integrals that
we consider below are computed on some compact subset of $\Gamma_D$.

Let $z$ be a local coordinate on $G$ with $|z|<1$ such that $z=0$ at
$\id$.
Let $\chi$ be a smooth positive function with compact support
in $\{|z|<1\}$ and equal to 1 in a neighbourhood of 0. Define 
$K_G:=\chi\log|z| (\ddc\log |z|)^{N-1}$. This is a negative current
 with support in $\{|z|<1\}$ and 
$\Omega_G:=-\ddc K_G+\delta_0$ is a smooth form.
We have $\|K_G(\cdot)\|_\infty \lesssim -\log|z|\cdot |z|^{2-2N}$  and
$\|\nabla K_G(\cdot)\|_\infty \lesssim |z|^{1-2N}$.

Observe that $\widetilde D:=\Pi_G^{-1}(\id)\cap \Gamma_D$ is
compact and is sent by $\Pi_X$
biholomorphically  to $D$. Therefore,
locally near
$\widetilde D$, one can find coordinates $(x_D,\rho_D,x_G)\in
\C^{n-p}\times\C^p\times \C^{N-p}$ such that
$\widetilde D=\{\rho_D=x_G=0\}$ and $\Pi_X(x_D,\rho_D,x_G)=(x_D,\rho_D)$. 
Define the negative form $K$ by
$$K:=(\Pi_X)_*(\Pi_G^*(K_G)).$$
So, $K$ is smooth outside $D$. 
Using the coordinates $(x_D,\rho_D,x_G)$ and that
$\Pi_G:\Gamma_D\rightarrow G$ is a trivial fibration, we obtain
$$\eta\circ\Pi_X\lesssim \log\dist(\cdot,\widetilde D)\lesssim
-\log|\Pi_G|.$$ 
This,
Lemmas \ref{lemma_dsh_singular} and the above estimates on $K_G$ imply that 
$$K\gtrsim \eta\ (\Pi_X)_*\big(\Pi_G^*(\Theta_G)\big),$$
where $\Theta_G:=\chi (\ddc\log|z|)^{N-1}$.

Define 
$$\Theta:=(\Pi_X)_*\big(\Pi_G^*(\Theta_G)\big).$$ 
Using the local coordinates
$(x_D,\rho_D,x_G)$ and that
$$\|\Pi_G^*(\Theta_G)\|_\infty \lesssim
\dist(\cdot,\widetilde D)^{2-2N}\lesssim (|\rho_D|^2+|x_G|^2)^{1-N}$$ 
on $\Gamma_D$, we obtain
\begin{eqnarray*}
\|\Theta(\cdot)\|_\infty & \lesssim &  \int_{|x_G|\leq 1}
{dx_G\over (|\rho_D|^2+|x_G|^2)^{N-1}}  \leq  \int_{|x_G|\leq 1}
{dx_G\over |\rho_D|^{2N-2}+|x_G|^{2N-2}} \\
& \simeq & \int_0^1 {x^{2N-2p-1}dx  \over |\rho_D|^{2N-2}+x^{2N-2}} 
\lesssim |\rho_D|^{2-2p} \int_0^\infty {ds\over 1+s^{2N-2}}
\lesssim |\rho_D|^{2-2p}.
\end{eqnarray*}
So, we have the  estimate $\|\Theta(\cdot)\|\lesssim
\dist(\cdot,D)^{2-2p}$. 

We then deduce the desired estimate on
$\|K(\cdot)\|_\infty$.  
We also have near $\widetilde D$
$$\|\nabla\Pi_G^*(K_G)(\cdot)\|_\infty \lesssim
\dist(\cdot,\widetilde D)^{1-2N}.$$
A similar computation as above gives that $\|\nabla K(\cdot)\|_\infty\lesssim
\dist(\cdot,D)^{1-2p}$. 
So, the singularities of $K$ satisfy the estimates in the proposition.
We have finally
\begin{eqnarray*}
\ddc K & = &  (\Pi_X)_*(\Pi_G^*(\ddc K_G))=  
(\Pi_X)_*(\Pi_G^*(\delta_\id -\Omega_G))\\
& = & (\Pi_X)_*(\Pi_G^*(\delta_\id))
- (\Pi_X)_*(\Pi_G^*(\Omega_G))\\
& = & [D]-(\Pi_X)_*(\Pi_G^*(\Omega_G))=:[D]-\Omega'.
\end{eqnarray*}
Because $\Omega_G$ is smooth,
$\Omega':=(\Pi_X)_*(\Pi_G^*(\Omega_G))$ is
also smooth. Since $\Omega$ and $\Omega'$ are both cohomologous
to $[D]$, there is a smooth real $(p-1,p-1)$-form $U$ such that
$\ddc U=\Omega-\Omega'$. Adding to $U$ a positive closed form large
enough allows to assume that $U$ is positive. Replacing $K$ by $K-U$
gives a negative form such that $\ddc K=[D]-\Omega$ with the desired
tame singularities.
\hfill $\square$ 

\

\noindent
{\bf Proof of Theorem \ref{thm_ddbar}.} 
We apply Proposition \ref{prop_bgs} to $X:=\P^k\times\P^k$,
$G:=\Aut(\P^k)\times \Aut(\P^k)$ and $D$ the diagonal of $X$. 
Since $\Aut(\P^k)\simeq {\rm PGL}(k+1,\C)$, we can identify
$\Aut(\P^k)$ to a Zariski open set in $\P^{k^2+2k}$ which is the
projectivization of the space of matrices of size $(k+1)\times (k+1)$.  The assumptions in
Proposition \ref{prop_bgs} are easily verified.
Let $(z,\xi)$ denote the homogeneous coordinates of $\P^k\times\P^k$
with $z=[z_0:\cdots:z_k]$ and $\xi:=[\xi_0:\cdots:\xi_k]$. The
diagonal $D$ is given by $\{z=\xi\}$. Choose 
$$\Omega(z,\xi):=\sum_{j=0}^k \omega(z)^j\wedge\omega(\xi)^{k-j}.$$
This form is cohomologous to $[D]$. Using the notation from Proposition
\ref{prop_bgs}, we define
\begin{equation*}
U(z):=\int_{\xi\not=z} R(\xi)\wedge K(z,\xi).
\end{equation*}
Observe that $K$ is smooth out of $D$ and that its coefficients have singularities like
$\log|z-\xi|\cdot |z-\xi|^{2-2k}$ near $D$ (there is an abuse of notation: we
should write $\log|z-\xi|\cdot |z-\xi|^{2-2k}$ on charts $\{z_i=\xi_i=1\}$ which
cover $D$). It follows that
the definition of $U$ makes sense for every current $R$ with measure
coefficients. This is a form with
coefficients in $\Lc^r$. An easy way to see that is to desintegrate
$R$ into currents with support at a point. 
The continuity with
respect to the $\Lc^r$-norm of $U$ and the weak topology on $\Cc_p$,
and the estimate on the
$\Lc^r$-norm of $U$ are easy to check. 

For the rest of the
proposition, by continuity, we can assume that $R$ is a smooth form in $\Cc_p$.
Denote by $\pi_1$ and $\pi_2$ the projections of $\P^k\times\P^k$ on
its factors. Observe that
$$U  =  (\pi_1)_*\big(\pi_2^*(R)\wedge K\big).$$
Hence, $U$ is negative since $K$ is negative and $R$ is positive. Since $R$
is closed, we also have
\begin{eqnarray*}
\ddc U  & = & (\pi_1)_*\big(\pi_2^*(R)\wedge \ddc K\big) \\
& = &  (\pi_1)_*\big(\pi_2^*(R)\wedge [D]\big) -
(\pi_1)_*\big(\pi_2^*(R)\wedge \Omega\big) \\
& = & R-\omega^p.
\end{eqnarray*}
Therefore, $U$ is a quasi-potential of $R$. We also have
$$dU=(\pi_1)_*\big(\pi_2^*(R)\wedge d K\big).$$
Since $dK$ has singularities like
$|z-\xi|^{1-2k}$ near $D$, it is clear that $\|dU\|_{\Lc^s}$
is bounded by a constant independent of $R$. 
\hfill $\square$

\begin{remark} \rm \label{rk_levine}
We call $U$ {\it the Green quasi-potential} of $R$. 
By Theorem \ref{thm_ddbar}, the mean $m$ of 
$U$ is bounded by a constant independent of $R$. So,
$U-m\omega^{p-1}$ is a quasi-potential of mean 0 of $R$. Its mass is
bounded uniformly on $R$. Note that $U$ depends on the choice of $K$.
\end{remark}

We now give some properties of Green quasi-potentials.

\begin{lemma} \label{lemma_levine_c1}
Let  $W'\Subset W$ be open subsets of $\P^k$ and $R$ be  a current in
$\Cc_p$. Assume that the restriction of $R$ to $W$ is a bounded form.
Then, there is a constant
$c>0$ independent of $R$ such that
$$\|U\|_{\Cc^1(W')}\leq c (1+\|R\|_{\infty, W}).$$
\end{lemma}
\proof
Observe that the derivatives of the coefficients of $K$ have
integrable singularities of order $|z-\xi|^{1-2k}$. This and the definition of $U$ imply the result.
\endproof

The precise estimate on the behavior of $U$ in the following
proposition will be needed for the dynamical applications.

\begin{proposition} \label{prop_levine_estimate}
Let $V$, $V_t$ and $\delta$ be as in Lemmas
\ref{lemma_dist_V} and \ref{lemma_cutoff}. 
Let $T_i$, $1\leq i\leq k-p+1$,  be positive closed $(1,1)$-currents on
$\P^k$, smooth on $\P^k\setminus V$. Assume that the quasi-potentials of $T_i$ are
$\alpha_i$-H{\"o}lder continuous with $0<\alpha_i\leq 1$. If $U$ is the Green quasi-potential of a
current $R\in\Cc_p$, then  
$$\Big|\int_{V_t\setminus V} U\wedge T_1\wedge \ldots \wedge T_{k-p+1}\Big|\leq c
t^\beta,\qquad \mbox{with}\quad \beta:=(20k^2\delta)^{-k}
\alpha_1\ldots\alpha_{k-p+1},$$
where $c>0$ is a constant independent of $R$ and of $t$.
\end{proposition}

We will use the notations from Theorem \ref{thm_ddbar} and Proposition \ref{prop_bgs}.
Define  $\eta_M:=\min\big(0,M+\eta\big)$ for $M\geq 0$. 
As in Lemma \ref{lemma_dsh_compose}, we can show that 
$\|\eta_M\|_\DSH$ is bounded independently
of $M$.
We have $\eta_M-M\leq \eta$. Define 
$K_M:=-M\Theta$ and  $K'_M:=\eta_M\Theta$.
Then, $K_M$ is negative closed and we have $K_M+K'_M\lesssim K$. Define also
$$U_M(z):=\int_\xi R(\xi)\wedge K_M(z,\xi)\quad \mbox{and}\quad 
U_M'(z):=\int_\xi R(\xi)\wedge K_M'(z,\xi).$$
The forms $U_M$ is
negative closed of mass $\simeq M$ and we  have 
$U_M+U_M'\lesssim U$. Choose $M:=t^{-\beta}$. We estimate $U_M$ and
$U_M'$ separately. Recall that $U$ is negative and that $\Theta$ has
singularities of order $\dist(z,\xi)^{2-2k}$.

\begin{lemma} We have 
$$\Big|\int_{V_t} U_M\wedge\omega^{k-p+1}\Big| \lesssim t.$$
\end{lemma}
\proof
We can assume $t<1/2$. 
We don't need that $R$ is closed.
So, we can assume that $R$ has
support at a point $a\in\P^k$. We define $U_M$ using the same integral
formula as above.
Then, the coefficients of $U_M$ have
singularities of type $M|x|^{2-2k}$ where $x$ are local
coordinates such that $x=0$ at $a$.
The problem is
local. Since $M\leq t^{-1/2}$, we can assume that $V$ is a hypersurface in a neighbourhood of
the unit ball $B$. It is sufficient to prove that
$$\int_{V_t\cap B}  |x|^{2-2k} (\ddc |x|^2)^k \lesssim
t^{3/2}.$$ 
Let $A$ be a maximal subset of $V\cap B$ such that the distance
between two points in $A$ is $\geq t$. The balls of radius
$2t$ with center in $A$ cover $V\cap B$ and the ones of radius $3t$
cover $V_t\cap B$.
Let $A_n$ be the set of points $p\in A$ such that
$nt\leq |p|< (n+1)t$ and $m_n$ the number of elements of $A_n$. 
Observe that the $m_0+\cdots+m_n$ balls of radius $t/2$ with center in
$A_1\cup\ldots\cup A_n$ are disjoint. They cover an open subset of
$V\cap \{|x|\leq (n+2)t\}$. 
Using Lelong's estimate in Example \ref{ex_Lelong}, see also \cite{Lelong},
 gives that
$$m_0+\cdots+m_n\lesssim n^{2k-2}.$$
Note that $m_0=0$ or 1 and the integral of $|x|^{2-2k}
(\ddc |x|^2)^k$ on a ball of radius $3t$ with center in $A_0$
is bounded by the integral of this function on the ball of
center 0 and of radius $4t$. Hence, it is 
of order $t^2$. For $n\geq 1$, it is clear that 
the integral of the considered form on a ball
with center in  $A_n$ is of order $n^{2-2k} t^2$. Using the estimates on $m_n$ and 
Abel's transform, one obtains
\begin{eqnarray*}
\int_{V_t\cap B} |x|^{2-2k} (\ddc |x|^2)^k
& \lesssim &
 t^2 +\sum_{1\leq n\leq 1/t} m_n n^{2-2k}t^2\\
& \lesssim &  t^2 +\sum_{1\leq n\leq 1/t}
\big[n^{2k-2}-(n-1)^{2k-2}\big] n^{2-2k} t^2\\
& \lesssim &   t^2 + t^2\sum_{1\leq n\leq 1/t} n^{-1}.
\end{eqnarray*}
This implies the lemma.
\endproof

We continue the proof of Proposition \ref{prop_levine_estimate}.
By continuity it is enough to consider the case where $R$ and $U$ are 
smooth. We also have that $U_M$ is smooth.

\begin{lemma} \label{lemma_v_t-estimate}
For every $0\leq l\leq k-p+1$ we have 
$$\Big|\int_{V_t} U_M\wedge T_1\wedge \ldots \wedge T_l\wedge
\omega^{k-p-l+1}\Big|\lesssim t^{\beta_l},\qquad \mbox{where}\quad
\beta_l:=(20 k^2\delta)^{-l} \alpha_1\ldots\alpha_l.$$ 
\end{lemma}
\proof
The proof is by induction. The previous lemma implies the case
$l=0$. Assume the lemma for $l-1$. 
Let $\chi_t$ be as in Lemma \ref{lemma_cutoff}.
We want to prove that
$$\int -\chi_t U_M\wedge T_1\wedge \ldots \wedge T_l\wedge
\omega^{k-p-l+1}\lesssim t^{\beta_l}.$$
Write $T_l=\omega+\ddc u$ with $u$ negative quasi-psh of
class $\Cc^{\alpha_l}$. By induction hypothesis, since
$\chi_t$ has support in $V_{A_1t^{1/\delta}}$, we obtain
$$\int -\chi_t U_M\wedge T_1\wedge \ldots \wedge T_{l-1}\wedge
\omega^{k-p-l+2} \lesssim   t^{\delta^{-1}\beta_{l-1}} \lesssim  t^{\beta_l}.$$
Therefore,  we only have to prove that 
$$\int -\chi_t U_M\wedge T_1\wedge \ldots \wedge T_{l-1}\wedge\ddc u\wedge
\omega^{k-p-l+1}\lesssim t^{\beta_l}.$$
By Proposition \ref{prop_regularization} and Lemma
\ref{lemma_deformation_c1}, 
there is a smooth function $u_\epsilon$
such that $\|u_\epsilon\|_{\Cc^2}\lesssim \epsilon^{-2k^2-4k-2}$ and
$\|u-u_\epsilon\|_\infty\lesssim \epsilon^{\alpha_l}$. Using
Stokes' theorem we can write the left hand side of the previous
inequality as
\begin{eqnarray*}
 \int -\chi_t U_M\wedge T_1\wedge \ldots \wedge T_{l-1}\wedge \ddc
u_\epsilon \wedge \omega^{k-p-l+1} \\
+\int -\ddc \chi_t \wedge U_M\wedge T_1\wedge
\ldots \wedge T_{l-1} (u-u_\epsilon)\wedge
\omega^{k-p-l+1}.
\end{eqnarray*}
By induction hypothesis, the previous estimates on
$\|u_\epsilon\|_{\Cc^2}$ and Lemma
\ref{lemma_cutoff}, we obtain
that the first term is of order at most equal to
$t^{\delta^{-1}\beta_{l-1}}\epsilon^{-2k^2-4k-2}$.
If we write $\ddc \chi_t=T^+-T^-$ with $T^\pm$
positive closed of bounded mass, the second term is of order less than
$$ \epsilon^{\alpha_l}\int T^+ \wedge U_M\wedge T_1\wedge
\ldots \wedge T_{l-1} \wedge
\omega^{k-p-l+1} + \epsilon^{\alpha_l}\int T^- \wedge U_M\wedge T_1\wedge
\ldots \wedge T_{l-1}\wedge
\omega^{k-p-l+1}.$$
These integrals can be computed cohomologically. The currents
$T^\pm$ have bounded mass. Since $K_M=-M\Theta$, we deduce from the
definition of $U_M$ that $-U_M$ is positive and closed of mass
$M=t^{-\beta}$. Therefore, the last sum is 
$\lesssim t^{-\beta}\epsilon^{\alpha_l}$.

Take
$\epsilon:=t^{\delta^{-1}(2k^2+4k+2+\alpha_l)^{-1}\beta_{l-1}}$. We
have
$$1-\frac{2k^2+4k+2}{2k^2+4k+2+\alpha_l}\geq \frac{\alpha_l}{10k^2},$$
then
$$t^{\delta^{-1}\beta_{l-1}}\epsilon^{-2k^2-4k-2}
\lesssim
t^{\delta^{-1}\beta_{l-1}(10k^2)^{-1}\alpha_l}\lesssim t^{\beta_l}$$
and
$$ t^{-\beta}\epsilon^{\alpha_l}\lesssim
t^{-\beta}t^{(10k^2\delta)^{-1}\beta_{l-1}\alpha_l} \lesssim
t^{-\beta} t^{2\beta_l}\lesssim t^{\beta_l}.$$
This implies the desired estimate.
\endproof

\begin{lemma}
We have $\|U'_M\|\lesssim \exp(-M/2).$ 
\end{lemma}
\proof
We can forget that $R$ is smooth and assume that $R$ has support at a
point $a$. The behavior of $\eta$ implies that $U_M'$ has
support in the ball of center $a$ of radius $\lesssim \exp(-M/2)$. The
coefficients of $U_M'$ have singularities $\lesssim -\log|x|\cdot |x|^{2-2k}$ for 
local coordinates $x$ with $x=0$ at $a$. Hence, $\|U'_M\|\lesssim \exp(-M/2).$ 
\endproof

The following lemma completes the proof of Proposition
\ref{prop_levine_estimate}, since $M=t^{-\beta}\gg |\log t|$.

\begin{lemma} For every $0\leq l\leq k-p+1$ we have 
$$\Big|\int U'_M\wedge T_1\wedge \ldots \wedge T_l\wedge
\omega^{k-p-l+1}\Big|\lesssim \exp(-(10k^2)^{-l}\alpha_1\ldots\alpha_lM/2).$$ 
\end{lemma}
\proof
The previous lemma implies the case $l=0$. Assume the lemma for $l-1$
and use the notation from the proof of Lemma \ref{lemma_v_t-estimate}. The integral to bound is equal to
$$\int -U'_M\wedge T_1\wedge \ldots \wedge T_{l-1}\wedge \ddc
u_\epsilon \wedge \omega^{k-p-l+1}+$$
$$+\int_{\P^k\times \P^k} - K'_M\wedge R(\xi)\wedge T_1(z)\wedge
\ldots \wedge T_{l-1}(z) \ddc \big(u(z)-u_\epsilon(z)\big)\wedge
\omega(z)^{k-p-l+1}.$$
Choose $\epsilon= \exp
(-(10k^2)^{-l}\alpha_1\ldots\alpha_{l-1}M)$.
Using the estimate on $\|u_\epsilon\|_{\Cc^2}$, by induction hypothesis, 
the first factor is of order at most equal to
$$\exp(-(10k^2)^{-l+1}\alpha_1\ldots\alpha_{l-1}M/2)\epsilon^{-2k^2-4k-2}\lesssim
\exp(-(10k^2)^{-l}\alpha_1\ldots\alpha_lM/2).$$
The second one is equal to
$$\int_{\P^k\times \P^k} - \ddc K'_M\wedge R(\xi)\wedge T_1(z)\wedge
\ldots \wedge T_{l-1}(z) \big(u(z)-u_\epsilon(z)\big)\wedge
\omega(z)^{k-p-l+1}.$$
Since the DSH-norm of $\eta_M$ in
the definition of $K_M'$ is bounded, 
the first term in the last integral can be bounded by a positive closed current with
bounded mass. So, this integral is
of order at most equal to
$$\|u-u_\epsilon\|_\infty\lesssim \epsilon^{\alpha_l}= \exp
(-(10k^2)^{-l}\alpha_1\ldots\alpha_{l-1}\alpha_lM).$$ 
This implies the result.
\endproof

We will use the following lemma in the study of deformation of currents.

\begin{lemma} \label{lemma_deform_pot}
Let $R$ be a current in $\Cc_p$ and $U$ be a quasi-potential of mean
$m$ of $R$. Let $R_{\theta y}=(\tau_{\theta y})_*(R)$ be defined as in Paragraph
\ref{section_topology}.
Then, there is a quasi-potential $U_{\theta y}'$ of $R_{\theta y}$ of
  mean $m$ such that $U_{\theta y}'-(\tau_{\theta y})_* (U)$ is a smooth form with 
$$\|U_{\theta y}'-(\tau_{\theta y})_* (U)\|_{\Cc^2}\leq c\|U\||\theta|$$
where $c>0$ is a constant independent of $R$, $U$, $\theta$ and $y$.
\end{lemma}
\proof
Since $\|(\tau_{\theta y})_*(\omega^p)-\omega^p\|_{\Cc^2}\lesssim
|\theta|$, there is a $(p-1,p-1)$-form $\Omega_{\theta
  y}$ such that $\|\Omega_{\theta y}\|_{\Cc^2}\lesssim |\theta|$ and 
$\ddc \Omega_{\theta y} = (\tau_{\theta
  y})_*(\omega^p)-\omega^p$. It is clear that the mean $m''$ of
$\Omega_{\theta y}$ is of order $\lesssim |\theta|$. Set 
 $U'_{\theta y} := (\tau_{\theta y})_*(U) +\Omega_{\theta y}$.
So, the mean $m'$ of $U'_{\theta y}$ satisfies
\begin{eqnarray*}
|m'-m| & = & \Big| \int (\tau_{\theta y})_*(U)\wedge
\omega^{k-p+1} +m'' -\int U\wedge \omega^{k-p+1}\Big| \\
& \leq & |m''| +\Big|\int U\wedge \big[(\tau_{\theta y})^*(\omega^{k-p+1})-\omega^{k-p+1}\big]\Big|.
\end{eqnarray*}
The last term is of order $\lesssim \|U\||\theta|$ since $\|(\tau_{\theta
  y})^*(\omega^{k-p+1})-\omega^{k-p+1}\|_\infty$ is of order $\lesssim
|\theta|$. 
Subtracting from
$U'_{\theta y}$ the form $(m'-m)\omega^{p-1}$, which is of order
$\lesssim |\theta|$,  gives a
quasi-potential satisfying the lemma.
\endproof


\subsection{Structural varieties in the spaces of currents}

\label{section_structural}

The notion of structural varieties of $\Cc_p$ was introduced
in \cite{DinhSibony7}, see also \cite{Dinh2}. In some sense, we consider
$\Cc_p$ as a space of infinite dimension admitting "complex
subvarieties" of finite dimension. 
The emphasis is that in order to connect two closed currents we use a
{\bf closed} current in higher dimension. 
Holomorphic families of analytic
cycles of codimension $p$ are examples of structural varieties in
$\Cc_p$. Other examples of structural varieties can be obtained by
deforming a given current in $\Cc_p$ using a holomorphic family of
automorphisms. The reader will find in Dujardin \cite{Dujardin} and in
\cite{DinhDujardinSibony} an application of
such a deformation to the dynamics of H\'enon-like maps. General
structural varieties are more flexible, and this is crucial in our
study.

Let $X$ be a complex manifold, $\pi_X: X\times
\P^k\rightarrow X$ and $\pi:X\times \P^k\rightarrow \P^k$ denote the canonical projections.
Consider a positive  closed $(p,p)$-current $\Rc$ in $X\times \P^k$.
By slicing theory \cite{Federer}, 
the slices $\langle
\Rc,\pi_X,x\rangle$ exist for almost every $x\in X$. Such a slice is a
positive  closed $(p,p)$-current on 
$\{x\}\times \P^k$ (following \cite{DinhSibony7}, 
we can prove that the slices exist for $x$ out of a
pluripolar set). We often identify $\langle \Rc,\pi_X,x\rangle$ with a $(p,p)$-current $R_x$ in $\P^k$.  

\begin{lemma}
The mass of $R_x$ does not depend on $x$.
\end{lemma}
\proof
Define $\Rc':=\Rc\wedge \pi^*(\omega^{k-p})$. Then, $\Rc'$ is  positive
closed on $X\times \P^k$ and $(\pi_X)_*(\Rc')$ is closed of
bidegree $(0,0)$ on $X$. Hence, it is a constant function. So, the
function
$$\varphi(x):=\|\langle \Rc',\pi_X,x\rangle\| =\int_{\P^k} R_x\wedge\omega^{k-p}=\|R_x\|$$
is constant. The lemma follows.
\endproof

We assume that the mass of $R_x$ is equal to 1.
The map
$x\mapsto R_x$ is defined almost everywhere on $X$ with
values in $\Cc_p$.

\begin{definition} \rm
We say that the
map $x\mapsto R_x$ or the family  $(R_x)_{x\in X}$ defines a
{\it structural variety} in $\Cc_p$. The positive
closed $(1,1)$-current 
$$\alpha_\Rc:=(\pi_X)_*(\Rc\wedge \pi^*(\omega^{k-p+1}))$$ 
on $X$ is called the {\it curvature} of the structural
variety, see Propositions \ref{lemma_sqp_qpsh} and \ref{prop_sqp_qpsh}
below. 
\end{definition}

\begin{definition} \rm
A structural variety associated to $\Rc$ is said to
be {\it special} if $R_x$ exists for every $x\in X$,  
$R_x$ depends continuously on $x$ and if the
curvature is a smooth form.
\end{definition}

In  order to simplify the argument, we restrict
to special structural
varieties or discs. The most useful structural discs in this work are
$(R_\theta)_{\theta\in\Delta}$, see Introduction and Lemma \ref{lemma_common_curvature} below.


\subsection{Deformation by automorphisms} \label{section_deformation}
 
Using the automorphisms of $\P^k$, we will construct some
special structural discs in $\Cc_p$ that we will use
later on. 
We first construct large
structural varieties parametrized by $X=\Aut(\P^k)$.

\begin{proposition} \label{prop_big_structural}
Let $R$ be a current in $\Cc_p$. 
Then, the map $h:\Aut(\P^k)\rightarrow \Cc_p$ 
with $h(\tau)=R_\tau:=\tau_*(R)$ defines a
special structural variety in $\Cc_p$. Moreover, its
curvature is bounded by a smooth positive $(1,1)$-form independent of $R$.
\end{proposition}
\proof
For any smooth test form $\Phi$, we have
$\langle R_\tau,\Phi\rangle = \langle R, \tau^*(\Phi)\rangle$.
So, clearly $\tau\mapsto R_\tau$ is continuous.
Consider the holomorphic map $H:\Aut(\P^k)\times \P^k\rightarrow\P^k$
defined by $H(\tau,z):=\tau^{-1}(z)$. The current $\Rc:=H^*(R)$
is  positive closed of
bidegree $(p,p)$. It is easy to check from the definition of slices that $R_\tau=\langle
\Rc,\pi_X,\tau\rangle$. Hence, $h$ defines a continuous structural variety.

Now, we have to show that the curvature 
$$\alpha_\Rc:=(\pi_X)_* \big(H^*(R)\wedge\pi^*(\omega^{k-p+1})\big)$$
is a smooth form. We prove this for any current $R$ of mass $\leq 1$ not
necessarily closed. Then, we can assume that $R$ is supported at a
point $a$, that is, there is a tangent $(k-p,k-p)$-vector $\Psi$ at
$a$ of
norm $\leq 1$ such that
$R=\delta_a\wedge \Psi$ (the general case is obtained using a
desintegration $R$ into currents of the previous type).
We have
$H^*(R)=[H^{-1}(a)]\wedge \widetilde\Psi$, where $\widetilde \Psi$ is
a $(k-p,k-p)$-vector field with support in $H^{-1}(a)$ such that
$H_*(\widetilde \Psi)=\Psi$. Because $H$ is a submersion, we can choose
$\widetilde\Psi$ smooth on $H^{-1}(a)$. 

Since $H^{-1}(a)$ is a
holomorphic graph
over $\Aut(\P^k)$, the form $\alpha_\Rc$ defined above is the direct
image of $[H^{-1}(a)]\wedge \widetilde\Psi\wedge \pi^*(\omega^{k-p+1})$
by $\pi_X$. So, $\alpha_\Rc$ is
smooth. Moreover, the $\Cc^s$-norm of $\alpha_\Rc$ on any fixed compact subset of $\Aut(\P^k)$
is uniformly bounded for every $s\geq 0$. The proposition follows. 
\endproof

\begin{remark} \rm
If $i:\Delta\rightarrow \Aut(\P^k)$ is a holomorphic map, then 
$x\mapsto i(x)_*R$, which is equal to $h\circ i$, defines a
special structural disc. We can also construct a structural disc
passing through $R$ and through the current of integration on a fixed plane of
codimension $p$ \cite{Dinh2}. So, $\Cc_p$ is connected by structural discs.  
\end{remark}

Let $R$ be a current in $\Cc_p$. The following lemma gives us 
a useful special structural disc passing through $R$.

\begin{lemma} \label{lemma_common_curvature}
Let $R_\theta$ be the currents constructed in Paragraph \ref{section_topology}. Then, the
family $(R_\theta)$ defines a special structural
  disc whose curvature is bounded by a smooth positive 
  $(1,1)$-form $\alpha$ which does not depend on $R$. 
\end{lemma}
\proof
By Proposition \ref{prop_big_structural}, 
for $|y|<1$, the family $(R_{\theta y})_{\theta\in\Delta}$ defines a
special disc in $\Cc_p$. Moreover, the $\Cc^s$-norm of its curvature  is bounded
uniformly with respect to $R$ and $y$. In particular, this curvature
is bounded by a positive form $\alpha$ which does not depend on
$R$ and $y$. 

Let $\Rc_y$ denotes the $(p,p)$-current on $\Delta\times \P^k$
associated to the structural disc $(R_{\theta y})$ and define
$\Rc:=\int \Rc_y d\rho(y)$. Recall that $R_\theta =
\int_y R_{\theta y} d\rho(y)$. Hence, $(R_\theta)$ is the family of 
slices of $\Rc$ and it defines a structural disc in $\Cc_p$.
We know that $R_\theta$ depends continuously on $\theta$. This and the above
properties of $(R_{\theta y})$ imply that the curvature of
$(R_\theta)$ is  bounded by $\alpha$.
\endproof


\section{Super-potentials of currents} \label{section_sp}

Consider a current $S$ in $\Cc_p$. We 
introduce a {\it super-potential} associated to $S$. It is an  affine upper
semi-continuous (u.s.c. for short) function $\Uc_S$ defined on $\Cc_{k-p+1}$
with values in $\R\cup\{-\infty\}$.


\subsection{Super-potentials of currents} \label{section_super_def}

Assume first that $S$ is a smooth form in $\Cc_p$. The general
case will be obtained using a regularization of $S$.
Consider an element $R$ of $\Cc_{k-p+1}$
and fix a real number $m$.
Define
\begin{equation} \label{eq_sqp}
\Uc_S(R):=\langle S, U_R \rangle, \quad U_R \mbox{ a quasi-potential
of mean } m \mbox{ of } R.
\end{equation}

\begin{lemma} \label{lemma_super_def}
The integral $\langle S, U_R \rangle$ does not depend on
   the choice of $U_R$ with a fixed mean
$m$. It defines an affine continuous function $\Uc_S$ 
on $\Cc_{k-p+1}$.
Moreover, if $U_S$ is a smooth
quasi-potential of $S$ with mean $m$, then $\Uc_S(R)=\langle
U_S,R\rangle$. In particular, we have $\Uc_S(\omega^{k-p+1})=m$. 
\end{lemma}
\proof
Let $U_S$ be a smooth quasi-potential of $S$ with mean $m$.
Using Stokes' formula, we obtain 
\begin{eqnarray*}
\Uc_S(R) & = & \langle S, U_R\rangle = \langle S-\omega^p, U_R\rangle
+\langle \omega^p, U_R\rangle \\
& = & \langle \ddc U_S,U_R\rangle +m = \langle U_S,
\ddc U_R\rangle +m\\
& = &  \langle U_S,
R-\omega^{k-p+1}\rangle +m =\langle U_S,R\rangle.
\end{eqnarray*}
This also shows that $\Uc_S(R)$ is independent of the choice of $U_R$
and it depends continuously on $R$. It is clear that $\Uc_S$ is
affine. 
\endproof
 
We say that $\Uc_S$ is the {\it super-potential of mean $m$} of
$S$. One obtains the super-potentials of mean $m'$ by adding $m'-m$ to
the super-potential of mean $m$.
We will see latter that the following lemma
holds also for an arbitrary current $S$ in $\Cc_p$  smooth or not,
see Corollary \ref{cor_def_sp} below.

\begin{lemma} \label{lemma_max_sqp}
There is a constant $c\geq 0$ independent of $S$ such
that if $\Uc_S$ is the super-potential of mean $m$ of $S$, then 
$\Uc_S\leq m+c$ everywhere.  
\end{lemma}
\proof
Without loss of generality we can assume $m=0$. 
Let $U_R'$ be the Green quasi-potential of $R$ which is a negative
 current
and $m'$ the mean of $U_R'$.  Then,
$U_R:=U_R'-m'\omega^{k-p}$ is a quasi-potential of mean $0$ of $R$.
By Lemma \ref{lemma_super_def}, since $U'_R$ is negative and $S$ is positive, we have
$$\Uc_S(R)=\langle S,U_R\rangle=\langle S,U'_R\rangle -m'\leq
-m'.$$
We have seen in Remark \ref{rk_levine} that $|m'|$ is bounded by a constant
independent of $S$. This implies the result. 
\endproof

As we have seen in the last paragraph, the convex set $\Cc_{k-p+1}$ 
can be considered as an infinite dimensional space admitting ``complex
subvarieties'' of finite
dimension. With this point of view, we can consider $\Uc_S$ as a
quasi-psh function on $\Cc_{k-p+1}$. 
More precisely, we will show that the restriction of $\Uc_S$ to a
special structural variety is a quasi-psh function, see Proposition
\ref{prop_sqp_qpsh} below.

We now extend the definition of $\Uc_S$ to an arbitrary current $S$ in $\Cc_p$.
For $R$ smooth, define $\Uc_S(R)$ as in (\ref{eq_sqp}) with $U_R$
smooth. Observe that $\Uc_S(R)$ depends continuously on $S$. 
We can show as in Lemma \ref{lemma_super_def} that the definition is independent
of the choice of $U_R$. We will extend $\Uc_S$
to a function on $\Cc_{k-p+1}$ with
values in $\R\cup\{-\infty\}$.
The reader can check that for $p=1$ we will obtain the same
super-potentials introduced in Paragraph \ref{section_quasi_psh}.

Let $(R_\theta)$ be the special structural disc in
$\Cc_{k-p+1}$ constructed  
in Paragraphs \ref{section_topology} and \ref{section_deformation} and
let $\alpha$ be as in
Lemma \ref{lemma_common_curvature}. Recall that  $R_\theta$ is smooth for
$\theta\not=0$.

\begin{lemma} \label{lemma_sqp_qpsh}
The function $u(\theta):=\Uc_S(R_\theta)$ defined on $\Delta^*$
can be extended as a quasi-subharmonic function on $\Delta$ such that
$\ddc u\geq -\alpha$. 
\end{lemma}
\proof
Proposition \ref{prop_regularization} implies that $u$ is continuous on $\Delta^*$. Lemma
\ref{lemma_max_sqp} holds for $S$ singular and $R$ smooth. So, $u$ is
bounded from above. Let $\Rc$ be the $(k-p+1,k-p+1)$-current in
$\Delta\times \P^k$ associated to $(R_\theta)$ and let $\pi_\Delta$, $\pi$
be as in Paragraph \ref{section_deformation}. Observe that $\Rc$ is smooth
on $\Delta^*\times \P^k$. 
If $U_S$ is a quasi-potential of mean $m$ of $S$, then by definition of
$\Uc_S$ we have
$$u=(\pi_\Delta)_*\big(\Rc\wedge \pi^*(U_S)\big)$$
in the sense of currents on $\Delta^*$. It follows that 
$$\ddc u = (\pi_\Delta)_*\big(\Rc\wedge \pi^*(\ddc U_S)\big)\geq
-(\pi_\Delta)_*\big(\Rc\wedge \pi^*(\omega^p)\big)
\geq -\alpha.$$
If $v$ is a smooth function such that $\ddc v=\alpha$, then $u+v$ is
subharmonic on $\Delta^*$. Since $u$ is bounded from above, $u+v$ can
be extended to a subharmonic function. The lemma follows. Observe that
if $R$ is a smooth form, then $u(\theta)$ is defined and is a
continuous function on $\Delta$. It is quasi-subharmonic and satisfies
$\ddc u\geq -\alpha$. 
\endproof

Recall that $S_\theta$ is defined 
as in Paragraphs \ref{section_topology} and \ref{section_deformation}
for $S$ instead of $R$. 
By Lemma \ref{lemma_deform} and Proposition \ref{prop_regularization}, 
$S_\theta$ is smooth and converges to $S$ when $\theta$ tends to 0.

\begin{proposition} \label{prop_sp_s_theta}
Let $\Uc_{S_\theta}$ denote the super-potential of mean $m$ of
$S_\theta$. Then, $\Uc_{S_\theta}(R)$ converge to
$u(0)$ when $\theta\rightarrow 0$.
In particular, if $R$ is a
smooth form, then $\Uc_{S_\theta}(R)$ converge to $\Uc_S(R)$.
\end{proposition}
\proof
When $R$ is smooth, we have $u(0)=\Uc_S(R)$.
So, we deduce easily the last assertion
from the first one. By Lemma \ref{lemma_sqp_qpsh}, there is a
constant $A>0$ independent of $R$ and $S$ such that $u(\theta)+A|\theta|^2$ is
subharmonic. Since this function is radial (recall here that $\rho$ is
radial, see Introduction), it decreases to $u(0)$
when $|\theta|$ decreases to 0.
Therefore, the proposition is deduced from Lemma \ref{lemma_super_error} below. 
\endproof
 
\begin{lemma} \label{lemma_super_error}
There is a constant $c>0$ independent of $R$ and $S$ such that
$$|\Uc_{S_\theta}(R)-\Uc_S(R_\theta)|=|\Uc_{S_\theta}(R)-u(\theta)|\leq
c|\theta|$$
for $\theta\in \Delta^*$.
\end{lemma}
\proof
Since $R$ can be approximated by smooth forms in $\Cc_{k-p+1}$, we can assume $R$
smooth. Then, we can also assume $S$ smooth. Indeed, the following estimates are
uniform on $R$ and $S$.
Let $U_S$ be a smooth quasi-potential of mean $m$ of $S$ with bounded mass. Define
$U_{\theta y}:= (\tau_{\theta y})^*U_S$. We have
$$\Uc_S(R_\theta)=\int_y \langle U_S,(\tau_{\theta y})_*R\rangle
d\rho(y) =\int_y \langle U_{\theta y},R\rangle d\rho(y).$$
As in Lemma \ref{lemma_deform_pot}, we show that there is a
quasi-potential $U'_{\theta y}$ of mean $m$ of $(\tau_{\theta y})^*(S)$
such that $\|U'_{\theta y} - U_{\theta y}\|_{\Cc^2}\lesssim |\theta|$.
We have
$$\Uc_{S_\theta}(R)=\int_y \langle U'_{\theta y},R\rangle d\rho(y).$$
The estimate on $U_{\theta y}'-U_{\theta y}$ implies that
$$|\Uc_{S_\theta}(R)-\Uc_S(R_\theta)|= \Big|\int_y
\langle U'_{\theta y}-U_{\theta y}, R\rangle d\rho(y) \Big|\lesssim
|\theta|.$$
The proof is complete.
\endproof

\begin{proposition} \label{prop_regularization_sqp}
There is a sequence
of smooth forms $(S_n)$ in $\Cc_p$ 
with super-potentials $\Uc_n$ of mean $m_n$ 
such that
\begin{enumerate}
\item[$\bullet$] $\supp(S_n)$ converge to $\supp(S)$;
\item[$\bullet$] $S_n$ converge to $S$ and $m_n\rightarrow m$;
\item[$\bullet$] $(\Uc_n)$ is a decreasing sequence;
\end{enumerate}
Moreover, if $S_n$, $m_n$ and $\Uc_n$ satisfy the last two properties,
then $\Uc_n(R)$ converge to $u(0)$. In particular, if $R$ is a
smooth form in $\Cc_{k-p+1}$, then $\Uc_n(R)$ converge to $\Uc_S(R)$.
\end{proposition}
\proof
Consider $S_n:=S_{\theta_n}$ where $(\theta_n)$ is a sequence in
$\Delta^*$ such that $|\theta_n|$ decrease to $0$ and that $\sum|\theta_n|$ is finite.
Define 
$$m_n:=m+A|\theta_n|^2+2c\sum_{i=n}^\infty|\theta_i|$$
where $c$ and $A$ are the constants introduced in Lemma \ref{lemma_super_error} and in the proof of
Proposition \ref{prop_sp_s_theta}. It is clear that
$S_n\rightarrow S$, $\supp(S_n)\rightarrow\supp(S)$ 
and $m_n\rightarrow m$. Define $\Uc_n:=\Uc_{S_n}+m_n-m$. This
is the
super-potential of mean $m_n$ of $S_n$. Lemma \ref{lemma_super_error} implies that
\begin{eqnarray*}
\Uc_n(R)-\Uc_{n+1}(R)  &\geq&
\Uc_{S_n}(R)-\Uc_{S_{n+1}}(R) +
A(|\theta_n|^2-|\theta_{n+1}|^2) +2c |\theta_n|\\ 
& \geq & 
\big[u(\theta_n)+A|\theta_n|^2\big] -  \big[u(\theta_{n+1})+A|\theta_{n+1}|^2\big].
\end{eqnarray*}
We have seen that $u(\theta)+A|\theta|^2$ is radial
subharmonic and decreases to $u(0)$ when $|\theta|$ decreases to
0. Hence, $(\Uc_n)$ is decreasing. This implies the first assertion
of the proposition.

For the second assertion, we show that $u_n(0)$ converge to
$u(0)$. Observe that
by definition, $\Uc_n$ converge to $\Uc_S$ on smooth forms $R$ in $\Cc_{k-p+1}$. Define
$u_n(\theta):=\Uc_n(R_\theta)$. Hence, $u_n$ converge to $u$ pointwise on $\Delta^*$. 
On the other hand, Lemma \ref{lemma_sqp_qpsh} implies that $(u_n+A|\theta|^2)$ is a
decreasing sequence of subharmonic functions for $A$ large enough. Hence, it converges
pointwise to a 
subharmonic function. We deduce that $u_n(0)$ converge to $u(0)$. 
This completes the proof.
\endproof

\begin{corollary} \label{cor_def_sp}
$\Uc_S$ can be extended in a unique way to an affine u.s.c function on
$\Cc_{k-p+1}$ with values in $\R\cup\{-\infty\}$, also denoted by $\Uc_S$, such that
$$\Uc_S(R)=\lim_{\theta\rightarrow 0}\Uc_{S_\theta}(R)=\lim_{\theta\rightarrow 0}\Uc_{S}(R_\theta).$$
In particular, we have 
$$\Uc_S(R)=\limsup_{R' \rightarrow R} \Uc_S(R') \mbox{ with } R' \mbox{ smooth}.$$
Moreover, if $c$ is the constant in Lemma \ref{lemma_max_sqp}, then
$\Uc_S\leq m+c$, independently of $S$.
\end{corollary}
\proof
Proposition \ref{prop_regularization_sqp} implies that the decreasing 
limit of $\Uc_{S_n}$ is an
extension of $\Uc_S$. Denote also by $\Uc_S$ this extension.
Since $\Uc_{S_n}$ are affine and continuous, 
$\Uc_S$ is affine and u.s.c. with values in $\R\cup\{-\infty\}$. 
In particular, we have
$$\Uc_S(R)\geq \limsup_{R' \rightarrow R}
\Uc_S(R') \mbox{ with } R' \mbox{ smooth}.$$
Proposition \ref{prop_regularization_sqp} implies also that
$\Uc_S(R)=u(0)$. By Proposition \ref{prop_sp_s_theta} and Lemma \ref{lemma_super_error}, 
we have
$$\Uc_S(R)=u(0)=\lim_{\theta\rightarrow 0} u(\theta)=
\lim_{\theta\rightarrow 0} \Uc_S(R_\theta)=
\lim_{\theta\rightarrow 0} \Uc_{S_\theta}(R).$$
The second limit is bounded above by 
$$\limsup_{R' \rightarrow R}
\Uc_S(R') \mbox{ with } R' \mbox{ smooth}.$$
It follows that 
$$\Uc_S(R)=\limsup_{R' \rightarrow R}
\Uc_S(R') \mbox{ with } R' \mbox{ smooth}.$$
The uniqueness of the extension of $\Uc_S$ is clear.
The inequality $\Uc_S\leq m+c$ is a consequence of  Lemma
\ref{lemma_max_sqp}.
\endproof

\begin{definition}\rm
We call $\Uc_S$ {\it the super-potential of mean $m$} of $S$. 
\end{definition}

It is clear that if $\Uc_S$ is the super-potential of mean
$m$ of $S$, then the super-potential of mean $m'$ of $S$ is equal to
$\Uc_S+m'-m$. The following result applied to $I=\varnothing$, shows that the super-potentials
determine the currents.

\begin{proposition}  \label{prop_unique_sp}
Let $I$ be a compact subset in $\P^k$ with $(2k-2p)$-dimensional
Hausdorff measure $0$.
Let $S$, $S'$ be currents in $\Cc_p$ and $\Uc_S$, $\Uc_{S'}$ be
super-potentials of $S$, $S'$. If
  $\Uc_S=\Uc_{S'}$ on smooth forms in $\Cc_{k-p+1}$ with compact support in  $\P^k\setminus
  I$, then $S=S'$. 
\end{proposition}
\proof
If $R$ is a current in $\Cc_{k-p+1}$ with compact support in
$\P^k\setminus I$, then $R_\theta$ has compact support
in $\P^k\setminus I$ for $\theta$ small enough. On the other hand, since
$R_\theta$ is smooth, we have
$$\Uc_S(R)=\lim_{\theta\rightarrow 0}
\Uc_S(R_\theta)=\lim_{\theta\rightarrow 0} \Uc_{S'}(R_\theta)=\Uc_{S'}(R).$$
Hence, $\Uc_S=\Uc_{S'}$ on every current $R$ with compact support in
$\P^k\setminus I$.
The hypothesis on the Hausdorff measure of $I$ implies that a generic projective subspace $P$ of dimension
$p-1$ does not intersect $I$. We can write $\omega^{k-p+1}$ as an
average of currents $[P]$. Since $\Uc_S=\Uc_{S'}$ at $[P]$ and
since $\Uc_S$ and $\Uc_{S'}$ are affine, they are equal at
$\omega^{k-p+1}$. Hence, $\Uc_S$ and $\Uc_{S'}$ have the same mean. We can
assume that this mean is 0.

If $K$ is compact in $\P^k\setminus I$, using an average of $[P]$,
we can construct a smooth form $R_1$ in $\Cc_{k-p+1}$ with  
 compact support in $\P^k\setminus I$ which
is strictly positive on $K$.
We show that $S=S'$ on $K$. Let $\Phi$ be a smooth $(k-p,k-p)$-form
with compact support on $K$. If $c>0$ is a constant large enough,
$cR_1+\ddc\Phi$ is a positive closed form of mass $c$ since it is
cohomologous to $cR_1$. We can write $cR_1+\ddc\Phi=cR_2$ with
$R_2\in\Cc_{k-p+1}$. We have $\Uc_S(R_1)=\Uc_{S'}(R_1)$ and  $\Uc_S(R_2)=\Uc_{S'}(R_2)$.
If $U_S$ is a quasi-potential of mean 0 of $S$, we have
\begin{eqnarray*}
\langle S,\Phi\rangle & = & \langle S-\omega^p,\Phi\rangle +\langle
\omega^p,\Phi\rangle
 =  \langle \ddc U_S,\Phi\rangle +\langle \omega^p,\Phi\rangle\\
& = & \langle U_S,\ddc \Phi\rangle +\langle \omega^p,\Phi\rangle
 =  \langle U_S,cR_2-cR_1\rangle +\langle \omega^p,\Phi\rangle\\
& = & c\Uc_S(R_2) -c\Uc_S(R_1)+\langle \omega^p,\Phi\rangle.
\end{eqnarray*}
The current $S'$ satisfies the same identity. We deduce that $\langle
S,\Phi\rangle=\langle S',\Phi\rangle$.
Hence, $S=S'$ on $K$. It follows that $S=S'$ on $\P^k\setminus I$. 
The hypothesis on the Hausdorff measure of $I$ implies that 
$S$ and $S'$ have no mass on $I$ \cite{HarveyPolking}. Therefore, $S=S'$ on $\P^k$.
\endproof


\subsection{Properties of super-potentials}

The following proposition extends Lemma \ref{lemma_sqp_qpsh}. It
shows that in some sense super-potentials can be considered
as quasi-psh functions on $\Cc_{k-p+1}$. In particular, they inherit
the compactness property of $\Cc_p$.

\begin{proposition} \label{prop_sqp_qpsh}
Let $(R_x)_{x\in X}$ be an
arbitrary special structural variety in $\Cc_{k-p+1}$
and $\alpha$ be the associated curvature.
Then, either $\Uc_S(R_x)=-\infty$ for every $x\in X$ or 
$x\mapsto \Uc_S(R_x)$ is a quasi-psh function on
$X$ such that $\ddc \Uc_S(R_x)\geq -\alpha$. 
\end{proposition}
\proof
By Proposition  \ref{prop_regularization_sqp}, it is enough to
consider the case where $S$ is smooth. The proof is the same as in
Lemma \ref{lemma_sqp_qpsh}.
Let $\Rc$, $\pi_X$ and $\pi$ be as in Paragraph
\ref{section_structural}. Then, $x\mapsto \Uc_S(R_x)$ is continuous and we have
$$\Uc_S(R_x)=(\pi_X)_* (\Rc\wedge \pi^*(U_S))$$
which implies that
$$\ddc \Uc_S(R_x)=(\pi_X)_*\big(\Rc\wedge \pi^*(\ddc
U_S)\big)\geq -(\pi_X)_*\big(\Rc\wedge \pi^*(\omega^{p})\big)
=-\alpha.$$
This completes the proof.
\endproof

The following result is the analogue of the classical Hartogs' lemma
for psh functions, see also Proposition \ref{prop_hartogs}.

\begin{proposition} \label{prop_hartogs_sqp}
Let $(S_n)$ be a sequence in $\Cc_p$  converging to a current $S$. Let $\Uc_{S_n}$
(resp. $\Uc_S$) be the super-potential of mean $m_n$ (resp. $m$) of $S_n$
(resp. $S$).
Assume that $m_n$ converge to $m$. Let $\Uc$ be a continuous function
on a compact subset $K$ of $\Cc_{k-p+1}$ such that $\Uc_S<\Uc$ on $K$. Then, for $n$ large enough we have 
$\Uc_{S_n}< \Uc$ on $K$. In particular, we have $\limsup \Uc_{S_n}\leq
\Uc_S$ on $\Cc_{k-p+1}$. 
\end{proposition}
\proof
Recall that $\Uc_S$ is u.s.c., $\Uc$ is continuous and $\Cc_{k-p+1}$ is
compact.  The proposition can be applied to $K=\Cc_{k-p+1}$.
Assume
there are currents $R_n$ in $K$ such that
$\Uc_{S_n}(R_n)\geq \Uc(R_n)$. Extracting a subsequence allows to
assume that $R_n$ converge to a current $R$ in $K$. 
Let $(R_{n,\theta})_{\theta\in\Delta}$ be the special structural disc associated to
$R_n$ constructed as in Paragraphs \ref{section_topology} and \ref{section_deformation}. Define 
$u_n(\theta):=\Uc_{S_n}(R_{n,\theta})$. Proposition \ref{prop_sqp_qpsh} implies
that $u_n$ is quasi-subharmonic and $\ddc u_n\geq -\alpha$ with
$\alpha$ as in Lemma \ref{lemma_common_curvature}. The first
assertion of Proposition \ref{prop_regularization}
implies that $u_n$ converge pointwise to $u(\theta):=\Uc_S(R_\theta)$ on $\Delta^*$. It
follows from the Hartogs' lemma for subharmonic functions that 
$$\Uc_S(R)=u(0)\geq\limsup_{n\rightarrow\infty}
u_n(0)=\limsup_{n\rightarrow\infty} 
\Uc_{S_n}(R_n)\geq \Uc(R).$$
This is a contradiction. The proof of the first assertion is
complete. Taking $K=\{R\}$ and $\Uc(R)=\Uc_S(R)+\epsilon$ gives the
second assertion.
\endproof

\begin{definition} \rm \label{def_cv_h}
Let $S_n$, $S$, $\Uc_{S_n}$, $\Uc_S$, $m_n$ and $m$ be as in
Proposition \ref{prop_hartogs_sqp}. If $\Uc_{S_n}\geq \Uc_S$ for every
$n$, we say that {\it $S_n$ converge to $S$ in the Hartogs'
  sense} or {\it $S_n$ H-converge to $S$} for short.  If a current $S'$ in $\Cc_p$ admits a
super-potential $\Uc_{S'}$ such that $\Uc_{S'}\geq \Uc_S$ we say that
$S'$ is {\it more H-regular than $S$} or simply {\it $S'$ is more
  diffuse than $S$}.
\end{definition}

\begin{remarks}\rm \label{rk_hartogs_cv}
By Lemma \ref{lemma_cv_sp_pointwise} below, the property that
$\Uc_{S_n}$ converge pointwise to $\Uc_S$ implies that $m_n\rightarrow m$ and $S_n\rightarrow S$. 
If $S_n$ H-converge to $S$ as in Definition
\ref{def_cv_h}, by Proposition \ref{prop_hartogs_sqp}, $\Uc_{S_n}\rightarrow\Uc_S$ pointwise.
If $\Uc_{S_n}$ decrease to $\Uc_S$, then $S_n$
H-converge to $S$, see also Corollary
\ref{cor_decreasing_sqp} below. 
We have seen in Proposition \ref{prop_regularization_sqp} that
$S_\theta$ H-converge to $S$  when $\theta\rightarrow 0$.
\end{remarks}

\begin{lemma} \label{lemma_cv_sp_pointwise}
Let $(S_n)$ be a sequence in $\Cc_p$ and $\Uc_{S_n}$
 be super-potentials of mean $m_n$  of $S_n$. Assume that $\Uc_{S_n}$
 converge to a finite function $\Uc$ on smooth forms in $\Cc_{k-p+1}$. Then, $m_n$ converge to a
 constant $m$, $S_n$ converge to a current $S$ and $\Uc$ is equal to
 the super-potential of mean $m$ of $S$ on smooth forms in $\Cc_{k-p+1}$. 
\end{lemma}
\proof
We have $m_n=\Uc_{S_n}(\omega^{k-p+1})$. Hence, $m_n$ converge to
$m:=\Uc(\omega^{k-p+1})$. Let $S$ and $S'$ be limit currents of
$(S_n)$. From the definition of super-potential, we deduce that the
super-potentials of mean $m$ of $S$ and of $S'$ are equal to $\Uc$ on smooth
forms in $\Cc_{k-p+1}$. By Proposition \ref{prop_unique_sp},
$S=S'$. Hence, $(S_n)$ is convergent.    
\endproof

We now give a compactness property of super-potentials.

\begin{proposition} \label{prop_compactness_sqp}
Let $\Uc_{S_n}$ be a super-potential of a current
$S_n$ in $\Cc_p$. Assume that $(\Uc_{S_n})$ is bounded from above 
and does not converge
uniformly to $-\infty$. Then, there is an increasing  sequence $(n_i)$ of
integers such that $S_{n_i}$ converge to a current $S$ and $\Uc_{S_{n_i}}$
converge  on smooth forms in $\Cc_{k-p+1}$ to a super-potential
$\Uc_S$ of $S$. 
Moreover, we have $\limsup \Uc_{S_{n_i}}\leq \Uc_S$.
\end{proposition}
\proof
By the last assertion in 
Corollary \ref{cor_def_sp}, since $(\Uc_{S_n})$ is bounded from above and does not converge to
$-\infty$, their means $m_n$ are bounded from above uniformly on $n$
and do not converge to $-\infty$. 
Extracting a subsequence
allows to assume that $S_n$ converge to a current $S$ and $m_n$
converge to a 
finite value $m$. So, we can assume $m_n=m=0$. Let $\Uc_S$ denote the
super-potential of mean 0 of $S$. By definition of
$\Uc_S(R)$ for $R$ smooth, we have
$\Uc_{S_n}(R)\rightarrow \Uc_S(R)$. The inequality  $\limsup \Uc_{S_{n_i}}\leq \Uc_S$
is a consequence of Proposition \ref{prop_hartogs_sqp}. 
\endproof

\begin{corollary} \label{cor_decreasing_sqp}
Let $\Uc_{S_n}$ be super-potentials of mean $m_n$ of $S_n$. Assume
that $\Uc_{S_n}$ decrease to a function $\Uc$ which is not identically $-\infty$. Then, $S_n$
converge to a current $S$, $m_n$ converge to a constant $m$ and $\Uc$
is the super-potential of mean $m$ of $S$.
\end{corollary}
\proof
By Lemma \ref{lemma_cv_sp_pointwise}, 
$S_n$ converge to a current $S$ and $m_n$ converge to a constant  $m$. 
Define $u(\theta):=\Uc(R_\theta)$ and
$u_n(\theta):=\Uc_{S_n}(R_\theta)$.
As in Proposition \ref{prop_regularization_sqp}, the functions $u_n$ are quasi-subharmonic
and decrease to $u$. 
Hence, $u$ is quasi-subharmonic. On the other hand, since $R_\theta$ is
smooth for $\theta\not=0$, we have $u(\theta)=\Uc_S(R_\theta)$ for
$\theta\not=0$ where $\Uc_S$ is the super-potential of mean $m$ of $S$.
The function $\theta\mapsto \Uc_S(R_\theta)$ is also quasi-subharmonic on
$\Delta$. So, we have necessarily $\Uc_S(R)=u(0)=\Uc(R)$.
This holds for every $R$ in $\Cc_{k-p+1}$.  Therefore, $\Uc$ is
the super-potential of mean $m$ of $S$.
\endproof

\begin{corollary} \label{cor_symetry_sqp}
Let $\Uc_S$ and $\Uc_R$ be super-potentials 
of the same mean $m$ of $S$ and $R$ respectively. Then,
$\Uc_S(R)=\Uc_R(S)$.
\end{corollary}
\proof
We have seen in the proof of Lemma \ref{lemma_super_def} that the corollary holds for $S$ smooth.
Let $S_n$ be smooth forms as in Proposition
\ref{prop_regularization_sqp}. The upper semi-continuity implies 
$$\Uc_S(R)=\lim_{n\rightarrow\infty} \Uc_{S_n}(R)=\lim_{n\rightarrow\infty} \Uc_R(S_n)\leq \Uc_R(S).$$
In the same way, we prove that $\Uc_R(S)\leq \Uc_S(R)$. 
\endproof

\begin{lemma} \label{lemma_sp_qp}
Let $S$, $S'$ be currents in $\Cc_p$ and let $\Uc_S$, $\Uc_{S'}$ be
 their super-potentials of mean $m$. Assume there is a
 positive $(p-1,p-1)$-current $U$ such that $\ddc U=S'-S$. Then, 
$\Uc_{S'}+\|U\|\geq \Uc_S$. In particular, if $S$ has bounded
super-potentials, then $S'$ has bounded super-potentials.
If $\Uc_R$ is a super-potential of a current $R\in\Cc_{k-p+1}$, then
$\Uc_R(S')+\|U\|\geq \Uc_R(S)$.
\end{lemma}
\proof
Let $U_S$ be a quasi-potential of mean $m$ of $S$. Then, $U_S+U$ is a
quasi-potential of mean $m+\|U\|$ of $S'$. 
For $R$ smooth, we have
$$\Uc_{S'}(R)+\|U\|=\langle U_S+U,R\rangle\geq \langle U_S,R\rangle =
\Uc_S(R).$$
Then, Corollaries \ref{cor_def_sp} and \ref{cor_symetry_sqp} imply the result.
\endproof

We have the following important result which can be considered as a
version of Lemma \ref{lemma_exp_estimate} for super-potentials. We can
apply it to $K=W=\P^k$.

\begin{proposition} \label{prop_estimate_sqp}
Let $W\subset \P^k$ be an open set and $K\subset W$ be a compact
set. Let $S$ be a current in $\Cc_p$ 
with support in $K$ and $R$ be a current in $\Cc_{k-p+1}$. Assume that
the restriction of $R$ to $W$ is a bounded form.
Then, the super-potential $\Uc_S$ of mean $0$ of $S$ satisfies 
$$|\Uc_S(R)|\leq
c\big(1+\log^+\|R\|_{\infty,W}\big)$$
where $c>0$ is a constant independent of $S$, $R$ and $\log^+:=\max(0,\log)$.
\end{proposition}
\proof
Recall that
$u(\theta):=\Uc_S(R_\theta)$ is a quasi-subharmonic function on
$\Delta$ such that $\ddc u\geq -\alpha$. Proposition \ref{prop_regularization} shows that
the family of these functions $u$ for $(S,R)\in\Cc_p\times\Cc_{k-p+1}$
is compact. So, Lemma \ref{lemma_exp_estimate} implies that $\|e^{-A
  u}\|_{\Lc^1(\Delta_{1/2})}\leq c$ for some positive constants $c$ and
$A$. 

Suppose the estimate in the lemma
is not valid. Recall that $\Uc_S$ is bounded from above by a constant
independent of $S$. Then, for $\epsilon>0$ arbitrary small
there is an $R$ such that $M:=\log\|R\|_{\infty,W}\gg 0$ and
$\Uc_S(R)\leq -2M/\epsilon$. 
It follows that
$u(0)=\Uc_S(R)\leq -2M/\epsilon$.
We will show that
$u(\theta)\leq -M/\epsilon$ on a disc of radius
$e^{-M}$ which contradicts the above estimate on
$e^{-Au}$ for $\epsilon$ small enough.

Let $U$ be the Green quasi-potential of $R$ and $m$ its mean. The mass
of $U$ is bounded by a constant independent of $R$.
By Lemma \ref{lemma_deform_pot}, there is a quasi-potential $U_{\theta
  y}'$ of $R_{\theta y}$ of mean $m$ such that 
$$\|U_{\theta y}'-(\tau_{\theta
  y})_*(U)\|_{\infty}\lesssim |\theta|.$$ 
We deduce that
$$|\Uc_S(R_\theta) -\Uc_S(R)|= \Big|\int_y \langle S, U_{\theta
  y}'-U\rangle d\rho(y)\Big|\lesssim |\theta|+\Big|\int_y \langle S, 
(\tau_{\theta y})_*(U)-U\rangle d\rho(y)\Big| .$$
Because $\theta$ is small, $\tau_{\theta y}^{-1} (K)\subset W'$ for
some fixed open set $W'\Subset W$. Since $\tau_{\theta y}$ is close to the
identity, using Lemma \ref{lemma_levine_c1}, 
we obtain 
$$\|(\tau_{\theta y})_*(U)-U\|_{\infty, K}\lesssim
|\theta|\|U\|_{\Cc^1(W')} \lesssim |\theta|e^M.$$
Therefore,
$$|u(\theta)-u(0)|=|\Uc_S(R_\theta) -\Uc_S(R)|\lesssim |\theta|e^M.$$
This implies the above claim and completes the proof.
\endproof


\subsection{Currents with regular super-potentials}

The PB or PC currents are introduced in \cite{DinhSibony1,
DinhSibony4, DinhSibony6} in the study of holomorphic dynamical
systems. They
correspond to currents with bounded or continuous
super-potentials. We recall first the definition of the space $\DSH^{k-p}(\P^k)$ of dsh currents. A real
$(k-p,k-p)$-current $\Phi$ of finite mass is {\it dsh} if there are positive closed
currents $R^\pm$ of bidegree $(k-p+1,k-p+1)$ such that\footnote{It is
  also useful to consider the space generated 
by such currents $\Phi$ which are negative. This
is necessary in order to defined the pull-back of DSH currents by
holomorphic maps.}
$\ddc\Phi=R^+-R^-$. Define 
$$\|\Phi\|_\DSH:=\|\Phi\|+\min \|R^\pm\|$$
with $R^\pm$ as above. 
We consider a {\it weak topology} on $\DSH^{k-p}(\P^k)$. A sequence
$(\Phi_n)$ converges to $\Phi$ in $\DSH^{k-p}(\P^k)$ if
$\Phi_n\rightarrow\Phi$ in the sense of currents and $\|\Phi_n\|_\DSH$
is uniformly bounded. A positive closed $(p,p)$-current $S$ is
said to be {\it PB} if there is a constant $c>0$ such that 
$$|\langle S,\Phi\rangle|\leq c\|\Phi\|_\DSH$$
for smooth real forms $\Phi$ of bidegree $(k-p,k-p)$. We say that
$S$ is {\it PC} if it can be extended to a linear form on
$\DSH^{k-p}(\P^k)$ which is continuous with respect to the weak
topology on $\DSH^{k-p}(\P^k)$.

\begin{proposition}
If a super-potential $\Uc_S$ of $S$ is finite everywhere, then it is
bounded.
A current $S$ is PB if and only if the super-potentials of $S$ are  bounded. 
A current $S$ is PC if and only if the super-potentials of $S$ are  continuous. 
\end{proposition}
\proof
Subtracting a constant from $\Uc_S$, we can assume $\Uc_S\leq 0$.
Assume that $\Uc_S$ is unbounded. Then, there are currents $R_n$ such
that $\Uc_S(R_n)\leq -2^n$. Define $R:=\sum 2^{-n} R_n$. Since $\Uc_S$
is affine and negative, we have $\Uc_S(R)\leq \sum_{n\leq N}
2^{-n}\Uc_S(R_n)$ for every $N$. Hence, $\Uc_S(R)=-\infty$. This is a
contradiction. So, $\Uc_S$ is bounded. Note that this property is
false for quasi-psh functions on $\P^k$. 

Assume that the super-potential $\Uc_S$ of mean 0 of $S$ satisfies 
$|\Uc_S|<M$ for some constant $M>0$. Consider a real smooth form $\Phi$ of bidegree $(k-p,k-p)$
and a constant $A\geq \|\Phi\|_\DSH$.  We will prove that $|\langle
S,\Phi\rangle|\leq  A(1+2C+2M)$ with $C>0$ independent of $S$. 
This implies that $S$ is PB. Since we can approximate $S$ in the
Hartogs' sense by smooth forms,
it is enough to prove this inequality for
$S$ smooth. Write $\ddc
\Phi=A(R^+-R^-)$ with $\|R^\pm\|=1$. By Remark \ref{rk_levine}, 
there are quasi-potentials $U^\pm$ of mean 0 of
$R^\pm$ such that $\|U^\pm\|_\DSH\leq C$ where $C>0$ is a
constant. Define $\Psi:=\Phi-AU^++AU^-$. We have $\ddc\Psi=0$
and 
$$\|\Psi\|\leq \|\Phi\|+A\|U^+\|+A\|U^-\|\leq A(1+2C).$$
Since $\ddc\Psi=0$ and since $S$ is cohomologous to $\omega^p$, we
have
$$|\langle S,\Psi\rangle|=|\langle \omega^p,\Psi\rangle|\leq
A(1+2C).$$
It follows that 
\begin{eqnarray*}
|\langle S,\Phi\rangle| & \leq & |\langle S,\Psi\rangle|+A|\langle
S,U^+\rangle|+A|\langle S,U^-\rangle|\\
& = & |\langle S,\Psi\rangle|+A|\Uc_S(R^+)|+A|\Uc_S(R^-)|\\
&\leq & A(1+2C+2M).
\end{eqnarray*}
Hence, $S$ is PB.

Conversely, if  $S$ is PB, we show that $\Uc_S$ is bounded. Consider
a smooth form $R$ in $\Cc_{k-p+1}$. Let $U_R$ be a 
quasi-potential of $R$ of mean 0 such that
$\|U_R\|_\DSH\leq C$. We have $\Uc_S(R)=\langle S,U_R\rangle$. Since
$S$ is PB, $\Uc_S(R)$ is bounded by a constant independent of
$R$. This implies that $\Uc_S$ is bounded. 

It is clear that if $S$ is PC, $\langle S,U_R\rangle$ for $R$ smooth 
can be extended to a continuous function on $\Cc_{k-p+1}$. Indeed, we can
choose $U_R$ depending continuously on $R$ with respect to the
weak topology in $\DSH^{k-p}(\P^k)$, see Theorem \ref{thm_ddbar} and
Remark \ref{rk_levine}. This implies that $\Uc_S$ is
continuous. Conversely, if $\Uc_S$ is continuous, we show that $S$ is
PC. If $\Phi$ and $R^\pm$ are smooth as above, we obtain
$$\langle S,\Phi\rangle =\langle
\omega^p,\Psi\rangle+A\Uc_S(R^+)-A\Uc_S(R^-).$$
The right hand side depends on $\Psi$ and on
$AR^+-AR^-=\ddc\Phi$ but not on the choice of $A,R^\pm$. Hence, since $\Psi$ and $\ddc\Phi$ depend
continuously on $\Phi$, we can
extend $S$ to a continuous linear form on $\DSH^{k-p}(\P^k)$. The
continuity is
with respect to the weak topology on $\DSH^{k-p}(\P^k)$.   This
completes the proof.
\endproof

\begin{lemma} 
If $S$ is a form of class $\Lc^s$ with $s>k$, then $S$
has continuous super-potentials.
\end{lemma}
\proof
Let $r$ be the positive number such that $1/r+1/s=1$. Then,
$r<k/(k-1)$. The Green quasi-potential $U_R$ of $R$ 
is a form of class $\Lc^r$. Moreover, with respect to the $\Lc^r$ topology, it
depends continuously on $R$, see Theorem \ref{thm_ddbar}. 
The mean $m_R$ of $U_R$ depends continuously on $R$.
On the other hand, the super-potential of mean 0 of $S$ satisfies  
$$\Uc_S(R)=\langle S, U_R\rangle-m_R$$
for $R$ smooth.
The right hand side  is defined for every $R$ and
depends continuously on $R$. Therefore, $\Uc_R$ is continuous.
\endproof

\begin{remark}\rm
$U_R$ is in the Sobolev space $W^{1,r}$ with $r<2k/(2k-1)$. So, we can
assume $S\in W^{-1,s}$ with $1/r+1/s=1$, and still $\Uc_S$ is continuous.
\end{remark}

\begin{proposition} \label{prop_compare_sp}
Let $S$ and $S'$ be currents in $\Cc_p$ such that $S'\leq c S$ for
some positive constant $c$. If $S$ has bounded super-potentials, then
$S'$ has bounded super-potentials.  If $S$ has continuous super-potentials, then
$S'$ has continuous super-potentials.
\end{proposition}
\proof
Write $S=\lambda S'+(1-\lambda)S''$ with $0<\lambda\leq 1$ and $S''$ a
current in $\Cc_p$. Let $\Uc_S$, $\Uc_{S'}$ and $\Uc_{S''}$ denote the super-potentials
of mean 0 of $S$, $S'$,  and $S''$. By definition of super-potentials,
we have 
$\lambda \Uc_{S'} + (1-\lambda) \Uc_{S''}=\Uc_S$ on smooth forms
$R$. Corollary \ref{cor_def_sp} 
implies that this equality holds for
every $R$. Since $\Uc_{S''}$ is bounded from above, if $\Uc_S$ is
bounded, it is clear that $\Uc_{S'}$ is bounded. If $\Uc_S$ is
continuous, since $\Uc_{S'}$ and $\Uc_{S''}$ are u.s.c., they are 
continuous.
\endproof

\begin{proposition} Let $S$ be a current with bounded
super-potentials. Then, $S$ has no mass on pluripolar sets of
$\P^k$. In particular, $S$ does not give mass to proper analytic subsets of $\P^k$.
\end{proposition}
\proof
Assume that $S$ has bounded super-potentials. Let $E\subset
\P^k$ be a pluripolar set and $u$ be a  quasi-psh
function such that $\ddc u\geq -\omega$ and $E\subset
\{u=-\infty\}$. Define $R:=(\ddc u+\omega)\wedge \omega^{k-p}$. This
is a current in $\Cc_{k-p+1}$. 

Let $(u_n)$ be a sequence of smooth
functions decreasing to $u$ and such that $\ddc u_n\geq -\omega$. Define
$R_n:=(\ddc u_n+\omega)\wedge \omega^{k-p}$. 
Observe that $u_n\omega^{k-p}$ are quasi-potentials of mean $m_n:=\int u_n\omega^k$
of $R_n$. 
If $\Uc_S$ is the super-potential of mean $m:=\int u\omega^k$ of $S$, then 
$\langle S, u_n\omega^{k-p}\rangle$ decrease to $\Uc_S(R)$.
Hence, $\Uc_S(R)=\langle S, u\omega^{k-p}\rangle$.
Since $S$ has bounded super-potentials, 
$\langle S,u\omega^{k-p}\rangle$ is finite. It follows that $S$ has no
mass on $\{u=-\infty\}$. 
\endproof

\begin{proposition} \label{prop_holder_sp}
Assume that $S$ admits a super-potential which is $\alpha$-H{\"o}lder
continuous with respect to the distance $\dist_1$ on $\Cc_{k-p+1}$ for
some exponent $\alpha\leq 1$. 
Let $\sigma_S$ denote the trace measure of $S$. There is a
constant $c>0$ such that if $B_r$ is a ball of radius $r$, then 
$\sigma_S(B_r)\leq c r^{2k-2p+\alpha}$. In particular,
$S$ has no mass on 
Borel subsets of $\P^k$ with Hausdorff dimension less than $2(k-p)+\alpha$.
\end{proposition}

Using Lemma \ref{lemma_compare_dist}, we deduce analogous results for a general distance
$\dist_\beta$ on $\Cc_{k-p+1}$.
Note that 
the last assertion in the proposition is deduced from the first one and some classical
arguments.
In order to prove the first assertion
It is enough to consider $r$ small.
So, we can assume that $B_r$ is a ball of center $0$ in an affine
chart $\C^k\subset\P^k$. It is sufficient to show that $\int_{\Delta_r^k} S\wedge
\omega^{k-p}\lesssim r^{2k-2p+\alpha}$. Let $z$ denote the
canonical coordinates in $\C^k$.

\begin{lemma}
There are positive constants $A$, $c$ independent of $r$, a positive
$(k-p,k-p)$-current $\Phi$
and two currents $R^\pm$ in $\Cc_{k-p+1}$ such that $\Phi\geq
(\ddc|z|^2)^{k-p}$ on $\Delta_r^k$,
$\|\Phi\| \leq A r^{2k-2p+2}$,
$\ddc\Phi = c
r^{2k-2p}(R^+-R^-)$ and $\dist_1(R^+,R^-)\leq A r$. 
\end{lemma} 
\proof
Observe that $(\ddc|z|^2)^{k-p}$ is a combination of the forms
$$(i dz_{i_1}\wedge d\overline z_{i_1})\wedge \ldots \wedge (i d
z_{i_{k-p}}\wedge d\overline z_{i_{k-p}}).$$
Without loss of generality, one only has to construct a $\Phi$,
$R^\pm$ satisfying the last three properties in the lemma and the inequality
$$\Phi\geq (i dz_1\wedge d\overline z_1)\wedge \ldots \wedge (i d
z_{k-p}\wedge d\overline z_{k-p})$$ 
on $\Delta_r^k$. Taking a combination of such currents 
gives currents satisfying the lemma.

Let $\chi$ be a smooth cut-off function
with compact support in $\Delta_2^k$, equal to 1 on $\Delta_1^k$. Let $v(z_{k-p+1})$ be a smooth function
with support in $\{|z_{k-p+1}|<2r\}$ such that $0\leq v\leq 1$,
$\|v\|_{\Cc^1}\lesssim r^{-1}$,
$\|v\|_{\Cc^2}\lesssim r^{-2}$ and $v=1$ on $\{|z_{k-p+1}|\leq r\}$.
Let $\pi:\C^k\rightarrow\C^{k-p}$ and $\pi':\C^k\rightarrow \C^{k-p+1}$
denote the canonical projections on the first factors of
$\C^k$. Consider  the
restriction $\Theta$ of $idz_1\wedge d\overline z_1\wedge
\ldots\wedge idz_{k-p}\wedge d\overline z_{k-p}$ to $\Delta_r^{k-p}$
and define
$$\Phi:=v(z_{k-p+1})\chi(z)
\pi^*(\Theta).$$
Then, $\Phi$ satisfies the desired lower estimate on $\Delta_r^k$. We have to check the
last three properties in the lemma.

 Since $\pi$ can be
extended to a rational map from $\P^k$ to $\P^{k-p}$,
$\pi^*(\Theta)$ can be extended to a positive closed current on $\P^k$
of mass $\|\Theta\|\simeq r^{2k-2p}$. Observe also
that Cauchy-Schwarz's inequality implies that 
$$\ddc [v(z_{k-p+1})\chi(z)]\lesssim r^{-2}
idz_{k-p+1}\wedge d\overline z_{k-p+1}+\omega.$$
Denote by $\Theta'$ the restriction of  $idz_1\wedge d\overline z_1\wedge
\ldots\wedge idz_{k-p+1}\wedge d\overline z_{k-p+1}$ to
$\Delta_r^{k-p}\times\Delta_{2r}$ and let
$\Omega^-:=\lambda{\pi'}^*(r^{-2}\Theta')+\lambda\omega\wedge
\pi^*(\Theta)$
with $\lambda>0$ large enough independent of $r$. Then, $\Omega^+:=\Omega^-+\ddc\Phi$ is
positive and closed. We have $\ddc\Phi=\Omega^+-\Omega^-$. The currents
$\Omega^\pm$ can be extended to positive closed currents on
$\P^k$. They have the same mass since they are cohomologous. This mass
is of order $r^{2k-2p}$ and we denote it by $cr^{2k-2p}$.
We obtain $\ddc\Phi=  c r^{2k-2p}(R^+-R^-)$ with
$R^\pm:=c^{-1} r^{2p-2k}\Omega^\pm$. The currents $R^\pm$ are in $\Cc_{k-p+1}$.
We want to bound $\dist_1(R,R')$. For any test form $\Psi$ with
$\|\Psi\|_{\Cc^1}\leq 1$, we have
$$|\langle R-R',\Psi\rangle| \simeq   r^{2p-2k} |\langle \ddc\Phi,
\Psi\rangle|
= r^{2p-2k} |\langle d^c\Phi, d\Psi\rangle| \lesssim 
r^{2p-2k} \|d^c\Phi\|.$$
On the other hand, we deduce from the definition of $\Phi$ that
$$\|\dc\Phi\|\lesssim r^{2k-2p}\|\dc v\|_{\Delta_{2r}}\lesssim r^{2k-2p+1}.$$
This implies the result.
\endproof

\noindent
{\bf End of the proof of Proposition \ref{prop_holder_sp}.}
Let $\Uc_S$ be a super-potential of $S$.
Since $\Uc_S$ is $\alpha$-H{\"o}lder continuous, we deduce from the previous lemma that
\begin{eqnarray*}
\int_{\Delta_r^k} S\wedge \omega^{k-p} & \leq & \langle S,\Phi\rangle 
= \langle \omega^p,\Phi\rangle +\langle \ddc U_S,\Phi\rangle \\
& = & \langle \omega^p,\Phi\rangle +\langle U_S,\ddc\Phi\rangle \\
& \lesssim &
\langle \omega^p,\Phi\rangle +r^{2k-2p}(\Uc_S(R^+)-\Uc_S(R^-))\\
& \lesssim & 
r^{2k-2p+\alpha}.
\end{eqnarray*}
This is the required
estimate.
\hfill $\square$


\subsection{Capacity of currents and super-polar sets}

We will define a notion of capacity for Borel subsets $E$ of $\Cc_{k-p+1}$.
This capacity does not describe how
``big'' is the set $E$ but rather how singular are the currents in
$E$. The definition mimicks the notion of capacity that we introduced in
\cite{DinhSibony6} for compact K{\"a}hler manifolds.
Let 
$$\Pc_p:=\big\{\Uc_S \mbox{ super-potential of } S\in\Cc_p, \quad
\max_{\Cc_{k-p+1}}\Uc_S=0\big\}.$$

\begin{definition} \rm
We define {\it the capacity} of $E$ to be the following quantity
$$\capacity(E):=\inf_{\Uc\in\Pc_p}\exp \big(\sup_{R\in E}\Uc(R)\big).$$
\end{definition}
It is clear that the capacity is increasing as a set function.
Propositions \ref{prop_regularization_sqp} and \ref{prop_hartogs_sqp} imply that when $E$ is
compact, in the previous definition we obtain
the same capacity if we only consider super-potentials of smooth forms. 
We also have  $\capacity(\Cc_{k-p+1})=1$ and it follows that the set
of smooth forms in $\Cc_{k-p+1}$ has capacity 1. 
Dense
subsets of smooth forms in $\Cc_{k-p+1}$ have also capacity 1. So, there is a
countable subset of $\Cc_{k-p+1}$ with capacity 1. 
\begin{definition}
\rm
We say that $E$ is {\it super-polar} or {\it complete
super-polar} in $\Cc_{k-p+1}$  if there is a
super-potential $\Uc_S$ of a current $S$ in $\Cc_p$ such that $E\subset \{\Uc_S=-\infty\}$
or $E=\{\Uc_S=-\infty\}$ respectively. 
\end{definition}

Let  $\widehat E$ be the {\it barycentric hull} of $E$, i.e. the set of
currents $\int R d\nu(R)$ where $\nu$ is a probability measure on
$\Cc_{k-p+1}$ such that $\nu(E)=1$.
Denote by $\widetilde E$ the
set of currents $cR+(1-c)R'$ with $R\in \widehat E$, $R'\in \Cc_{k-p+1}$ and $0<c\leq
1$. Then, $\widetilde E$ and $\widehat
E$ are convex.

\begin{proposition} \label{prop_convex_hull}
The following properties are equivalent
\begin{enumerate}
\item $E$ is super-polar in $\Cc_{k-p+1}$.
\item $\widehat E$ is super-polar in $\Cc_{k-p+1}$.
\item  $\widetilde E$ is super-polar in $\Cc_{k-p+1}$.
\item $\capacity(E)=0$.
\end{enumerate}
Moreover,  a countable union of super-polar sets is super-polar, complete super-polar sets are convex and 
$\capacity(E)=\capacity(\widehat E)$.
\end{proposition}
\proof
Since every function $\Uc$ in $\Pc_p$ is affine and negative, 
if $\Uc$ is equal to $-\infty$ on $E$, it is also equal to
$-\infty$ on $\widehat E$ and $\widetilde E$. 
Therefore, the first three properties are equivalent. 
We also deduce that  if $E$ is complete super-polar, then $E$ is
convex and $E=\widetilde E$.
Moreover, for any $\Uc$ we have $\sup_E\Uc=\sup_{\widehat E}\Uc$. This implies that
$\capacity(E)=\capacity(\widehat E)$.

It is clear that if $E$ is
super-polar, then $\capacity(E)=0$.  Assume that
 $\capacity(E)=0$. We show that $E$ is
super-polar.
There are super-potentials $\Uc_{S_n}$ of
 $S_n$ such that $\max \Uc_{S_n}=0$ and $\Uc_{S_n}\leq -2^n$ on
$E$. Corollary \ref{cor_def_sp} implies that the means of $\Uc_{S_n}$
are bounded. 
This and Corollary \ref{cor_decreasing_sqp} imply that 
$\Uc=\sum_{n\geq 1} 2^{-n} \Uc_{S_n}$ is a super-potential of 
$\sum_{n\geq 1} 2^{-n} S_n$. It is equal to $-\infty$ on $E$. 
Hence, $E$ is super-polar. A similar argument implies that a countable union
of super-polar sets is super-polar. 
\endproof

\begin{proposition}
Let $E\subset \Cc_{k-p+1}$ be a compact set. Then, 
$E$ has positive capacity if and only if its barycentric hull contains a current with bounded
super-potentials. Moreover, 
for every $\epsilon>0$ there is a current $R$ in the barycentric hull
$\widehat E$ of $E$  such that its super-potential of mean $0$ satisfies
$$\Uc_R \geq \log\capacity(E) -\epsilon \mbox{ on } \Cc_p.$$
\end{proposition}
\proof
If $R$ is a current with bounded super-potentials, then by symetry $\Uc(R)\not=-\infty$
for every $\Uc\in\Pc_p$. Proposition \ref{prop_convex_hull} implies that $\{R\}$ is not
super-polar. Hence,
if $\widehat E$ contains
a current with bounded super-potentials, $\widehat E$ has
positive capacity. Proposition \ref{prop_convex_hull} also implies that $E$
has positive capacity. 
Now, assume that $E$ has positive capacity. 
We show that $\widehat E$ contains
a current with bounded super-potentials.
In what follows, the symbol $\Uc$ denotes a super-potential of mean 0.
We have
$$\inf_{S\in\Cc_p} \sup_{R\in\widehat E} \Uc_S(R)\geq M:=\log\capacity(E).$$
The function $\Uc_S(R)$ is affine in both variables $R$ and $S$. Hence, for every
convex compact set $\Cc$ of continuous forms in $\Cc_p$, the
minimax theorem \cite{Sion} implies 
$$\sup_{R\in \widehat E}\inf_{S\in\Cc}
\Uc_S(R)=\inf_{S\in\Cc}\sup_{R\in\widehat E} \Uc_S(R)\geq M.$$
Observe that the family of all the convex compact sets $\Cc$ is
ordered by inclusion and if $\Cc,\Cc'$ are such two sets, we
have $\max(\Cc,\Cc')=\widehat{\Cc\cup\Cc'}$.
Define
$$E_\Cc:=\big\{R\in\widehat E,\ \Uc_S(R)\geq M-\epsilon\quad \mbox{for every } S\in \Cc\big\}.$$ 
So, $(E_\Cc)$ is an ordered family of compact sets in $\widehat E$ and
we can apply Zorn's lemma. Take an element $R$ in a minimal set $\Cc^0$.
Then, 
$\Uc_R(S)=\Uc_S(R)\geq M-\epsilon$ for
every continuous $S$, if not we consider $\widehat{\Cc^0\cup\{S\}}$. 
This completes the proof.
\endproof

Consider the set of the super-potentials $\Uc$ of mean 0 of currents in
$\Cc_p$ and define  
$c_{k,p}:=\sup_{\Cc_p}\max \Uc$.
Corollary \ref{cor_def_sp} implies that this constant is finite.

\begin{corollary} \label{cor_cap_r}
For every current $R$ in $\Cc_{k-p+1}$, if $\Uc_R$ is the
super-potential of mean $0$ of $R$, then 
$$\log\capacity(R)\geq 
\inf_{\Cc_p}-c_{k,p} +\Uc_R.$$
\end{corollary}
\proof
Let $\Uc_S$ be the super-potential of mean 0 of $S$. By definition of
capacity and of $c_{k,p}$, we have
$$\log\capacity(R)\geq \big[
\inf_{S\in\Cc_p}\Uc_S(R)-c_{k,p}\big].$$
Corollary \ref{cor_symetry_sqp} implies the result.
\endproof

\begin{corollary} For every $r>k$, there is a constant $c>0$ such that
if $R$ is a form in $\Cc_{k-p+1}$ with coefficients in  $\Lc^r$, then
$$\log\capacity(R)\geq -c_{k,p}-c\|R\|_{\Lc^r}.$$
\end{corollary}
\proof
Let $s$ be the positive number such that $1/r+1/s=1$. Then, $s<k/(k-1)$.
Let $U_S$ be the Green quasi-potential of $S$. This is a negative
form with $\Lc^s$ norm bounded uniformly on $S$. Hence
$$\Uc_R(S)\geq \langle U_S,R\rangle \geq -c \|R\|_{\Lc^r}$$
for some constant $c>0$. We obtain the result from Corollary \ref{cor_cap_r}.
\endproof

The following result is a consequence of Proposition \ref{prop_estimate_sqp}.

\begin{corollary} There are constants $c>0$ and $\lambda>0$ such that
  for every bounded form $R$ in $\Cc_{k-p+1}$
$$\capacity(R)\geq c\|R\|_\infty^{-\lambda}.$$
\end{corollary}


\section{Theory of intersection of currents} \label{section_intersection}

In this paragraph, we develop the theory of intersection for positive
closed currents of arbitrary bidegree. The method can be extended to
currents on compact K{\"a}hler manifolds or in some local
situation, see also \cite{DinhSibony7}. Here, for simplicity, we only consider 
currents in the projective space.


\subsection{Some universal super-functions}

Let $p$ be a positive integer such that
$1\leq p \leq k$. Define a universal function $\Uc_p$ on
$\Cc_p\times\Cc_{k-p+1}$ by
$$\Uc_p(S,R):=\Uc_S(R)=\Uc_R(S)$$
where $\Uc_S$ and $\Uc_R$ are super-potentials of mean 0 of $S$ and
$R$, see Corollary \ref{cor_symetry_sqp}.
We have seen that when $S$ is fixed $\Uc_p$ is quasi-psh on
special varieties of $\Cc_{k-p+1}$ and when $R$ is fixed, it  is quasi-psh on
special varieties of $\Cc_p$.

\begin{lemma} \label{lemma_universal}
The function $\Uc_p$ is
 is u.s.c. on $\Cc_p\times\Cc_{k-p+1}$. 
\end{lemma}
\proof
Let $S_n$ be currents in $\Cc_p$ converging to $S$ and $R_n$ in
$\Cc_{k-p+1}$ converging to $R$. Let $\Uc_{S_n}$  denote
the super-potential of mean $0$ of $S_n$.
Choose $\Uc$ continuous with $\Uc_S<\Uc$. By Proposition \ref{prop_hartogs_sqp}, for $n$ large enough,
$\Uc_{S_n}<\Uc$ and hence $\Uc_{S_n}(R_n)< \Uc(R_n)$. We then
get 
$$\limsup_{n\rightarrow\infty} \Uc_{S_n}(R_n)\leq \Uc(R).$$
Since $\Uc$ is arbitrary, we deduce that
$$\limsup_{n\rightarrow\infty} \Uc_{S_n}(R_n)\leq \Uc_S(R).$$
This proves the lemma.
\endproof

\begin{lemma} \label{lemma_compare_sf}
Let $S'$, $R'$ be currents in $\Cc_p$ and $\Cc_{k-p+1}$, and $\Uc_{S'}$,
$\Uc_{R'}$ be their super-potentials of mean $0$. Assume there are constants
$a,b$ such that $\Uc_{S'}+a\geq \Uc_S$ and $\Uc_{R'}+b\geq \Uc_R$. Then,
$\Uc_p(S',R')\geq \Uc_p(S,R)-a-b.$
\end{lemma}
\proof
We have 
$$\Uc_p(S,R')=\Uc_{R'}(S)\geq \Uc_R(S)-b=\Uc_p(S,R)-b$$
and
$$\Uc_p(S',R')=\Uc_{S'}(R')\geq \Uc_S(R')-a=\Uc_p(S,R')-a.$$
This implies the result.
\endproof

\begin{lemma} \label{lemma_decreas_sf}
Let $(S_n)_{n\geq 0}$ and  $(R_n)_{n\geq 0}$ 
be sequences of currents in $\Cc_p$ 
and $\Cc_{k-p+1}$ H-converging to $S$ and $R$ respectively.
Then, $\Uc_p(S_n,R_n)$ converge to $\Uc_p(S,R)$.
Moreover, if $\Uc_p(S,R)$ is finite, then
$\Uc_p(S_n,R_n)$ is finite for every $n$.
\end{lemma}
\proof
Let $\Uc_{S_n}$ and  $\Uc_{R_n}$  be the
super-potentials of mean 0 of $S_n$ and $R_n$. The H-convergence 
implies the existence of
constants $a_n$ and $b_n$ with limit 0, such that
$\Uc_{S_n}+a_n\geq \Uc_S$ and $\Uc_{R_n}+b_n\geq \Uc_R$.
It follows from Lemma \ref{lemma_universal} that
$$\limsup_{n\rightarrow\infty}  \Uc_p(S_n,R_n)\leq \Uc_p(S,R).$$
It is sufficient to prove that 
$$\Uc_p(S_n,R_n)\geq \Uc_p(S,R)-a_n-b_n.$$
This is a  consequence of Lemma \ref{lemma_compare_sf}.
\endproof


\subsection{Intersection of currents} \label{section_wedge}

Let $p_i$, $1\leq i\leq l$, be positive integers such that
$p_1+\cdots+p_l\leq k$. 
Let $R_i$ be currents in $\Cc_{p_i}$ with $1\leq i\leq l$. We want to
define the wedge-product $R_1\wedge \ldots\wedge R_l$, as a current. In general, one
cannot define this product in a consistent way; for example, when
$R_1$ and $R_2$ are currents of integration on the same projective line
of $\P^2$.  We will define the intersection of the $R_i$
when they satisfy a quite natural condition.
Consider first the case of two currents, i.e. $l=2$.

\begin{proposition} \label{prop_hyp_wedge}
The following conditions are equivalent and are symmetric on $R_1$ and
$R_2$:
\begin{enumerate}
\item $\Uc_{p_1}(R_1,R_2\wedge \Omega)$ is finite for at least one smooth
      form $\Omega$ in $\Cc_{k-p_1-p_2+1}$.
\item $\Uc_{p_1}(R_1,R_2\wedge \Omega)$ is finite for every smooth
      form $\Omega$ in $\Cc_{k-p_1-p_2+1}$.
\item There are sequences $(R_{i,n})_{n\geq 0}$ in
      $\Cc_{p_i}$ converging to 
      $R_i$  and a smooth form $\Omega$ in
      $\Cc_{k-p_1-p_2+1}$ such that $\Uc_{p_1}(R_{1,n},R_{2,n}\wedge\Omega)$ is bounded.
\end{enumerate}
\end{proposition}
\proof
It is clear that the second condition implies the third one: we can
choose $R_{i,n}=R_i$; and the
third condition implies the first one because $\Uc_{p_1}$ is u.s.c. 
Assume the first condition. We show that
$\Uc_{p_1}(R_1,R_2\wedge \Omega')$ 
is finite for every smooth
      form $\Omega'$ in $\Cc_{k-p_1-p_2+1}$. Write
      $\Omega'-\Omega=\ddc U$ with $U$ smooth. Adding
      to $U$ a large positive closed form, we can assume that $U$
      is positive. If $V$ is a quasi-potential of $R_2\wedge \Omega$,
      then the quasi-potential  $V+R_2\wedge U$ of $R_2\wedge \Omega'$ is
      larger than $V$.
Lemmas \ref{lemma_sp_qp} and \ref{lemma_compare_sf} imply that
      $\Uc_{p_1}(R_1,R_2\wedge \Omega')$ is finite. Therefore, the three
      previous conditions are equivalent.

It remains to prove that the first condition is symmetric. 
We can assume $\Omega=\omega^{k-p_1-p_2+1}$. 
Consider the case where $R_1$ 
is smooth. If $U_2$ is a quasi-potential of mean
$0$ of $R_2$, then $U_2\wedge \Omega$ is a quasi-potential of mean 0 of
$R_2\wedge\Omega$. We have
$$\Uc_{p_1}(R_1,R_2\wedge \Omega)=\langle
R_1,U_2\wedge\Omega\rangle=\langle U_2, R_1\wedge\Omega\rangle =\Uc_{p_2}(R_2,R_1\wedge\Omega).$$ 
Suppose now that $R_1$ is arbitrary. Let $R_{1,\theta}$ be the smooth forms
constructed in Paragraph \ref{section_topology}, starting with the current $R_1$, we have
\begin{eqnarray*}
\Uc_{p_1}(R_1,R_2\wedge \Omega) & = & \lim_{\theta\rightarrow 0} \Uc_{p_1}(R_{1,\theta},R_2\wedge \Omega)\\
& = & \lim_{\theta\rightarrow
  0}\Uc_{p_2}(R_2,R_{1,\theta}\wedge\Omega)\\
& \leq & 
\Uc_{p_2}(R_2,R_1\wedge\Omega)
\end{eqnarray*}
since $\Uc_{p_2}$ is u.s.c. In the same way, we obtain
$\Uc_{p_2}(R_2,R_1\wedge\Omega)\leq \Uc_{p_1}(R_1,R_2\wedge
\Omega)$. Hence, $\Uc_{p_2}(R_2,R_1\wedge\Omega)= \Uc_{p_1}(R_1,R_2\wedge
\Omega)$. This implies the symmetry of the first condition in the proposition.
\endproof

\begin{definition} \rm
We say that $R_1$ and $R_2$ are {\it wedgeable} if they satisfy the
conditions in Proposition \ref{prop_hyp_wedge}.
\end{definition}

Note that for $R_1$ fixed, the set of $R_2$ such that $R_1$ and $R_2$
are not wedgeable is a super-polar set in $\Cc_{p_2}$. Indeed, this is
the set of $R_2$ such that $\Uc(R_2)=-\infty$ where $\Uc$ is a
super-potential of $R_1\wedge \omega^{k-p_1-p_2+1}$. So, $R_1$ is
wedgeable for every $R_2$ if and only if $R_1\wedge
\omega^{k-p_1-p_2+1}$ has bounded super-potentials.

\begin{proposition} \label{prop_criteria_wedge}
Let $R_i$ and $R_i'$ be currents in $\Cc_{p_i}$. Assume that $R_1$ and
$R_2$ are wedgeable.
Then,  $R_1'$ and $R_2'$ are wedgeable
 in the following cases:
\begin{enumerate}
\item $R_i'$ is more diffuse than $R_i$ for $i=1,2$.
\item  There is a constant $c>0$ such that $R_i'\leq c R_i$ for $i=1,2$. 
\end{enumerate}
\end{proposition}
\proof
The first assertion is a consequence of Lemma \ref{lemma_compare_sf}. For the second one,
it is enough to show that $R_1$ and $R_2'$ are wedgeable.
Then, in the same way, $R_1'$ and $R_2'$ are wedgeable.
Write $R_2=\lambda R_2'+(1-\lambda)R_2''$ with
$0<\lambda\leq 1$ and $R_2''\in\Cc_{p_2}$. We obtain from the fact
that $\Uc_{p_1}$ is affine that
\begin{eqnarray*}
\lambda \Uc_{p_1}(R_1,R_2'\wedge\omega^{k-p_1-p_2+1}) & = & 
\Uc_{p_1}(R_1,R_2\wedge \omega^{k-p_1-p_2+1})\\
& & -(1-\lambda)\Uc_{p_1}(R_1,R_2''\wedge\omega^{k-p_1-p_2+1})\\
& \not= & -\infty
\end{eqnarray*}
since $\Uc_{p_1}(R_1,R_2\wedge\omega^{k-p_1-p_2+1})\not=-\infty$ and
$\Uc_{p_1}$ is bounded from above. This proves the property.
\endproof

Assume that  $R_1$ and $R_2$ are wedgeable.
We define the wedge-product (or the intersection)
$R_1\wedge R_2$. 
This will be a current of bidegree $(p_1+p_2,p_1+p_2)$.
For every  smooth real form $\Phi$ of bidegree 
$(k-p_1-p_2,k-p_1-p_2)$, write
$\ddc\Phi=c(\Omega^+-\Omega^-)$
where $\Omega^\pm$ are smooth forms in $\Cc_{k-p_1-p_2+1}$
and $c$ is a positive constant.
First, consider the case where $R_1$ or $R_2$ is smooth. 
So, $R_1\wedge R_2$ is defined.
Let $U_1$ be a quasi-potential of mean 0 of $R_1$. Choose $U_1$
smooth if $R_1$ is smooth.
We have
\begin{eqnarray*}
\langle R_1\wedge R_2,\Phi\rangle & = & \langle \omega^{p_1}\wedge R_2,\Phi\rangle +
\langle (R_1-\omega^{p_1})\wedge R_2,\Phi\rangle\\
& = & \langle R_2,\omega^{p_1}\wedge\Phi\rangle +
\langle \ddc(U_1\wedge R_2),\Phi\rangle\\
& = & \langle R_2,\omega^{p_1}\wedge\Phi\rangle +
\langle U_1\wedge R_2,\ddc\Phi\rangle\\
& = & \langle R_2,\omega^{p_1}\wedge\Phi\rangle +
c\Uc_{p_1}(R_1,R_2\wedge\Omega^+)-c\Uc_{p_1}(R_1,R_2\wedge\Omega^-).
\end{eqnarray*}
We deduce that the last expression is independent of the choice of $c$
and $\Omega^\pm$. 
This formally justifies the following formula for wedgeable $R_1$ and $R_2$.
Define
\begin{equation} \label{eq_wedge}
\langle R_1\wedge R_2,\Phi\rangle  :=   \langle R_2,\omega^{p_1}\wedge\Phi\rangle +
c\Uc_{p_1}(R_1,R_2\wedge\Omega^+)-c\Uc_{p_1}(R_1,R_2\wedge\Omega^-).
\end{equation}
The following theorem justifies our definition.

\begin{theorem} \label{th_wedge}
Assume that $R_1$ and $R_2$ are wedgeable.
Then, the right hand
side of (\ref{eq_wedge}) is
independent of the choice of $c$, $\Omega^\pm$ and depends
linearly on $\Phi$. Moreover,
$R_1\wedge R_2$ defines
a positive closed $(p_1+p_2,p_1+p_2)$-current of mass $1$ with support in
$\supp(R_1)\cap\supp(R_2)$ 
which depends linearly on each $R_i$ and is symmetric with respect
to the variables.
\end{theorem}
\proof
First, observe that the linear dependence of $\Phi$ and of $R_i$ are easily deduced from the
properties of $\Uc_{p_1}$.
Write $\ddc\Phi=\widetilde c (\widetilde\Omega^+-\widetilde\Omega^-)$
with $\widetilde c\geq 0$ and $\widetilde\Omega^\pm$ smooth in
$\Cc_{k-p_1-p_2+1}$. We have
$$c\Omega^+-c\Omega^-=  \widetilde c\widetilde\Omega^+- \widetilde
c\widetilde\Omega^-.$$
Since $\Uc_{p_1}$ is affine on each variable, we have 
$$c\Uc_{p_1}(R_1,R_2\wedge\Omega^+)-c\Uc_{p_1}(R_1,
  R_2\wedge\Omega^-) =  \widetilde c\Uc_{p_1}(R_1, R_2\wedge\widetilde \Omega^+)-\widetilde
c\Uc_{p_1}(R_1,R_2\wedge \widetilde \Omega^-).$$
So, the right hand side of (\ref{eq_wedge}) does not change if we
replace $c$, $\Omega^\pm$ by $\widetilde c$, $\widetilde \Omega^\pm$. 

Let $R_{i,\theta}$ be the currents constructed in Paragraph
\ref{section_topology} starting with the currents $R_i$; they are smooth for $\theta\not=0$.
Lemma \ref{lemma_decreas_sf} implies that $\Uc_{p_1}(R_{1,\theta_1},R_{2,\theta_2}\wedge\Omega^\pm)$ converge to
$\Uc_{p_1}(R_1,R_2\wedge\Omega^\pm)$ when $\theta_i\rightarrow 0$, see also Remarks \ref{rk_hartogs_cv}. It
follows that when $\theta_i\rightarrow 0$ and
$(\theta_1,\theta_2)\not=(0,0)$, the currents
$R_{1,\theta_1}\wedge R_{2,\theta_2}$ converge to $R_1\wedge
R_2$. Hence, $R_1\wedge R_2$ is a positive closed current of mass
1. Since $\supp(R_{i,\theta})\rightarrow\supp(R_i)$, $R_1\wedge R_2$
has support in $\supp(R_1)\cap\supp(R_2)$.
We also have that $R_{1,\theta_1}\wedge R_{2,\theta_2}=R_{2,\theta_2}\wedge
R_{1,\theta_1}$, hence $R_1\wedge R_2=R_2\wedge R_1$. 
\endproof

\begin{lemma} \label{lemma_wedge_regular}
Let $R_i$ and $R_i'$ be currents in $\Cc_{p_i}$. Assume
  that $R_1$ and $R_2$ are wedgeable.
If $R_i'$ is  more diffuse than $R_i$ for $i=1,2$, then
$R_1'\wedge R_2'$ is more diffuse than $R_1\wedge R_2$. 
\end{lemma}
\proof
By Proposition \ref{prop_criteria_wedge}, $R_1'$ and $R_2$ are
wedgeable. 
Theorem \ref{th_wedge} shows that $R_1\wedge R_2$, $R_1'\wedge R_2$,
$R_1\wedge R_2'$ and $R_1'\wedge R_2'$ are well-defined. 
We show that $R_1'\wedge R_2$ is more diffuse
than $R_1\wedge R_2$. In the same way, we will get that $R_1'\wedge R_2'$
is more diffuse than $R_1' \wedge R_2$ which will complete
the proof.

The symbols $U$ and $\Uc$ below denote quasi-potentials and 
super-potentials of mean 0.
By hypothesis, there is a constant $a$ such that
$\Uc_{R_1'}+a\geq \Uc_{R_1}$.
Consider a smooth form $R$ in
$\Cc_{k-p_1-p_2+1}$ and choose $U_R$ smooth. Since $\ddc U_R=R-\omega^{k-p_1-p_2+1}$,
we deduce from (\ref{eq_wedge}) that
\begin{eqnarray*}
\Uc_{R_1'\wedge R_2}(R) & = &  \langle R_1'\wedge R_2,U_R\rangle\\
&=& \langle R_2,\omega^{p_1}\wedge U_R\rangle +
\Uc_{R_1'}(R_2\wedge R)-\Uc_{R_1'}(R_2\wedge\omega^{k-p_1-p_2+1}).
\end{eqnarray*}
The same identity for $R_1\wedge R_2$ and the inequality
$\Uc_{R_1'}+a\geq \Uc_{R_1}$ imply 
$$\Uc_{R_1'\wedge R_2}(R) - \Uc_{R_1\wedge R_2}(R) \geq 
- a -\Uc_{R_1'}(R_2\wedge\omega^{k-p_1-p_2+1})+\Uc_{R_1}(R_2\wedge\omega^{k-p_1-p_2+1}).$$ 
The last expression is finite and independent of $R$. Hence, using the
regularization $R_\theta$ of $R$ for an arbitrary $R$ in
$\Cc_{k-p_1-p_2+1}$, we deduce that $\Uc_{R_1'\wedge
  R_2}-\Uc_{R_1\wedge R_2}$ is bounded
below by a constant. So, $R_1'\wedge R_2$ is more diffuse than
$R_1\wedge R_2$. 
\endproof

The following continuity result shows that the
wedge-product is the right extension to currents of the wedge-product of smooth forms. 

\begin{proposition} \label{prop_cv_wedge_hartogs}
Let $R_1$, $R_2$ be wedgeable currents as
  above and $R_{i,n}$ be currents in $\Cc_{p_i}$ H-converging to $R_i$.
Then, $R_{1,n}$, $R_{2,n}$ are wedgeable
  and $R_{1,n}\wedge R_{2,n}$ H-converge to $R_1\wedge R_2$.
\end{proposition}
\proof
Let $\Uc_{i,n}$ and $\Uc_i$ denote the super-potentials of mean 0 of
$R_{i,n}$ and $R_i$. Let $a_{i,n}$ be constants converging to 0 such
that $\Uc_{i,n}+a_{i,n}\geq \Uc_i$.  Define
$$\epsilon_n:=\Uc_{1,n}(R_2\wedge \omega^{k-p_1-p_2+1}) - \Uc_1(R_2\wedge \omega^{k-p_1-p_2+1}).$$
We have $\epsilon_n\geq -a_{1,n}$.
Since $\Uc_1(R_2\wedge \omega^{k-p_1-p_2+1})$ is finite, Proposition
\ref{prop_hartogs_sqp} implies that $\limsup \epsilon_n\leq 0$. 
So, $\epsilon_n\rightarrow 0$.
Define 
$$K:=\{R_{1,1},R_{1,2},\ldots\} \cup \{R_1\}$$ 
and
$$\delta_n:=\sup_{S\in K} |\Uc_{2,n}(S\wedge \omega^{k-p_1-p_2+1}) - \Uc_2(S\wedge \omega^{k-p_1-p_2+1})|.$$
We first show  that $\delta_n\rightarrow 0$. Since 
$\Uc_{2,n}-\Uc_2\geq -a_{2,n}$, it is enough to prove that $\limsup
\delta_n'\leq 0$ where 
$$\delta_n':=\sup_{S\in K} \Big(\Uc_{2,n}(S\wedge
\omega^{k-p_1-p_2+1}) - 
\Uc_2(S\wedge \omega^{k-p_1-p_2+1})\Big).$$
Because $R_{1,n}\rightarrow R_1$, $K$ is
compact. Since
$\Uc_{1,m}\rightarrow \Uc_1$ pointwise, we have
\begin{eqnarray*}
\lefteqn{\Uc_2(R_{1,m}\wedge \omega^{k-p_1-p_2+1})=\Uc_{1,m}(R_2\wedge
\omega^{k-p_1-p_2+1})}\\
&& \rightarrow \Uc_1(R_2\wedge \omega^{k-p_1-p_2+1})=\Uc_2(R_1\wedge \omega^{k-p_1-p_2+1}).
\end{eqnarray*}
So, $\Uc_2$, restricted to $K$, is continuous. Proposition
\ref{prop_hartogs_sqp} applied to $\Uc_{2|K}+\epsilon$, implies that
$\limsup \delta_n'\leq 0$. Therefore, $\delta_n\rightarrow 0$.

Proposition \ref{prop_criteria_wedge} implies that
$R_{1,n}$, $R_{2,n}$ are wedgeable and  $R_{1,n}$,
$R_2$ are wedgeable. Let $\Uc_n$ $\Uc_n'$ and $\Uc$ denote the super-potentials of mean 0 of
$R_{1,n}\wedge R_{2,n}$, $R_{1,n}\wedge R_2$ and $R_1\wedge R_2$. We obtain as in Lemma
\ref{lemma_wedge_regular} 
for $R$ smooth that $\Uc_n(R)$ and $\Uc_n'(R)$ converge to $\Uc(R)$. Moreover,
$$\Uc'_n(R) -\Uc(R)  \geq  -|a_{1,n}| - |\epsilon_n|$$
and
$$\Uc_n(R) -\Uc_n'(R)  \geq  -|a_{2,n}| -  \delta_n.$$
Hence,
$$\Uc_n(R) \geq \Uc(R)  -|a_{1,n}| - |a_{2,n}|  -|\epsilon_n|-\delta_n$$
for $R$ smooth. Using the approximation of $R$ by $R_\theta$, we
deduce this inequality for arbitrary $R$. The super-potentials
$\Uc_n+|a_{1,n}| + |a_{2,n}|  +|\epsilon_n|+\delta_n$ are larger than
$\Uc$ and converge to
$\Uc$. Hence, the sequence $R_{1,n}\wedge R_{2,n}$ H-converges to $R_1\wedge R_2$.
\endproof

\begin{lemma} \label{lemma_wedge_tau}
Let $R_1$ and $R_2$ be currents in $\Cc_{p_i}$. Then,
for $\tau\in \Aut(\P^k)$ outside some pluripolar set, 
$R_1$ and $\tau_*(R_2)$  are wedgeable. 
Moreover, if $R_1,R_2$ are wedgeable, then 
$R_1\wedge \tau_*(R_2)$ converge to $R_1\wedge R_2$ when 
$\tau\rightarrow\id$ in the
fine topology on $\Aut(\P^k)$, i.e. the coarsest topology for which quasi-psh functions
are continuous.
\end{lemma}
\proof
Let $\Uc_{R_1}$ be a super-potential of $R_1$. 
Recall that $\Uc_{R_1}$ is an affine function which is finite on
smooth forms $R$ in $\Cc_{k-p_1+1}$. On the other hand, 
using an average of $\tau_*(R_2)\wedge \omega^{k-p_1-p_2+1}$ we can
obtain a smooth form $R$ in $\Cc_{k-p_1+1}$. Therefore,
the function $\tau\mapsto \Uc_{R_1}\big(\tau_*(R_2)\wedge
\omega^{k-p_1-p_2+1}\big)$  is not identically $-\infty$.
So, it is a quasi-psh function on
$\Aut(\P^k)$ and is finite outside a pluripolar set. Hence, $R_1$
and $\tau_*(R_2)$ are wedgeable for $\tau$ outside this pluripolar
set. 

Assume now that $R_1$ and $R_2$ are wedgeable.
Let $\Phi$ be a real smooth form
of bidegree $(k-p_1-p_2,k-p_1-p_2)$.  By (\ref{eq_wedge}), $\langle
R_1\wedge\tau_*(R_2),\Phi\rangle$ can be written as a difference
of quasi-psh functions on $\Aut(\P^k)$. Hence, in the fine topology on
$\Aut(\P^k)$,
$\langle R_1\wedge\tau_*(R_2),\Phi\rangle$ converge to $R_1\wedge
R_2$ when $\tau\rightarrow \id$. The lemma follows.
\endproof

In order to define the wedge-product of several currents, we need the
following result.

\begin{lemma} \label{lemma_wedge_associative}
Assume that $R_1$ and $R_2$ are wedgeable and that $R_1\wedge
R_2$ and $R_3$ are wedgeable. Then, $R_2$ $R_3$ are
wedgeable and $R_1$, $R_2\wedge R_3$ are wedgeable. Moreover, we have 
$$(R_1\wedge R_2)\wedge R_3=R_1\wedge
(R_2\wedge R_3).$$
\end{lemma}
\proof
We use the symbols $U$ and $\Uc$ for quasi-potentials and super-potentials of mean 0.
Since $\omega^{p_1}$ is more diffuse than $R_1$, by Lemma
\ref{lemma_wedge_regular}, 
$\omega^{p_1}\wedge R_2$ is  more diffuse than $R_1\wedge
R_2$. Proposition \ref{prop_criteria_wedge} implies that $\omega^{p_1}\wedge R_2$ and
$R_3$ are wedgeable. Hence,
$\Uc_{R_3}(\omega^{k-p_2-p_3+1}\wedge R_2)$ is finite. It follows that
$R_2$ and $R_3$ are wedgeable.

We show that $R_1$ and $R_2\wedge R_3$
are wedgeable. By Proposition \ref{prop_cv_wedge_hartogs} and Remark
\ref{rk_hartogs_cv}, $R_{2,\theta}\wedge R_{3,\theta}$ H-converge to
$R_2\wedge R_3$. Using Lemma \ref{lemma_decreas_sf},
we obtain for $p=p_1+p_2+p_3$
\begin{eqnarray*}
\lefteqn{\Uc_{R_1}(R_2\wedge R_3\wedge \omega^{k-p+1}) =  \lim_{\theta\rightarrow
  0}\Uc_{R_1}(R_{2,\theta}\wedge R_{3,\theta}\wedge \omega^{k-p+1})}\\
& = &\lim_{\theta\rightarrow
  0} \langle U_{R_1}, R_{2,\theta}\wedge
R_{3,\theta}\wedge\omega^{k-p+1}\rangle
 =  \lim_{\theta\rightarrow  0} \langle R_{3,\theta}, U_{R_1}\wedge
R_{2,\theta}\wedge\omega^{k-p+1}\rangle\\
& = &  \lim_{\theta\rightarrow  0} \Uc_{R_{3,\theta}}(R_1\wedge
R_{2,\theta}\wedge \omega^{k-p+1})+\langle \omega^{p_3},U_{R_1}\wedge
R_{2,\theta}\wedge\omega^{k-p+1}\rangle\\
&& -\Uc_{R_3}(R_2\wedge\omega^{k-p_2-p_3+1})\\
& = & \Uc_{R_3}(R_1\wedge R_2\wedge \omega^{k-p+1})+
\Uc_{R_1}(R_2\wedge \omega^{k-p_1-p_2+1})-\Uc_{R_3}(R_2\wedge\omega^{k-p_2-p_3+1}).
\end{eqnarray*}
The last sum is finite. Hence, by Proposition \ref{prop_hyp_wedge}, $R_1$ and $R_2\wedge R_3$
are wedgeable.

We now prove 
the identity in the lemma. Proposition \ref{prop_cv_wedge_hartogs} and Remarks \ref{rk_hartogs_cv} imply that  
$R_{1,\theta}\wedge (R_{2,\theta}\wedge R_{3,\theta})$ converge to $R_1\wedge (R_2\wedge
R_3)$ and
$(R_{1,\theta}\wedge
R_{2,\theta})\wedge R_{3,\theta}$ converge to $(R_1\wedge R_2)\wedge
R_3$. For $\theta\not=0$, since $R_{i,\theta}$ are smooth, we have $(R_{1,\theta}\wedge
R_{2,\theta})\wedge R_{3,\theta}=R_{1,\theta}\wedge
(R_{2,\theta}\wedge R_{3,\theta})$. Letting $\theta\rightarrow 0$
gives the result.
\endproof

\begin{definition}\rm
We say that $R_1,\ldots, R_l$ are {\it wedgeable} if $R_1\wedge \ldots\wedge
R_m$ and $R_{m+1}$ are wedgeable for $m=1,\ldots,l-1$. 
\end{definition}

Lemma \ref{lemma_wedge_associative} implies that this property
and the wedge-product $R_1\wedge \ldots \wedge R_l$ are
symmetric with respect to $R_i$. The wedge-product is a positive closed current of
mass 1. 
Applying inductively Proposition \ref{prop_cv_wedge_hartogs} gives the following result.

\begin{theorem} \label{th_wedge_hartogs}
Let $(R_{i,n})_{n\geq 0}$ be sequences of currents in $\Cc_{p_i}$
H-converging to $R_i$.
Assume that  $R_1,\ldots, R_l$ are wedgeable. 
Then, $R_{1,n},\ldots, R_{l,n}$ are wedgeable and
$R_{1,n}\wedge \ldots\wedge R_{l,n}$ converge  to
$R_1\wedge\ldots\wedge R_l$ in the Hartogs' sense. 
\end{theorem}

\begin{definition} \rm
Let $S$ and $R$ be wedgeable currents in $\Cc_p$ and 
$\Cc_{k-p}$ respectively. Let $a$ be a point in $\P^k$. 
We let $\nu_R(S,a)$ denote the mass of $S\wedge R$ at $a$ and we refer
to it as {\it the Lelong number of $S$ at $a$ relatively to $R$}.
\end{definition}

This notion is related to the directional Lelong numbers of $S$
developed in 
 \cite{Demailly3}.
Consider a classical example.

\begin{example}\rm
Let $S$ be a current in $\Cc_1$ and $u$ be a quasi-potential of
$S$. We have $S=\omega+\ddc u$. If $R$ is the current of integration on
a projective line $D$ which is not contained in $\{u=-\infty\}$,
then $S$ and $[D]$ are wedgeable and
$\nu_{[D]}(S,a)$ exists for every $a$. It is equal to the mass of
$S\wedge [D]=\ddc (u[D])+\omega\wedge [D]$ at $a$, i.e. to the mass of
$\ddc(u[D])$ at $a$. 
\end{example}

We will see in Proposition \ref{prop_good_support} below that if $R$
is locally bounded in a neighbourhood of a hypersurface,
then $\nu_R(S,a)$ exists for every $S$. For the classical case, when $R$
is locally bounded out of $a$, see \cite{Demailly3}.


\subsection{Intersection with currents with regular potentials}

In this paragraph, we will give sufficient conditions for currents to
be wedgeable.

\begin{proposition}
Let $R_i$ be currents in $\Cc_{p_i}$ with $1\leq i\leq l$. Assume that
$R_i$ have 
bounded super-potentials for $1\leq i\leq l-1$. Then, 
$R_1,\ldots,R_l$ are wedgeable.
If moreover $R_l$ has
bounded super-potentials, then $R_1\wedge \ldots \wedge R_l$ has
bounded super-potentials. 
\end{proposition}
\proof
Consider $R'_i:=\omega^{p_i}$. Their super-potentials of mean 0 vanish identically.
It is clear that $R'_1,\ldots,R'_{l-1},R_l$ are wedgeable. 
Since $R_i$ have bounded super-potentials, they are  more
diffuse than $R_i'$. Proposition \ref{prop_criteria_wedge}
implies that $R_1,\ldots,R_l$ are wedgeable.

Assume  that the super-potentials of $R_l$ are bounded. Then, $R_l$ are
more diffuse than $R_l'$. Lemma \ref{lemma_wedge_regular} implies that 
$R_1\wedge \ldots \wedge R_l$ is more diffuse than 
$R_1'\wedge \ldots \wedge R_l'$. It follows that 
$R_1\wedge \ldots \wedge R_l$ has bounded super-potentials.
\endproof

\begin{proposition}
Let $R_i$ be currents in $\Cc_{p_i}$ with $1\leq i\leq l$. Assume that
$R_i$ have continuous super-potentials for $1\leq i\leq l-1$. Then, 
$R_1\wedge \ldots \wedge R_l$ depends continuously on $R_l$. If moreover $R_l$ has
continuous super-potentials, then $R_1\wedge \ldots \wedge R_l$ has
continuous super-potentials. 
\end{proposition}
\proof
We only have to consider the
case where $l=2$. Since $R_1$ has continuous super-potentials, it
follows from (\ref{eq_wedge}) that $R_1\wedge R_2$ depends
continuously on $R_2$. Assume that $R_2$ has also continuous
super-potentials. Let $\Uc_{R_1\wedge R_2}$ and $\Uc_{R_i}$ denote the super-potentials of mean 0 of
$R_1\wedge R_2$ and of $R_i$.
Applying (\ref{eq_wedge}) to a smooth quasi-potential $U_R$ of mean 0
of a smooth form $R$ in $\Cc_{k-p_1-p_2+1}$ gives
\begin{eqnarray*}
\Uc_{R_1\wedge R_2}(R) & = &\langle R_1\wedge R_2,U_R\rangle \\
& = &  
\Uc_{R_2}(\omega^{p_1}\wedge R) +
\Uc_{R_1}(R_2\wedge R)-\Uc_{R_1}(R_2\wedge\omega^{k-p_1-p_2+1}).
\end{eqnarray*}
Since $\Uc_{R_i}$ are continuous and $R_2\wedge R$
 depends continuously on $R$, the last expression can be extended
 continuously to $R$ in $\Cc_{k-p_1-p_2+1}$. Hence, $R_1\wedge R_2$
 has continuous super-potentials.
\endproof

\begin{definition} \label{def_psc}
\rm
A compact subset $K$ of $\P^k$ is {\it $(p+1)$-pseudoconvex} if there is a
current in $\Cc_{k-p}$ with compact support in $\P^k\setminus K$, see
also \cite{FornaessSibony2}. 
\end{definition}

Observe that one can approximate the previous current by smooth
elements of $\Cc_{k-p}$ with compact support in $\P^k\setminus K$. So, there
is a smooth positive closed $(k-p,k-p)$-form $\Theta$ with compact support in
$\P^k\setminus K$. If the $2(k-p)$-dimensional Hausdorff measure of $K$
vanishes, then $K$ is $(p+1)$-pseudoconvex. Indeed, generic projective
planes of dimension $p$ do not intersect $K$. In particular, analytic
sets of pure codimension $p$ are $p$-pseudoconvex.

To explain the terminology, observe that we can assume that $\Theta$
has mass 1 and there is a smooth $(k-p-1,k-p-1)$-form $\Phi$
such that $\ddc \Phi=-\Theta+\omega^{k-p}$.
So, $\ddc \Phi$ is strictly positive on $K$. Adding to $\Phi$ a large
positive closed form allows to assume that $\Phi$ is positive on $\P^k$, compare with
Definition \ref{def_loc_psc} for $X=\P^k$.

\begin{proposition} \label{prop_good_support}
Let $R_i$ be currents in $\Cc_{p_i}$.
Assume that $R_i$ are locally bounded forms on open sets $W_i\subset
\P^k$ such that $\P^k\setminus (W_1\cup W_2)$ is 
$(p_1+p_2)$-pseudoconvex. Then, $R_1$ and $R_2$ are wedgeable.
\end{proposition}
\proof
Let $\Theta$ be a smooth form in $\Cc_{k-p_1-p_2+1}$ with compact
support in $W_1\cup W_2$. Fix open sets $W_i'\Subset W_i$ such that
$\supp(\Theta)\subset W_1'\cup W_2'$.  
Reducing $W_i$ if necessary, we can assume that $R_i$ are bounded on
$W_i$.
Proposition \ref{prop_hyp_wedge} implies that it suffices to show that 
$$\Uc_{p_1}(R_1,R_2\wedge\Theta)\geq
-A(1+\|R_1\|_{\infty, W_1}+\|R_2\|_{\infty, W_2})$$
where $A>0$ is independent of $R_i$. This estimate is uniform on
$R_i$, we can then use a regularization 
and assume that $R_i$ are smooth.

Let $U_i$ denote the Green quasi-potentials of $R_i$ and $m_i$
their means. Lemma
\ref{lemma_levine_c1} implies that $\|U_i\|_{\Cc^1(W_i')}\leq
c(1+\|R_i\|_{\infty, W_i})$ and
$|m_i|\leq c$ for $c>0$ independent of $R_i$. Let
$\chi_i$ be positive smooth functions with compact support in $W_i'$
such that $\chi_1+\chi_2=1$ on $\supp(\Theta)$. We have
\begin{eqnarray*}
\Uc_{p_1}(R_1,R_2\wedge\Theta) & = & \langle U_1,R_2\wedge
\Theta\rangle -m_1 \\
& = & \langle \chi_2 U_1 ,R_2\wedge \Theta\rangle + \langle \chi_1 U_1 ,R_2\wedge
\Theta\rangle -m_1. 
\end{eqnarray*}
Since $\chi_1U_1$ is bounded,
we only have to estimate the first integral. By Stokes' formula, it is
equal to the sum of $\langle
\chi_2U_1,\omega^{p_2}\wedge\Theta\rangle$ which is bounded, and the integral
\begin{eqnarray*}
\langle \chi_2 U_1 , \ddc U_2\wedge \Theta\rangle  
& = & \langle \chi_2 \ddc U_1, U_2 \wedge \Theta\rangle+ \langle
d\chi_2 \wedge\dc U_1, U_2 \wedge \Theta\rangle \\
& &  - \langle \dc \chi_2 \wedge d U_1, U_2 \wedge \Theta\rangle+
\langle U_1\wedge\ddc \chi_2,U_2\wedge \Theta\rangle\\
& = & \langle \chi_2 R_1, U_2 \wedge \Theta\rangle -\langle\chi_2\omega^{p_1},U_2\wedge\Theta\rangle
- \langle d\chi_1\wedge \dc U_1, U_2 \wedge \Theta\rangle\\
& & + \langle \dc \chi_1 \wedge d U_1, U_2 \wedge \Theta\rangle
-\langle U_1\wedge \ddc\chi_1,U_2\wedge\Theta\rangle.
\end{eqnarray*}
We used $d\chi_2=-d\chi_1$ and $\ddc \chi_2=-\ddc\chi_1$ on $\supp(\Theta)$. It is clear that the last
sum is of order at most equal to $1+\|R_1\|_{\infty,
  W_1}+\|R_2\|_{\infty, W_2}$. Indeed, we have $\|U_i\|\leq c$ and each integral is taken on
a domain where we can use the estimates on $\|U_i\|_{\Cc^1(W_i')}$.
\endproof

\begin{remark} \rm
It is enough to assume that $R_i$ are in $\Lc^s_\loc(W_i)$ with $s>2k$. 
\end{remark}

We deduce from Proposition \ref{prop_good_support} and Lemma \ref{lemma_levine_c1} the following results.

\begin{corollary} \label{cor_good_supprt}
Let $R_i$ be currents in $\Cc_{p_i}$. Assume for
  $i=2,\ldots,l$ that
  the intersection of the supports of $R_1,\ldots,R_i$ is 
  $(p_1+\cdots+ p_i)$-pseudoconvex. Then,
  $R_1,\ldots,R_l$ are wedgeable.
\end{corollary}

\begin{corollary}
Let $V_i$, $1\leq i\leq l$, be analytic subsets of pure codimension $p_i$ in
$\P^k$. Assume that their intersection is of pure codimension
$p_1+\cdots+p_l$. Let $I_n$ denote the components of
$V_1\cap\ldots\cap V_l$ and $m_n$ their multiplicities. 
Then, the currents of integration on $V_i$ are
wedgeable and we have
$$[V_1]\wedge\ldots\wedge [V_l]=\sum m_n[I_n].$$
\end{corollary}
\proof
It is clear that $V_1\cap\ldots\cap V_i$ is of pure codimension
$p_1+\cdots+p_i$. Hence, it is 
$(p_1+\cdots+p_i)$-pseudoconvex. Corollary \ref{cor_good_supprt} implies that
$V_1,\ldots,V_l$ are wedgeable  and
$[V_1]\wedge\ldots\wedge [V_l]$ has support in 
$V_1\cap\ldots\cap V_l$ which
is of pure codimension $p_1+\cdots+p_l$. It follows that $[V_1]\wedge\ldots\wedge [V_l]$ is
a combination of $[I_n]$. For the identity in the corollary, by
induction, it is enough to prove it for
 $l=2$. Since $\sum m_n[I_n]$ depends continuously on $V_1$ and
$V_2$, Lemma \ref{lemma_wedge_tau} implies that it is enough to prove the
corollary for $V_1$ and $\tau(V_2)$ where $\tau$ is a generic
automorphism close enough to the identity. So, we can assume that
$m_n=1$ for every $n$. Hence, for a generic point $a$ in $V_1\cap
V_2$, $a$ belongs to the regular parts of $V_1$, $V_2$ and $V_1$,
$V_2$ intersect transversally at $a$. 
It is enough to prove that $[V_1]\wedge [V_2]=[V_1\cap V_2]$ in a
neighbourhood of $a$.
In this
neighbourhood, the $\theta$-regularization $[V_2]_\theta$ of $[V_2]$
is an average of currents of integration on manifolds $\tau(V_2)$ where
$\tau$ is an automorphism close to the identity. Observe that $\tau(V_2)$ is
close to $V_2$ and it intersects $V_1$ transversally on a manifold close
to $V_1\cap V_2$. Hence, $[V_1]\wedge [V_2]_\theta$ is an average 
of $[V_1\cap \tau(V_2)]$. When $\theta$ tends to 0, this mean
converges to $[V_1\cap V_2]$. On the other hand, we have seen in Proposition \ref{prop_cv_wedge_hartogs} that
$[V_1]\wedge [V_2]_\theta$ converge to $[V_1]\wedge [V_2]$. Therefore,
$[V_1]\wedge [V_2]=[V_1\cap V_2]$. The corollary follows.
\endproof


\subsection{Intersection with bidegree (1,1) currents}

Consider now the case where $p_2=\cdots=p_l=1$. For $2\leq i\leq l$,
there is a quasi-psh function $u_i$ on $\P^k$ such that
$$\ddc u_i=R_i- \omega.$$
We have the following lemma.

\begin{lemma} \label{lemma_wedge_equi}
The currents $R_1,\ldots, R_l$ are wedgeable if and only if 
for every $2\leq i\leq l$, $u_i$ is integrable with respect to the trace measure
      of $R_1\wedge\ldots\wedge R_{i-1}$. 
In particular, the last condition is symmetric with respect to
$R_2,\ldots, R_l$. 
\end{lemma}
\proof
It is enough to consider the case $l=2$.
We can assume that $u_2$ is of mean 0. Let $u_{2,\theta}$ be the
quasi-potential of mean 0 of $R_{2,\theta}$.
Since $R_{2,\theta}$ H-converge to
$R_2$, there are constants $a_\theta$ converging to 0 such that
$u_{2,\theta}+a_\theta\geq u_2$ and $u_{2,\theta}$ converge pointwise to $u_2$.
If $\Uc_{R_1}$ is the super-potential of mean 0 of $R_1$, then
$$\Uc_{R_1}(R_2\wedge \omega^{k-p_1}) = \lim_{\theta\rightarrow 0} \Uc_{R_1}(R_{2,\theta}\wedge
\omega^{k-p_1}) =  \lim_{\theta\rightarrow 0} \langle R_1,
u_{2,\theta}\omega^{k-p_1}\rangle =\langle R_1,u_2
\omega^{k-p_1}\rangle.$$
Therefore, $\Uc_{R_1}(R_2\wedge \omega^{k-p_1})$ is finite if and only if   
$u_2$ is integrable with respect to the trace measure
$R_1\wedge\omega^{k-p_1}$ of $R_1$. This
implies the lemma.
\endproof

If $R_2$ has a quasi-potential integrable with respect to $R_1$, it is classical to
define the wedge-product $R_1\wedge R_2$ by
$$R_1\wedge R_2:=\ddc (u_2 R_1)+ \omega\wedge R_1.$$
One defines $R_1\wedge \ldots\wedge R_l$ by
induction.

\begin{lemma} The previous definition coincides with the definition
  given in Paragraph \ref{section_wedge}.
\end{lemma}
\proof
Proposition \ref{prop_cv_wedge_hartogs}
implies that $R_1\wedge
R_{2,\theta}$ converge to $R_1\wedge R_2$ when $\theta\rightarrow
0$. Since $R_{2,\theta}$ is smooth, we have
$$R_1\wedge R_{2,\theta}=R_1\wedge (\ddc u_{2,\theta}+\omega) = \ddc
(u_{2,\theta}R_1)+ \omega\wedge R_1.$$
It is clear that the last expression converge to $\ddc
(u_2R_1)+\omega\wedge R_1$.
\endproof


\section{Complex dynamics in higher dimension}

Super-potentials allow us to construct and to study invariant currents
in complex dynamics. We will give here some applications of this new notion.


\subsection{Pull-back of currents by meromorphic maps} \label{section_pull_back}

The results in this paragraph hold for meromorphic correspondences, in
particular for 
the inverse of a dominant meromorphic map.
For simplicity, we only consider meromorphic maps on $\P^k$. 
Recall that a meromorphic map $f:\P^k\rightarrow\P^k$ is holomorphic
outside an analytic subset $I$ of codimension $\geq 2$ in $\P^k$. Let
$\Gamma$ denote the
closure of the graph of the restriction of $f$ to $\P^k\setminus
I$. This is an irreducible analytic set of dimension $k$ in
$\P^k\times\P^k$. 

Let $\pi_1$ and $\pi_2$ denote the canonical
projections of $\P^k\times \P^k$ on the factors. 
The {\it indeterminacy locus} $I$ of $f$ is the set of points $z\in\P^k$ such that $\dim
\pi_1^{-1}(z)\cap\Gamma\geq 1$. We assume that $f$ is {\it dominant},
that is, $\pi_2(\Gamma)=\P^k$. The {\it second indeterminacy set} of
$f$ is the set $I'$  of points $z\in\P^k$ such that $\dim
\pi_2^{-1}(z)\cap\Gamma\geq 1$. Its codimension is also at least equal
to $2$.
If $A$ is a subset of $\P^k$, define
$$f(A):=\pi_2(\pi_1^{-1}(A)\cap\Gamma)\quad \mbox{and}\quad
f^{-1}(A):=\pi_1(\pi_2^{-1}(A)\cap\Gamma).$$

Define formally for a current $S$ on $\P^k$, not necessarily positive
or closed, the pull-back $f^*(S)$ by
\begin{equation} \label{eq_pullback_def}
f^*(S):=(\pi_1)_*\big(\pi_2^*(S)\wedge [\Gamma]\big)
\end{equation}
where $[\Gamma]$ is the current of integration of $\Gamma$. This 
makes sense if the wedge-product $\pi_2^*(S)\wedge
[\Gamma]$ is well-defined, in particular, when $S$ is smooth. Note that
when $S$ is smooth $f^*(S)$ is an $\Lc^1$ form.
Consider now the case of positive closed currents.
We need some preliminary results.

\begin{lemma} \label{lemma_pull_smooth}
Let $S$ be a current in $\Cc_p$. Assume that the restriction of $S$ to
a neighbourhood of $I'$ is a smooth form. Then, formula
(\ref{eq_pullback_def}) defines 
 a positive closed $(p,p)$-current.
Moreover, the mass $\lambda_p$ of $f^*(S)$ does not depend on $S$.
\end{lemma}
\proof
Since $\pi_{2|\Gamma}$ is a finite
map outside $\pi_2^{-1}(I')\cap\Gamma$, the current $\pi_2^*(S)\wedge
[\Gamma]$ is well-defined there and depends continuously on $S$, see \cite{DinhSibony8}. So, if $S$ is
smooth in a neighbourhood of $I'$,  $\pi_2^*(S)\wedge
[\Gamma]$ is well-defined in a neighbourhood of  $\pi_2^{-1}(I')\cap\Gamma$,
hence, $f^*(S)$ is well-defined and is positive.
Let $U$ be the Green quasi-potential of $S$. This is a negative
form such that $S-\omega^p=\ddc U$. By \cite{DinhSibony8}, $\pi_2^*(U)\wedge
[\Gamma]$ is well-defined outside $\pi_2^{-1}(I')$. Lemma \ref{lemma_levine_c1}
implies that $U$ is continuous in a
neighbourhood of $I'$. Hence, as for $S$, we obtain that 
$f^*(U)$ is well-defined. We 
have $f^*(S)-f^*(\omega^p)=\ddc
f^*(U)$. It follows that $f^*(S)$ and $f^*(\omega^p)$ are
cohomologous. Therefore, they 
 have the same mass.
\endproof

The operator $f_*$ is formally defined by 
\begin{equation}  \label{eq_pushforward_def}
f_*(R):=(\pi_2)_*\big(\pi_1^*(R)\wedge [\Gamma]\big).
\end{equation}

\begin{lemma} \label{lemma_push_smooth}
Let $R$ be a current in $\Cc_{k-p+1}$ which is smooth in a
neighbourhood of $I$. Then, the formula (\ref{eq_pushforward_def})
defines
a positive closed $(k-p+1,k-p+1)$-current.
Moreover, the mass of $f_*(R)$ does not depend on $R$ and is equal to $\lambda_{p-1}$.
\end{lemma}
\proof
We obtain the first part as in Lemma \ref{lemma_pull_smooth}.
Since $f_*(\omega^{k-p+1})$ and $f^*(\omega^{p-1})$ have $\Lc^1$ coefficients,
we also have
$$\|f_*(R)\|=\|f_*(\omega^{k-p+1})\|=\int f_*(\omega^{k-p+1})\wedge
\omega^{p-1} =\int \omega^{k-p+1}\wedge f^*(\omega^{p-1}) =
\lambda_{p-1},$$
which proves the last assertion in the lemma.
\endproof

In order to define $f^*(S)$ we need to define
$\pi_2^*(S)\wedge [\Gamma]$. For this purpose, we can introduce the
notion of super-potential in $\P^k\times\P^k$ and study the
intersection of currents there. We avoid this here.
We call $\lambda_p$ {\it the intermediate
  degree of order $p$} of $f$. 
Denote for simplicity $L:=\lambda_p^{-1}f^*$ and
$\Lambda:=\lambda_{p-1}^{-1}f_*$. 
With this normalization, for $S\in\Cc_p$, $R\in\Cc_{k-p+1}$, the
currents $L(S)$ and $\Lambda(R)$ have mass 1 when they are well-defined.

\begin{lemma} \label{lemma_pull_sp}
Let $S$ be a smooth form in $\Cc_p$ and $\Uc_S$ be a
super-potential of $S$. If $\Uc_{L(\omega^p)}$ is a super-potential of
$L(\omega^p)$, then $\lambda_p^{-1}\lambda_{p-1}\Uc_S\circ\Lambda+\Uc_{L(\omega^p)}$ is equal to a
super-potential of $L(S)$ on the currents $R\in\Cc_{k-p+1}$ which are smooth on
a neighbourhood of $I$.
\end{lemma}
\proof
We can assume that $\Uc_S$ and $\Uc_{L(\omega^p)}$ are of mean 0. Let $\Uc_{L(S)}$ be the
super-potential of mean 0 of $L(S)$.
Let $U_S$ be a smooth
quasi-potential of mean 0 of $S$ and $U_R$ be a quasi-potential of
mean 0 of $R$
which is smooth in a neighbourhood of $I$. 
Since $L(S)$ and $L(\omega^p)$ are smooth outside $I$, the following
computation holds
\begin{eqnarray*}
\Uc_{L(S)}(R) & = & \langle L(S),U_R\rangle=\lambda_p^{-1}\langle
S,f_*(U_R)\rangle\\
& = & \lambda_p^{-1}\langle S-\omega^p,f_*(U_R)\rangle + \lambda_p^{-1}\langle \omega^p,f_*(U_R)\rangle\\
& = & \lambda_p^{-1}\langle \ddc U_S,f_*(U_R)\rangle + \lambda_p^{-1}\langle
f^*(\omega^p),U_R\rangle\\
& = & \lambda_p^{-1}\langle U_S,f_*(\ddc U_R)\rangle + \Uc_{L(\omega^p)}(R)\\
& = & \lambda_p^{-1}\langle U_S,f_*(R)\rangle - \lambda_p^{-1}\langle
U_S,f_*(\omega^{k-p+1})\rangle + \Uc_{L(\omega^p)}(R)\\
& = & \lambda_p^{-1}\lambda_{p-1}\Uc_S(\Lambda(R)) - \lambda_p^{-1}\langle
U_S,f_*(\omega^{k-p+1})\rangle + \Uc_{L(\omega^p)}(R).
\end{eqnarray*}
This implies the result  since the second term in the last line is
independent of $R$. 
\endproof

\begin{definition} \rm
 We say that a current $S$ in $\Cc_p$ 
is {\it $f^*$-admissible} if  there is a current $R_0$ in $\Cc_{k-p+1}$,
which is smooth on a neighbourhood of $I$, 
such that the super-potentials of $S$ are finite at $\Lambda(R_0)$.
\end{definition}

\begin{lemma} \label{lemma_admis_compare}
Let $S$ be an $f^*$-admissible current in $\Cc_p$. Then, 
the super-potentials of $S$ are
finite at $\Lambda(R)$ for every $R$ smooth in $\Cc_{k-p+1}$. In
particular, if $S'\in\Cc_p$ such that $S'\leq cS$ for some positive
constant $c$ or if $S'$ is more diffuse than $S$, then
$S'$ is also $f^*$-admissible.
\end{lemma}
\proof
Since $R$ admits a smooth quasi-potential, we can find a positive
current $U$ such that $\ddc U=R-R_0$ and $U$ smooth in a neighbourhood
of $I$. We have $\ddc\Lambda(U)=\Lambda(R)-\Lambda(R_0)$ and 
by Lemma \ref{lemma_sp_qp}, 
$$\Uc_S(\Lambda(R))\geq
\Uc_S(\Lambda(R_0))-\|\Lambda(U)\|.$$ 
This implies
the first assertion. 
When $S'\leq cS$, as in Proposition \ref{prop_compare_sp}, we obtain $\Uc_{S'}(\Lambda(R_0))>-\infty$.
This also holds when $S'$ is more diffuse than $S$.
Hence, $S'$ is $f^*$-admissible.
\endproof

\begin{lemma} \label{lemma_pull_adm}
Let $S$ be an $f^*$-admissible current in $\Cc_p$. Let $S_n$ be  smooth forms in $\Cc_p$ H-converging
to $S$. Then, $f^*(S_n)$ H-converge 
to a positive closed $(p,p)$-current of mass $\lambda_p$ which does
not depend on the choice of $S_n$.
\end{lemma}
\proof
Let $\Uc_{S_n}$ and $\Uc_S$ be super-potentials of mean 0 of $S_n$ and
$S$. Let $c_n$ be constants converging to 0 such that
$\Uc_{S_n}+c_n\geq \Uc_S$. Recall that $\Uc_{S_n}$ converge pointwise
to $\Uc_S$. If $R$ is smooth in a neighbourhood of $I$, we have
$$\lambda_p^{-1}\lambda_{p-1}\Uc_{S_n}(\Lambda(R))
+\Uc_{L(\omega^p)}(R)
\rightarrow
\lambda_p^{-1}\lambda_{p-1}\Uc_S(\Lambda(R))+\Uc_{L(\omega^p)}(R).$$ 
Lemma \ref{lemma_admis_compare} implies that the last
sum is not identically $-\infty$. 

Lemmas \ref{lemma_pull_sp} and 
\ref{lemma_cv_sp_pointwise} imply that $L(S_n)$ converge  to 
a positive closed current $S'$ of bidegree $(p,p)$. Lemma \ref{lemma_pull_smooth} implies that 
the mass of $S'$ is 
$\lambda_p$. Moreover, $\lambda_p^{-1}\lambda_{p-1}\Uc_{S_n}\circ \Lambda+\Uc_{L(\omega^p)}$ 
(resp. $\lambda_p^{-1}\lambda_{p-1}\Uc_S\circ
\Lambda+\Uc_{L(\omega^p)}$) is 
equal on smooth forms $R$ to some super-potential of $L(S_n)$ (resp. of
$S'$). Denote by $\Uc_{L(S_n)}$ and $\Uc_{S'}$ these super-potentials. We have
$\Uc_{L(S_n)}+\lambda_p^{-1}\lambda_{p-1}c_n\geq \Uc_{S'}$ on smooth forms
$R$. Corollary \ref{cor_def_sp} implies that this inequality holds for every
$R$. Therefore, $L(S_n)\rightarrow S'$ in the Hartogs' sense.

Finally, observe that if $S_n'$ are smooth forms in $\Cc_p$ H-converging 
to $S$, then $S_1$, $S_1'$, $S_2$, $S_2'$, $\ldots$ H-converge also to $S$.
It follows that  $L(S_1)$, $L(S_1')$, $L(S_2)$, $L(S_2')$, $\ldots$
converge. We deduce that the limit $S'$ does not depend on the choice of $S_n$.
We can also obtain the result using that $\Uc_{S'}(R)$ does not depend
on the choice of $S_n$. 
\endproof

\begin{definition}\rm
Let $S$ and $S_n$ be as in Lemma \ref{lemma_pull_adm}. The limit of $f^*(S_n)$
is denoted by $f^*(S)$ and is called {\it the pull-back of $S$ under
  $f$}. We say that $S$ is {\it invariant under $f^*$} or $S$ is
{\it $f^*$-invariant} if $S$ is
$f^*$-admissible and  $f^*(S)=\lambda_pS$.  
\end{definition}

The following result extends Lemmas \ref{lemma_pull_sp} and 
\ref{lemma_pull_adm} when $S$ and $S_n$ are not necessarily smooth.

\begin{proposition} \label{prop_pull_general}
Let $S$ be an $f^*$-admissible current in $\Cc_p$. Let $\Uc_S$,
$\Uc_{L(\omega^p)}$ be super-potentials of $S$ and  $L(\omega^p)$.
Let $S_n$ be currents in $\Cc_p$ H-converging
to $S$. Then, $S_n$ are $f^*$-admissible and $f^*(S_n)$ H-converge
towards $f^*(S)$. Moreover, 
$\lambda_p^{-1}\lambda_{p-1}\Uc_S\circ\Lambda+\Uc_{L(\omega^p)}$ is
equal to a super-potential of $L(S)$ for
$R\in\Cc_{k-p+1}$, smooth in a neighbourhood of $I$. 
\end{proposition}
\proof
If $\Uc_{S_n}$ are super-potentials of mean 0 of $S_n$,
there are constants $c_n$ converging to 0 such that $\Uc_{S_n}+c_n\geq
\Uc_S$.
The last assertion in the proposition was already obtained in the proof of Lemma
\ref{lemma_pull_adm}.
Let $\Uc_{L(S)}$ denote
the super-potential of $L(S)$ which is equal to
$\lambda_p^{-1}\lambda_{p-1}\Uc_S\circ\Lambda+\Uc_{L(\omega^p)}$ for
$R$ smooth in $\Cc_{k-p+1}$.
Let $\Uc_{L(S_n)}$ denote the analogous super-potentials of
$L(S_n)$. Since $\Uc_{S_n}\rightarrow \Uc_S$ pointwise, 
$\Uc_{L(S_n)}\rightarrow \Uc_{L(S)}$ on smooth forms in
    $\Cc_{k-p+1}$.  As in Lemma \ref{lemma_pull_adm}, we obtain
$\Uc_{L(S_n)}+\lambda_p^{-1}\lambda_{p-1}c_n\geq \Uc_{L(S)}$ and this implies that
$L(S_n)$ H-converge towards $L(S)$.
\endproof

In the same way, we have the following.

\begin{definition} \rm
 We say that a current $R$ in $\Cc_{k-p+1}$ 
is {\it $f_*$-admissible} if  the super-potentials of $R$ are finite
at $L(S_0)$
for at least one current $S_0$ in $\Cc_p$
which is smooth in a neighbourhood of $I'$ (or equivalently, for every
$S_0$ smooth in $\Cc_p$).
\end{definition}

If $R'\in\Cc_{k-p+1}$ such that $R'\leq cR$ for some positive
constant $c$ or if $R'$ is more diffuse than $R$, then
$R'$ is also $f_*$-admissible.

\begin{lemma} \label{lemma_push_adm}
Let $R$ be an $f_*$-admissible current in $\Cc_{k-p+1}$. Let $R_n$ be 
 smooth forms in $\Cc_{k-p+1}$ H-converging
to $R$. Then,  $R_n$ are $f_*$-admissible and 
$f_*(R_n)$ H-converge 
to a positive closed $(k-p+1,k-p+1)$-current of mass $\lambda_{p-1}$ which does
not depend on the choice of $R_n$.
\end{lemma}

\begin{definition}\rm
Let $R$ and $R_n$ be as in Lemma \ref{lemma_push_adm}. The limit of $f_*(R_n)$
is denoted by $f_*(R)$ and is called {\it the push-forward of $R$
  under $f$}.
 We say that $R$ is {\it invariant under $f_*$} or  $R$ is {\it
   $f_*$-invariant} if $R$ is
$f_*$-admissible and if $f_*(R)=\lambda_{p-1} R$.  
\end{definition}

\begin{proposition} \label{prop_push_gen}
Let $R$ be an $f_*$-admissible current in $\Cc_{k-p+1}$. 
Let $\Uc_R$, $\Uc_{\Lambda(\omega^{k-p+1})}$ be super-potentials of
$R$ and  $\Lambda(\omega^{k-p+1})$. 
Let $R_n$ be $f_*$-admissible currents in $\Cc_{k-p+1}$ H-converging
to $R$. Then, $f_*(R_n)$ H-converge 
to $f_*(R)$. Moreover, 
$\lambda_p\lambda_{p-1}^{-1}\Uc_R\circ
L+\Uc_{\Lambda(\omega^{k-p+1})}$ is equal to a super-potential of
$\Lambda(R)$ on 
$S\in\Cc_p$, smooth in a neighbourhood of $I'$.
\end{proposition}

Note that if an analytic subset $H$ of pure dimension in $\P^k$, of
a given degree, is 
generic in the Zariski sense, then
$[H]$ is $f^*$- and $f_*$-admissible. One
can check that $f^*[H]$ and $f_*[H]$ depend continuously on $H$.


\subsection{Pull-back by maps with small singularities} \label{section_pullback_reg}

We will give in this paragraph
sufficient conditions, easy to check, in order to define the pull-back and push-forward
operators. 
We need some preliminary results. In what follows, $X$ is a complex manifold
of dimension $k$ and $\omega_X$ is a Hermitian form on $X$.

\begin{definition} \label{def_loc_psc}
\rm
A compact subset $K$ of $X$ is {\it weakly $p$-pseudoconvex} if  there is a positive
smooth $(k-p,k-p)$-form
$\Phi$ on $X$ such that $\ddc\Phi$ is strictly
positive on $K$.
\end{definition}

Note that using a cut-off function, we can assume that $\Phi$ has
compact support in $X$. It follows from the discussion after
Definition \ref{def_psc} that $p$-pseudoconvex sets in $\P^k$ are weakly $p$-pseudoconvex.

\begin{lemma} \label{lemma_psc_haus}
If the $(2k-2p+1)$-dimensional Hausdorff measure of $K$ is zero, then
$K$ is weakly $p$-pseudoconvex.
\end{lemma}
\proof
Consider a point $a$ in $K$. We construct a  positive
smooth $(k-p,k-p)$-form $\Phi_a$  such that $\ddc
\Phi_a$ is positive on $K$ and strictly positive at $a$. Since $K$ is
compact, there is a finite sum $\Phi$ of such forms satisfying
Definition \ref{def_loc_psc}. 
Consider local coordinates $z=(z_1,\ldots,z_k)$ such
that $z=0$ at $a$.  Define $z':=(z_1,\ldots,z_{k-p})$ and $z'':=(z_{k-p+1},\ldots,z_k)$.
The hypothesis on the measure of $K$ allows to choose $z$ so that
$K$ does not intersect the set 
$\{|z'|\leq 1 \mbox{ and } 1-\epsilon\leq |z''|\leq 1\}$
where $\epsilon>0$ is a constant. Let $\Theta$ be a positive
$(k-p,k-p)$-form with compact support in the unit ball $\{|z'|<1\}$ of
$\C^{k-p}$, strictly positive at 0. Let $\varphi$ be a positive
function with compact support in the unit ball of $\C^p$ such that
$\varphi=|z''|^2$ for $|z''|\leq 1-\epsilon$. Let $\pi$ denote the
projection $z\mapsto z'$ and define
$\Psi_a:=\varphi(z'')\pi^*(\Theta)$. It is clear that $\Psi_a$ is
positive with compact support in $X$ and $\ddc \Psi_a\geq 0$ on
$K$. Nevertheless, $\ddc\Psi_a$ is not strictly positive at $0$, but it
does not vanish at 0. Observe that if $\tau$ is a linear automorphism
of $\C^k$ close enough to the identity, then $\tau^*(\Psi_a)$ satisfies
the same properties as $\Psi_a$ does. Taking a finite sum of
$\tau^*(\Psi_a)$ gives a form $\Phi_a$ which is strictly positive at 0.
\endproof

The following result is a version of the Oka's inequality, see \cite{FornaessSibony2}.

\begin{proposition} \label{prop_oka}
Let $K$ be a weakly $p$-pseudoconvex compact subset of $X$. 
Let $T$ be a positive $(p,p)$-current on $X$, not necessarily closed. 
Then, for every negative $(p-1,p-1)$-current $U$ on $X$ with
$\ddc U\geq -T$, we have 
$$\|U\|_X\leq c(1+\|U\|_{X\setminus K})$$
where $c>0$ is a constant independent of $U$.
\end{proposition}
\proof
Since $\|U\|_X=\|U\|_{X\setminus K}+\|U\|_K$, we only 
have to bound the mass of $U$ on $K$. Let $\Phi$ be as in
Definition \ref{def_loc_psc} with compact support. Without loss of generality,
we can assume $\ddc\Phi\geq \omega_X^{k-p+1}$ on $K$. We have for some
positive constant $c'$
\begin{eqnarray*}
\|U\|_K & = & -\int_K U\wedge \omega_X^{k-p+1} \leq -\int_K U\wedge
\ddc\Phi\\
& = &   \int_{X\setminus K} U\wedge \ddc\Phi-\int_X U\wedge \ddc\Phi 
\\
& \leq & c'\|U\|_{X\setminus K}-\int_X \ddc U \wedge \Phi
 \leq  c'\|U\|_{X\setminus K}+\int_X T \wedge \Phi.
\end{eqnarray*}
This implies the result since $T$ is fixed.
\endproof

Let $\widetilde\Sigma'$ denote the analytic subset of the points $x$
in $\Gamma$ such that
$\pi_2$ restricted to $\Gamma$ is not locally finite at $x$.
Define $\Sigma':=\pi_1(\widetilde\Sigma')$. We have
$\widetilde\Sigma'\subset \pi_2^{-1}(I')\cap\Gamma$ and
$\Sigma'\subset f^{-1}(I')$.  
The following proposition gives a sufficient condition in order to define
the pull-back of a $(p,p)$-current, see also Lemma
\ref{lemma_critera_iterate} below. 
The result can be applied to a generic meromorphic map in $\P^k$, see
Proposition \ref{prop_generic_mer_map} below.
Note that the hypothesis is satisfied for $p=1$ and in this case the result is
due to M{\'e}o \cite{Meo}. 

\begin{proposition} \label{prop_critera_pullback}
Assume that $\dim \Sigma'\leq k-p$. Then, every
positive closed $(p,p)$-current $S$ is $f^*$-admissible. Moreover, the
pull-back operator $S\mapsto f^*(S)$ is continuous with respect to the
weak topology on currents. 
\end{proposition}
\proof
Let $S_n$ be smooth forms in $\Cc_p$ converging to
$S$. Let $\Uc_{S_n}$ denote the super-potentials of mean 0 of $S_n$.
It is sufficient to prove that for $R$ smooth in
$\Cc_{k-p+1}$,
$\Uc_{S_n}(\Lambda(R))$ converge to a finite number. 
Propositions \ref{prop_pull_general} and \ref{prop_hartogs_sqp} will imply that 
$S$ is $f^*$-admissible.
The convergence implies also that
the limit does not depend
on the choice of $S_n$, see the last argument in Lemma
\ref{lemma_pull_adm}, and that $f^*$ is continuous.

Let $U_{S_n}$ denote the Green quasi-potentials of $S_n$ which are
smooth negative forms such that $\ddc U_{S_n}\geq -\omega^p$. 
These forms converge in $\Lc^1$ to the Green
quasi-potential $U_S$ of $S$. Hence, the means $c_{S_n}$ of $U_{S_n}$
converge to the mean $c_S$ of $U_S$.
Since $U_{S_n}$ and $R$ are smooth, we have
$$\Uc_{S_n}(\Lambda(R))=\langle U_{S_n},\Lambda(R)\rangle -c_{S_n}= 
\lambda_{p-1}^{-1}\langle f^*(U_{S_n}), R\rangle -c_{S_n}.$$
So, it is enough to prove that $f^*(U_{S_n})$
converge in the sense of currents.

The restriction of $\pi_2$ to
$\Gamma\setminus\widetilde\Sigma'$ is a finite map. Under this hypothesis, it was proved in
\cite{DinhSibony8} that  $\pi_2^*(U_{S_n})\wedge
[\Gamma]$ converge in $\P^k\times\P^k$ outside
$\widetilde\Sigma'$. It follows that
$f^*(U_{S_n})$ converge outside $\Sigma'$. Hence,
the mass of $f^*(U_{S_n})$ outside a small neighbourhood $V$ of
$\Sigma'$ is bounded uniformly on $n$.
By Lemma \ref{lemma_psc_haus}, $\Sigma'$ is weakly $p$-pseudoconvex in
$\P^k$. Hence, since $V$ is small, $\overline
V$ is also $p$-pseudoconvex. Using that $\ddc f^*(U_{S_n})\geq -f^*(\omega^p)$,
Proposition \ref{prop_oka} gives 
$$\|f^*(U_{S_n})\|\leq c(1+\|f^*(U_{S_n})\|_{\P^k\setminus V})$$ 
with $c>0$ independent of $S_n$. Therefore, the mass of $f^*(U_{S_n})$
is bounded uniformly on $n$. We can extract
from $f^*(U_{S_n})$ convergent subsequences. In order to prove the convergence
of $f^*(U_{S_n})$ in $\P^k$, it remains
to check that the limit values $U$ of $f^*(U_{S_n})$ have no mass on  $\Sigma'$.

Let $W$ be a small open set in $\P^k$. Write $f^*(\omega^p)=\ddc \Phi$ with
$\Phi$ negative on $W$. So, $\Phi$ and
$U':=U+\Phi$ are negative currents with $\ddc$ positive. Since 
the currents
$U$, $\Phi$ are of bidimension $(k-p+1,k-p+1)$ and $\dim \Sigma'\leq k-p$, it follows from
a result of Alessandrini-Bassanelli 
\cite[Thm 5.10]{AlessandriniBassanelli} that $\Phi$ and $U'$ have no mass on
$\Sigma'$.  
This implies the result. 
\endproof 

\begin{remark}\rm \label{rk_pullback_U}
Assume that $\dim \Sigma'\leq k-p$.
The previous proof gives a definition of $f^*(U_S)$
which depends continuously on $U_S$. The definition can be extended to
negative currents $U$ such that $\ddc U$ is bounded below by a negative closed
current of bounded mass. We still have that $f^*(U)$ depends
continuously on $U$. 
\end{remark}
 
\begin{proposition} \label{prop_pb_push}
Under the hypothesis of Proposition \ref{prop_critera_pullback}, if
$R$ is a current  in $\Cc_{k-p+1}$ with 
bounded (resp. continuous) super-potentials,
then $R$ is $f_*$-admissible and
$\Lambda(R)$ is a current in $\Cc_{k-p+1}$ with bounded (resp. continuous)
super-potentials. 
\end{proposition}
\proof
Assume that the super-potentials of $R$ are bounded. 
It is clear that $R$ is $f_*$-admissible.
Proposition
\ref{prop_push_gen} implies that $\Lambda(R)$ admits a
super-potential equal to $\lambda_p\lambda_{p-1}^{-1}\Uc_R\circ
L+\Uc_{\Lambda(\omega^{k-p+1})}$ 
on smooth $S\in\Cc_p$. The first term is bounded. By Proposition
\ref{prop_critera_pullback}, it
can be extended to a continuous function on $\Cc_p$ if $R$ has
continuous super-potentials. 
So, it is sufficient to prove that the
super-potential $\Uc_{\Lambda(\omega^{k-p+1})}$ of mean 0 of 
$\Lambda(\omega^{k-p+1})$ is continuous. Let $U_S$
be the Green quasi-potential of $S$ and $c_S$ be its
mean. Recall that
$U_S-c_S\omega^{p-1}$ is a quasi-potential of mean 0 of $S$ and $c_S$
depends continuously on $S$. For $S$
smooth, we have 
$$\Uc_{\Lambda(\omega^{k-p+1})}(S)=\langle
U_S-c_S\omega^{p-1},\Lambda(\omega^{k-p+1})\rangle =
\lambda_{p-1}^{-1} \langle f^*(U_S)-c_S
f^*(\omega^{p-1}),\omega^{k-p+1}\rangle.$$
By Remark \ref{rk_pullback_U}, the left hand side can be extended
continuously to $S$ in $\Cc_p$. So, $\Uc_{\Lambda(\omega^{k-p+1})}$
is continuous.
\endproof

If $g:\P^k\rightarrow \P^k$ is a dominant meromorphic map, 
the composition $g\circ f$ is well-defined on a
Zariski dense open set. We extend it as a meromorphic map by
compactifying the graph. {\it The iterate of order $n$ of $f$} is 
the map $f^n:=f\circ\cdots\circ f$ ($n$ times). 
The inverse of $f^n$ is denoted by $f^{-n}$. It should be distinguished
from $f^{-1}\circ \cdots\circ f^{-1}$.
Define $I_n$, $I_n'$ and $\Sigma_n'$ as above for $f^n$ instead of $f$. 
The following lemma will be useful in our dynamical study. 

\begin{lemma} \label{lemma_critera_iterate}
The following conditions are equivalent
\begin{enumerate}
\item  $\dim\Sigma'\leq k-p$.
\item $\dim f^{-1}(A)\leq k-p$ for every analytic subset $A$ of $\P^k$
with $\dim A \leq k-p$. 
\item $\dim \Sigma_n'\leq k-p$ for every $n\geq 1$.
\end{enumerate}
\end{lemma}
\proof
It is easy to check that the first condition implies the second one
and the third condition implies the first one. 
Suppose the second condition. We prove that the first one is satisfied. If not, we can find an
irreductible analytic subset $A$ of $I'$, of minimal dimension, such that
$\dim \pi_1(\pi_2^{-1}(A)\cap\widetilde\Sigma')>k-p$. The second condition in the lemma implies that $\dim
A>k-p$. Let $\widetilde A$ be an irreducible component of
$\pi_2^{-1}(A)\cap\widetilde \Sigma'$ such that
$A':=\pi_1(\widetilde A)$ has dimension $>k-p$. 
By definition of $\widetilde \Sigma'$, we have
$\dim\widetilde A\geq \dim A+1 \geq k-p+2$.

Choose a dense Zariski
open set $\Omega$ of $\widetilde A$ such that 
$\pi_1:\Omega\rightarrow A'$ and $\pi_2:\Omega\rightarrow
A$ are locally submersions. Denote by $\tau_1$ and $\tau_2$ these maps.
If $H$ is a hypersurface of $A$ then $\widetilde
H:=\tau_2^{-1}(H)$ is a hypersurface of
$\Omega$. It has dimension $\geq k-p+1$. The minimality of $\dim A$
implies that $\dim \tau_1(\widetilde H)\leq
k-p<\dim \widetilde H$. Hence, the fibers of $\tau_1$ are of positive  dimension.
Moreover, $\tau_1(\widetilde H)$ has positive codimension in
$A'$. Therefore, since $\widetilde H$ is a hypersurface in $\widetilde
A$, it should be a union of
components of the fibers of $\tau_1$. This
holds for every $H$. Hence,
the fibers of $\tau_2$, which can be obtained as intersections of such
$\widetilde H$, 
are unions of components of the fibers of $\tau_1$. The
intersection of a fiber of $\tau_1$ and a fiber of $\tau_2$ contains
at most 1 point. We deduce that
$\tau_1$ is locally finite, which is a contradiction.

Now, assume the first two conditions. It remains to check that $\dim
\Sigma_n\leq k-p$ for $n\geq 2$.
Using inductively the second condition we get that 
$f^{-1}\circ\cdots\circ f^{-1}(\Sigma')$ has dimension $\leq k-p$.
Observe that $\Sigma_n'$ is the union of the components of dimension $\geq 1$
in the fibers $f^{-n}(x)$. So, $\Sigma_n'$ is contained in the union of
$f^{-1}\circ\cdots\circ f^{-1}(\Sigma')$.
This gives the result.
\endproof


\subsection{Green super-functions for algebraically stable maps} \label{section_green_stable}

Consider a dominant meromorphic map $f$ on $\P^k$ of algebraic degree
$d\geq 2$ and the associated sets
 $I$, $I'$, $I_n$, $I_n'$, $\Sigma'$, $\Sigma_n'$
as in Paragraphs \ref{section_pull_back} and
\ref{section_pullback_reg}. 
Some results in this paragraph can be easily extended to the case of
correspondences, in particular to $f^{-1}$ instead of $f$.
Let $\lambda_p$ denote the
intermediate degree of order $p$ of $f$ and $\lambda_p(f^n)$ the
intermediate degree of order $p$ of $f^n$. Note that
$\lambda_1(f)=d$. 
We have the following elementary lemma, see \cite{DinhSibony2, DinhSibony3} for a
more general context.

\begin{lemma}
The sequence of intermediate degrees $\lambda_p(f^n)$ is sub-multiplicative,
i.e. $\lambda_p(f^{m+n})\leq \lambda_p(f^m)\lambda_p(f^n)$. 
We also have $\lambda_{p+q}(f^n)\leq \lambda_p(f^n)\lambda_q(f^n)$ and $\lambda_p(f^n)\leq d^{pn}$. 
\end{lemma}
\proof
Observe that $(f^{m+n})^*(\omega^p)$ has no mass on
analytic sets. Let $S_j$ be smooth positive closed
forms of mass $\lambda_p(f^n)$ converging locally uniformly to
$(f^n)^*(\omega^p)$ on a Zariski open set. Then, the currents
$(f^m)^*(S_j)$ are of mass $\lambda_p(f^m)\lambda_p(f^n)$ and converge
to  $(f^{m+n})^*(\omega^p)$ on a Zariski open set. 
If $S$ is a limit of $(f^m)^*(S_j)$ in $\P^k$, 
it is of mass $\lambda_p(f^m)\lambda_p(f^n)$ and it satisfies
$S\geq (f^{m+n})^*(\omega^p)$. Hence,
$\|S\|\geq \|(f^{m+n})^*(\omega^p)\|$.
The first inequality in the lemma follows. 

In the same way, we approximate $(f^n)^*(\omega^p)$ and
$(f^n)^*(\omega^q)$ locally uniformly on a suitable Zariski open
set by smooth forms $S_j$ and $S_j'$. If $S$ is a limit current  of
$S_j\wedge S_j'$ in $\P^k$, it satisfies $S\geq  (f^n)^*(\omega^{p+q})$. This
implies   $\lambda_{p+q}(f^n)\leq \lambda_p(f^n)\lambda_q(f^n)$.
For $p=1$ the first assertion in the lemma implies
$\lambda_1(f^p)\leq d^p$. Applying inductively the second inequality
for $q=1$ gives $\lambda_p(f^n)\leq d^{pn}$.
\endproof

The previous lemma implies that the limit
$$d_p:=\lim_{n\rightarrow\infty} \lambda_p(f^n)^{1/n}=\inf_n\lambda_p(f^n)^{1/n}.$$
exists. It is called {\it the dynamical degree of order $p$} of
$f$. We have
$d_p\leq d^p$ for every $p$. The last dynamical degree $d_k$ is also
called {\it the topological degree} of $f$. It is equal to the number
of points in a generic fiber of $f$ and we have
$\lambda_k(f^n)=d_k^n$. 
In general, $\lambda_p(f^n)$ is the degree of $f^{-n}(H)$ where $H$ is
a generic projective plane of codimension $p$. So,  $\lambda_p(f^n)$
is an integer.
A result by Gromov \cite[Theorem 1.6]{Gromov}
implies that $p\mapsto \log \lambda_p(f^n)$ is concave in $p$. 
It follows that $p\mapsto \log d_p$ is also concave in $p$. 
If $f$ is holomorphic, we have $d_p=\lambda_p=d^p$. If $f$ is not
holomorphic, it is easy to prove that $d_k<d^k$. Indeed, if $a$ is the
intersection of generic hyperplanes $H_1,\ldots, H_k$, then
$f^{-1}(a)\subset f^{-1}(H_1)\cap\ldots\cap f^{-1}(H_k)\setminus I$. 
By B\'ezout's theorem, the last set has cardinal $\leq d^k-1$ since all the hypersurfaces
$f^{-1}(H_i)$ contain $I$. 

\begin{definition} \rm
We say that $f$ is {\it algebraically
  $p$-stable} if 
$\lambda_p(f^n)=\lambda_p^n$
for every $n\geq 1$. 
\end{definition}

For such a map we have $d_p=\lambda_p$.
For $p=1$,  the
algebraic 1-stability coincides with the notion introduced by
Forn\ae ss and the second author \cite{Sibony}, i.e.  no
hypersurface is sent by $f^n$ to $I$, see also
 \cite{Nguyen} and Lemma \ref{lemma_fn_variety} below.

\begin{lemma} \label{lemma_small_sigma}
Assume that  $\dim\Sigma'\leq k-p$. Then,
$f$ is algebraically $p$-stable if and only if $(f^*)^n= (f^n)^*$ on $\Cc_p$. 
\end{lemma}
\proof
Recall that by Proposition \ref{prop_critera_pullback} and 
Lemma \ref{lemma_critera_iterate}, $(f^n)^*$ is well-defined and is
continuous on $\Cc_p$.
If $(f^*)^n=(f^n)^*$ on $\Cc_p$, it is clear that
$$\lambda_p(f^n)=\|(f^n)^*(\omega^p)\|=\|(f^*)^n(\omega^p)\|
=\lambda_p^n.$$
Hence, $f$ is   algebraically $p$-stable.
Conversely, by continuity, it is enough to prove the identity
$(f^*)^n=(f^n)^*$ on smooth forms $S$ in $\Cc_p$. Observe that
$(f^*)^n(S)=(f^n)^*(S)$ on a Zariski dense open set $V$ such that $V$,
$f(V)$, $\ldots$, $f^{n-1}(V)$ do not intersect $I$.
As we observed after the definition (\ref{eq_pullback_def}), since $S$ is smooth,
$(f^n)^*(S)$ has no mass on analytic sets. So, $(f^*)^n(S)\geq
(f^n)^*(S)$. When $f$ is algebraically $p$-stable, $(f^*)^n(S)$
and $(f^n)^*(S)$ have mass $\lambda_p^n$ and $\lambda_p(f^n)$ 
 which are equal. It follows that $(f^*)^n(S)= (f^n)^*(S)$.
\endproof

\begin{lemma} \label{lemma_fn_variety}
Assume that $\dim\Sigma'\leq k-p$. For every analytic subset $A_0$ of
$\P^k$ of dimension $k-p$, define by induction
$A_n:=f(A_{n-1}\setminus I)$, and assume 
that $A_n$ is not contained in $I$
for every $n\geq 0$. 
Then, $f$ is algebraically $l$-stable for $l\leq p$.
\end{lemma}
\proof
It is enough to show that
$(f^*)^n(\omega^l)=(f^n)^*(\omega^l)$. We have seen that
the identity holds outside $A:=I\cup f^{-1}(I)\cup\ldots\cup
(f^{-1})^n(I)$ 
and that $(f^*)^n(\omega^l)\geq
(f^n)^*(\omega^l)$. 
The hypothesis implies that $A$ is of dimension $<k-p$. 
Hence, $(f^*)^n(\omega^l)$ has no mass on $A$ because
$(f^*)^n(\omega^l)$ is of bidimension $(k-l,k-l)$. This
completes the proof.
\endproof

\begin{proposition} \label{prop_small_sing}
If $\dim \Sigma'<k-p$, then $f$ is algebraically $l$-stable for $l\leq
p$. In particular, if $f$ is finite, i.e. $I'=\varnothing$, then $f$ is
algebraically $p$-stable for every $p$.
\end{proposition} 
\proof
When $\dim \Sigma'<k-p$, by Proposition \ref{prop_pb_push} applied to
$l+1$ instead of $p$,
$(f_*)^n(\omega^{k-l})$ is well-defined and has no mass on analytic
sets. We deduce as in Lemma \ref{lemma_fn_variety} that
$(f_*)^n(\omega^{k-l})=(f^n)_*(\omega^{k-l})$ and that $f$ is
algebraically $l$-stable.
\endproof

Let $f$ be a finite map. We have $f^{-n}=f^{-1}\circ\cdots\circ
f^{-1}$, $n$ times, therefore, $I_n=I\cup \ldots\cup
f^{-n+1}(I)$. So, the dimension of $I_n$ 
is independent of $n$.  It is not
difficult to prove that $d_p=d^p$ for $p< k-\dim I$. Indeed, for such
$p$, we have
$f^*(\omega^p)=f^*(\omega)\wedge\ldots\wedge f^*(\omega)$, $p$ times.  
The following proposition implies that
generic maps in $\Mc_d(\P^k)\setminus\Hc_d(\P^k)$ are algebraically $p$-stable.

\begin{proposition} \label{prop_generic_mer_map}
The family of finite meromorphic maps of algebraic degree $d\geq 2$
on $\P^k$ whose dynamical degrees $d_s$ satisfy $d_1<\cdots <d_k$, 
contains a Zariski dense open set of $\Mc_d(\P^k)\setminus\Hc_d(\P^k)$. 
\end{proposition} 
\proof
Denote for simplicity $\Mc:=\Mc_d(\P^k)\setminus \Hc_d(\P^k)$ and
recall that this is an irreducible hypersurface of $\Mc_d(\P^k)$ \cite{GuelfandKapranov}.
We can check easily using the coefficients of $f$ that the set 
$\Mc'$ of maps $f$ in $\Mc$ which are finite and of (maximal) topological degree
$d^k-1$ is a Zariski open set in $\Mc$. 
For such a map, we have $d_{k-1}\leq d^{k-1}<d_k$
and since $p\mapsto \log d_p$ is concave, we obtain $d_1<\cdots<d_k$.
It remains to check that $\Mc'$ is not empty.

Consider the map defined on homogeneous coordinates  by
$$f[z_0:\cdots:z_k]:=[z_0^{d-1}z_1:z_0^{d-1}z_2-z_1^d:\cdots:z_0^{d-1}z_k-z_{k-1}^d:
z_0^{d-1}z_1-z_k^d].$$
The indeterminacy set is the common zero set of the components of
$f$. So, $I$ contains only the point $[1:0:\cdots:0]$. The 
map $f$ is not holomorphic, hence $d_k\leq d^k-1$. On the other hand, if
$t$ is a root of order $d^k-1$ of the unity,
$[1:t:t^d:\cdots:t^{d^{k-1}}]$ is sent by $f$ to $I$. Hence, $d_k=d^k-1$.
We show that $f$ is {\it finite}, i.e. $I'$ is
empty. If not, there is $(a_0,\ldots,a_k)\not=0$ in $\C^{k+1}$ such
that the equations
$$z_0^{d-1}z_1=a_0,\quad z_0^{d-1}z_2-z_1^d=a_1,\quad \ldots \quad,
z_0^{d-1}z_1-z_k^d=a_k$$
define an algebraic set of positive dimension. Consider a sequence of
solutions $z^{(n)}=(z_0^{(n)},\ldots,z_k^{(n)})$ such that $|z^{(n)}|$
tend to infinity and that $z_j^{(n)}/|z^{(n)}|$ converge to some
values $x_j$. We have $|x|=1$ and
$$x_0^{d-1}x_1=0,\quad x_0^{d-1}x_2-x_1^d=0,\quad \ldots \quad,
x_0^{d-1}x_1-x_k^d=0.$$
Hence, $|x_0|=1$ and $x_1=\cdots=x_k=0$. Therefore, we can assume
that $z_0^{(n)}$ tends to infinity and is strictly large than the other $z_j^{(n)}$. 
Extracting a subsequence allows to assume that
for some index  $m\geq 1$, $z_m^{(n)}$ is the largest coordinate
between  $z_1^{(n)},\ldots,z_k^{(n)}$. The equation
$z_0^{d-1}z_m-z_{m-1}^d=a_m$ implies that $z_m^{(n)}\rightarrow
0$. Hence, $z_j^{(n)}\rightarrow 0$ for every $j\geq 1$. 
On the other hand, we deduce from the considered equations that $z_k^d=a_0-a_k$. So,
$a_k=a_0$ and $z_k^{(n)}=0$. Using the given equations and the fact
that $z_j^{(n)}\rightarrow 0$, we obtain inductively that $z_j^{(n)}=0$ for
$j\geq 1$ and then $a_j=0$ for every $j\geq 0$. This is a contradiction. 
\endproof

\begin{theorem} \label{th_green_current}
Let $f:\P^k\rightarrow\P^k$ be an algebraically $p$-stable meromorphic map of dynamical
degrees $d_s$ and $\Sigma'$ be defined as above. 
Assume that $\dim \Sigma'\leq k-p$ and $d_{p-1}<d_p$.
Let $S_n$ be currents in $\Cc_p$ and $\Uc_{S_n}$ be 
super-potentials of $S_n$ such that 
$\|\Uc_{S_n}\|_\infty=o\big(d_{p-1}^{-n}d_p^n\big)$. Then, 
$d_p^{-n}(f^n)^*(S_n)$ H-converge 
to an $f^*$-invariant current $T$ in $\Cc_p$ which does
not depend on $S_n$.
\end{theorem}

We call $T$ {\it the Green $(p,p)$-current} associated to $f$.
Define for simplicity $L:=d_p^{-1}f^*$ and $\Lambda:=d_{p-1}^{-1}f_*$. 
Proposition \ref{prop_small_sing} implies that $f$ is algebraically $(p-1)$-stable. Hence,
$\lambda_{p-1}=d_{p-1}<d_p$. 
We have seen that $L:\Cc_p\rightarrow\Cc_p$ is continuous and
$L^n=d_p^{-n}(f^n)^*$ on $\Cc_p$. It follows that the convex set of $f^*$-invariant
currents $S$ in $\Cc_p$ is not empty. Indeed, it contains all the
limit values of the Ces{\`a}ro means
$${1\over N} \sum_{j=0}^{N-1} L^j(\omega^p).$$
Let $\Cc^b_{k-p+1}$ denote the set of the currents $R$ in
$\Cc_{k-p+1}$ with bounded super-potentials.
By Proposition \ref{prop_pb_push}, the operator
$\Lambda:\Cc^b_{k-p+1}\rightarrow\Cc^b_{k-p+1}$
is well-defined.
Consider a current $S$ in $\Cc_p$, a
super-potential $\Uc_S$ of $S$ and a negative super-potential $\Uc_{L(\omega^p)}$ of
$L(\omega^p)$.

\begin{lemma} \label{lemma_pullback_sqp}
The current $L(S)$ admits a super-potential which is
equal to $d_{p-1} d_p^{-1} \Uc_S\circ \Lambda +\Uc_{L(\omega^p)}$
on $\Cc^b_{k-p+1}$. If $S_0$ is an $f^*$-invariant current in $\Cc_p$, 
then it admits a super-potential $\Uc_{S_0}$
satisfying $\Uc_{S_0}=d_{p-1} d_p^{-1} \Uc_{S_0}\circ \Lambda +\Uc_{L(\omega^p)}$
on $\Cc^b_{k-p+1}$.
\end{lemma}
\proof
We prove the first assertion.
By Proposition \ref{prop_pull_general}, we can assume that $S$ is
smooth.  Moreover, there is 
a super-potential $\Uc_{L(S)}$ of $L(S)$ which is
equal to $d_{p-1} d_p^{-1} \Uc_S\circ \Lambda +\Uc_{L(\omega^p)}$
on smooth forms in $\Cc_{k-p+1}$. Consider a current $R$ in
$\Cc^b_{k-p+1}$ and smooth forms $R_n$ in $\Cc_{k-p+1}$ H-converging  to
$R$.
We have $\Uc_{L(S)}(R_n)\rightarrow
\Uc_{L(S)}(R)$ and $\Uc_{L(\omega^p)}(R_n)\rightarrow
\Uc_{L(\omega^p)}(R)$. 
By Proposition \ref{prop_push_gen}, $\Lambda(R_n)\rightarrow
\Lambda(R)$. 
Since $\Uc_S$ is continuous, we deduce that
$\Uc_S(\Lambda(R_n))\rightarrow \Uc_S(\Lambda(R))$. Therefore, 
 $\Uc_{L(S)}=d_{p-1} d_p^{-1} \Uc_S\circ \Lambda
+\Uc_{L(\omega^p)}$ at $R$.

For the second assertion, if $\Uc$ is a super-potential of $S_0$,
since $L(S_0)=S_0$, the first assertion implies that
$\Uc=d_{p-1} d_p^{-1} \Uc\circ \Lambda +\Uc_{L(\omega^p)}+c$
on $\Cc^b_{k-p+1}$, where $c$ is a constant. The super-potential
$\Uc_{S_0}:=\Uc - cd_p(d_p-d_{p-1})^{-1}$ satisfies the lemma.
We use here the property that $d_p\not=d_{p-1}$. 
\endproof

\noindent
{\bf Proof of Theorem \ref{th_green_current}.}
Replacing $\Uc_{S_n}$ by $\Uc_{S_n}+\|\Uc_{S_n}\|_\infty$ allows to
assume that $\Uc_{S_n}$ are positive.
We apply inductively 
Lemma \ref{lemma_pullback_sqp} for $S=L^j(S_n)$. We obtain that
$L^n(S_n)$ admits a
super-potential $\Uc_{L^n(S_n)}$ satisfying
$$\Uc_{L^n(S_n)}= d_{p-1}^nd_p^{-n} \Uc_{S_n}\circ \Lambda^n +
\sum_{j=0}^{n-1} d_{p-1}^j d_p^{-j} \Uc_{L(\omega^p)}\circ
\Lambda^j$$ 
on $\Cc^b_{k-p+1}$. By hypothesis, the first term converges to 0. 
Since $\Uc_{L(\omega^p)}$ is negative, the second term decreases to
$$\Uc:=\sum_{j=0}^\infty d_{p-1}^j d_p^{-j}
\Uc_{L(\omega^p)}\circ\Lambda^j.$$
Hence, $\Uc_{L^n(S_n)}$ converge pointwise in $\Cc^b_{k-p+1}$ to $\Uc$. 
We show that $\Uc$ is not identically $-\infty$.
 
Let $S_0$ be an $f^*$-invariant current in $\Cc_p$ and $\Uc_{S_0}$ be
a super-potential as in Lemma \ref{lemma_pullback_sqp}. We have 
$$\Uc_{S_0}=d_{p-1}d_p^{-1} \Uc_{S_0}\circ \Lambda +\Uc_{L(\omega^p)}.$$
on $\Cc_{k-p+1}^b$. Iterating this identity gives
$$\Uc_{S_0}=d_{p-1}^nd_p^{-n} \Uc_{S_0}\circ \Lambda^n
+\sum_{j=0}^{n-1} d_{p-1}^j d_p^{-j} \Uc_{L(\omega^p)}\circ
\Lambda^j.$$
Since $\Uc_{S_0}$ is bounded from above and since $d_{p-1}< d_p$,
letting $n\rightarrow\infty$ gives $\Uc\geq \Uc_{S_0}$. So, $\Uc$ is not identically $-\infty$.

We deduce from Propositions \ref{prop_unique_sp} and
\ref{prop_compactness_sqp} 
that $L^n(S_n)$ converge to a
current $T$ which admits a super-potential equal to $\Uc$ on
$\Cc_{k-p+1}^b$. The fact that $\Uc$ does not depend on $S_n$ implies that $T$ is
also independent of $S_n$.  
Because $\Uc_{S_n}$ are positive, the convergence is in the Hartogs' sense.
We have
$$L(T)=L(\lim_{n\rightarrow\infty} L^n(S_n))=\lim_{n\rightarrow\infty} L^{n+1}(S_n)=T.$$
Hence, $T$ is $f^*$-invariant.
\hfill $\square$.

\begin{theorem} \label{th_extremal_green_pull}
Let $f$ be as in Theorem \ref{th_green_current}.
Then, the Green $(p,p)$-current $T$ of $f$ 
is the most diffuse current in $\Cc_p$ which is
$f^*$-invariant. In particular, $T$ is extremal in the convex
set of $f^*$-invariant currents in $\Cc_p$.
\end{theorem}
\proof
We have seen in the proof of Theorem \ref{th_green_current} that
$T$ admits a super-potential $\Uc_T$ which is equal to $\Uc$ on
$\Cc_{k-p+1}^b$. It follows that
$$\Uc_T=d_{p-1}d_p^{-1}\Uc_T\circ\Lambda + \Uc_{L(\omega^p)}$$
on $\Cc_{k-p+1}^b$. It is clear that $\Uc_T$ is the unique
super-potential of $T$ satisfying this identity.
Let $S_0$ and $\Uc_{S_0}$ be as above.
We have seen that
$\Uc_T\geq \Uc_{S_0}$ on
$\Cc_{k-p+1}^b$. By Corollary \ref{cor_def_sp}, this inequality holds on
$\Cc_{k-p+1}$. Hence, $T$ is the most diffuse current in $\Cc_p$ which is
$f^*$-invariant. 

We now prove that $T$ is extremal among $f^*$-invariant currents in $\Cc_p$. Assume $T={1\over 2}
(T_1+T_2)$ with $T_i$ in $\Cc_p$ invariant under $f^*$. 
By Lemma \ref{lemma_pullback_sqp}, the $T_i$ admit
super-potentials $\Uc_{T_i}$ such that
$$\Uc_{T_i}=d_{p-1}d_p^{-1} \Uc_{T_i}\circ\Lambda +\Uc_{L(\omega^p)}$$ 
on $\Cc^b_{k-p+1}$. 
This and the uniqueness of $\Uc_T$ imply that $\Uc_T={1\over 2}(\Uc_{T_1}+\Uc_{T_2})$.
On the other hand, we have $\Uc_T\geq \Uc_{T_i}$. Hence,
$\Uc_T=\Uc_{T_i}$ and $T_i=T$. This completes the proof.
\endproof

\begin{theorem} \label{th_push_green}
Let $f:\P^k\rightarrow\P^k$ be a dominant meromorphic map of dynamical
degrees $d_s$ and $\Sigma'$ be defined as above. Assume that $\dim\Sigma'\leq
k-p$ and that $d_p<d_{p-1}$. 
Let $R_n$ be currents in $\Cc_{k-p+1}$ and $\Uc_{R_n}$ be 
super-potentials of $R_n$ such that 
$\|\Uc_{R_n}\|_\infty=o\big((d_p+\epsilon)^{-n}d_{p-1}^n\big)$ for
some constant $\epsilon>0$. Then, 
$d_{p-1}^{-n}(f^n)_*(R_n)$ H-converge 
to an $f_*$-invariant current $T'$ in $\Cc_{k-p+1}$ which does
not depend on $R_n$ and has continuous super-potentials.
\end{theorem}
\proof
Proposition \ref{prop_small_sing} implies that 
$f$ is algebraically $(p-1)$-stable. Hence, $\lambda_{p-1}=d_{p-1}$.
It follows from 
Proposition \ref{prop_pb_push} that the operator $\Lambda:\Cc_{k-p+1}^b\rightarrow
\Cc^b_{k-p+1}$ is well-defined. 
By Proposition \ref{prop_critera_pullback}, $L:\Cc_p\rightarrow\Cc_p$ is well-defined and
is continuous, but we do not have necessarily that $L^n=d_p^{-n}(f^n)^*$.
Replacing $f$ by an iterate $f^N$ allows to assume that
$\lambda_p<d_{p-1}$ and that
$\|\Uc_{R_n}\|_\infty=o\big(\lambda_p^{-n}d_{p-1}^n\big)$. We can also
assume that $\Uc_{R_n}$ are positive.
Let $\Uc_{\Lambda(\omega^{k-p+1})}$ be a negative super-potential of
$\Lambda(\omega^{k-p+1})$. By Proposition \ref{prop_pb_push},
$\Uc_{\Lambda(\omega^{k-p+1})}$ is continuous.
Proposition
\ref{prop_push_gen} implies that $\Lambda^n(R_n)$ admits a super-potential
which is equal to
$$\lambda_p^{n}d_{p-1}^{-n} \Uc_{R_n}\circ L^n+
\sum_{j=0}^{n-1}\lambda_p^jd_{p-1}^{-j}
\Uc_{\Lambda(\omega^{k-p+1})}\circ L^j$$  
on smooth forms in $\Cc_p$. Letting $n\rightarrow\infty$, the first
term tends to 0, the second term decreases to a continuous function
on $\Cc_p$ since $\Uc_{\Lambda(\omega^{k-p+1})}$ and $L$
are continuous and $\lambda_p <d_{p-1}$. This function does not depend
on $R_n$. We deduce that
$\Lambda^n(R_n)$
converge to a current $T'$ which is independent of $R_n$. 
The convergence is in the Hartogs' sense because $\Uc_{R_n}$ are positive.
Moreover,
$T'$ admits a super-potential $\Uc_{T'}$ such that
$$\Uc_{T'}:=\sum_{j= 0}^\infty\lambda_p^jd_{p-1}^{-j}
\Uc_{\Lambda(\omega^{k-p+1})}\circ L^j$$  
on smooth forms in $\Cc_p$. We have seen that the right hand side
defines a continuous function on $\Cc_p$. Hence, $\Uc_{T'}$ is continuous
and the last identity holds on $\Cc_p$. It follows from the convergence
of $\Lambda^n(R_n)$ that $T'$ is $f_*$-invariant.
\endproof

\begin{theorem} \label{th_extremal_as}
Let $f$ and $T'$ be as in Theorem \ref{th_push_green}. Then, $T'$ is
the only $f_*$-invariant current in $\Cc_{k-p+1}$ which has bounded super-potentials.
Moreover, it is extremal in
the convex set of $f_*$-invariant currents in $\Cc_{k-p+1}$.
\end{theorem}
\proof
Let $R$ be a 
 current in $\Cc_{k-p+1}$ with bounded super-potentials.
Theorem
\ref{th_push_green} implies that $\Lambda^n(R)\rightarrow T'$. So, if $R$ is
$f_*$-invariant, then $R=T'$. This implies the first assertion. We
deduce from this and Proposition \ref{prop_compare_sp} the extremality of $T'$.
\endproof


\subsection{Equidistribution problem for endomorphisms}

Consider a holomorphic map $f:\P^k\rightarrow \P^k$
of algebraic degree $d\geq 2$. Recall that $f^*$ acts continuously on
positive closed currents of any bidegree \cite{Meo,DinhSibony8}, see
also Paragraphs \ref{section_pull_back} and \ref{section_pullback_reg}. It is well-known that
$d^{-n}(f^n)^*(\omega)$ converge to a positive closed
$(1,1)$-current $T$ with H{\"o}lder continuous quasi-potentials. One
deduces from the intersection theory of currents that
$d^{-pn}(f^n)^*(\omega^p)$ converge to $T^p$, see
\cite{Sibony,Fornaess} for the first stages of the theory. 
The current $T^p$ is  {\it the Green current of order $p$} 
and its super-potentials are  {\it the Green super-functions of
  order $p$} of $f$.
In the following result, we give a new construction and new properties
of $T^p$.

\begin{theorem} \label{th_green_sp_hol} 
Let  $f:\P^k\rightarrow \P^k$ be a holomorphic map
of algebraic degree $d\geq 2$. 
Then, the Green super-potentials of $f$ are 
H{\"o}lder continuous. Moreover, $T^p$ is extremal in the convex set
of $f^*$-invariant currents
$S$ in $\Cc_p$.
If $S_n$ are currents in $\Cc_p$ of super-potentials
$\Uc_{S_n}$ such that 
$\|\Uc_{S_n}\|_\infty = o(d^n)$, then $d^{-pn} (f^n)^*(S_n)$ H-converge
to $T^p$. 
\end{theorem}

We will see that the proof also gives that $(f,R)\mapsto \Uc_{T^p}(R)$
is locally H\"older continuous on $\Hc_d(\P^k)\times \Cc_{k-p+1}$.
The following lemma is a special case of \cite[Proposition 2.4]{DinhSibony4}. For
the reader's convenience, we give here the proof.

\begin{lemma} \label{lemma_holder}
Let $K$ be a metric space with finite diameter and  $\Lambda:K\rightarrow K$ be a
Lipschitz map: $\|\Lambda(a)-\Lambda(b)\|\leq A\|a-b\|$ with $A>0$. 
Let $\Uc$ be an $\alpha$-H{\"o}lder continuous function on $K$. Then, 
$\sum_{n\geq 0} d^{-n}\Uc\circ\Lambda^n$ converges pointwise to a
function which is $\beta$-H{\"o}lder continuous on $K$ for every $\beta$
such that $\beta<\alpha$ and $\beta\leq \log d/\log A$. 
\end{lemma}
\proof
Here, $\|a-b\|$ denotes the distance between
two points $a$, $b$ in $K$. 
 Since $K$ has
finite diameter (it is enough to assume that $\Uc$ is bounded), it is sufficient to consider $\|a-b\|\ll 1$.
By hypothesis, there is a constant  $A'>0$ such that
 $|\Uc(a)-\Uc(b)|\leq A'\|a-b\|^\alpha$.
Define $A'':=\|\Uc\|_\infty$. Since $K$ has finite diameter, $A''$ is finite.
If $N$ is an integer, we have
\begin{eqnarray*}
\lefteqn{\big|\sum_{n\geq 0} d^{-n}\Uc\circ\Lambda^n(a)-\sum_{n\geq 0}
d^{-n}\Uc\circ\Lambda^n(b)\big|}\\
& \leq & \sum_{0\leq n\leq N} d^{-n}|\Uc\circ\Lambda^n(a)-
\Uc\circ\Lambda^n(b)| + \sum_{n>N} d^{-n}|\Uc\circ\Lambda^n(a)-
\Uc\circ\Lambda^n(b)|\\
& \leq &  A' \sum_{0\leq n \leq N} d^{-n}\|\Lambda^n(a)-\Lambda^n(b)\|^\alpha
+ 2A''\sum_{n>N} d^{-n} \\
&\lesssim & \|a-b\|^\alpha\sum_{0\leq n \leq N} d^{-n}A^{n\alpha} + d^{-N}.
\end{eqnarray*}
If $A^\alpha\leq d$, the last sum is of order at most equal to
$N\|a-b\|^\alpha+d^{-N}$. For a given $0<\beta<\alpha$, choose
$N\simeq -\beta\log\|a-b\|/\log d$. So, the last expression is
$\lesssim \|a-b\|^\beta$. In this case, the function is
$\beta$-H{\"o}lder continuous for every $0<\beta<\alpha$.
When $A^\alpha>d$, the sum is $\lesssim
d^{-N}A^{N\alpha}\|a-b\|^\alpha+d^{-N}$.
If $N\simeq -\log\|a-b\|/\log A$, the last expression is $\lesssim
\|a-b\|^\beta$ with $\beta:=\log d/\log A$. Therefore, the
function is $\beta$-H{\"o}lder continuous. 
\endproof

Define  $L:=d^{-p} f^*$ and
$\Lambda:=d^{-p+1} f_*$.
Recall that $L:\Cc_p\rightarrow \Cc_p$  and
$\Lambda:\Cc_{k-p+1}\rightarrow\Cc_{k-p+1}$ are well-defined 
and are continuous.

\begin{lemma} \label{lemma_push_lipschitz}
The operator $\Lambda$ is
Lipschitz with respect to the distance $\dist_\alpha$ on $\Cc_{k-p+1}$
for $\alpha>0$.
\end{lemma}
\proof
If $\Phi$ is a $\Cc^\alpha$ test $(p-1,p-1)$-form such that
$\|\Phi\|_{\Cc^\alpha}\leq 1$, it is clear that
$\|f^*(\Phi)\|_{\Cc^\alpha}\leq c_\alpha$ for a constant $c_\alpha>0$
independent of $\Phi$. If $R$ and $R'$ are currents in
$\Cc_{k-p+1}$, we have
$$|\langle \Lambda(R)-\Lambda(R'),\Phi\rangle|=|\langle R-R',
d^{-p+1}f^*(\Phi)\rangle|
\leq c_\alpha\dist_\alpha(R,R').$$
The lemma follows. Observe that the estimates are locally uniform on
$f\in\Hc_d(\P^k)$. 
\endproof

\noindent
{\bf Proof of Theorem \ref{th_green_sp_hol}.}  
Theorems \ref{th_green_current} and \ref{th_extremal_green_pull} imply that
$L^n(S_n)$ H-converge to a
current $T_p$
which does not depend on $S_n$ and is extremal among $f^*$-invariant
currents in $\Cc_p$. For $S_n=\omega^p$ and $\Uc_{S_n}=0$,
the computation in those theorems  shows that  $T_p$ admits a
super-potential $\Uc_{T_p}$ satisfying
$$\Uc_{T_p}=\sum_{j=0}^\infty d^{-j}\Uc_{L(\omega^p)}\circ \Lambda^j$$
on smooth forms in $\Cc_{k-p+1}$.
Since $L(\omega^p)$ is smooth,
$\Uc_{L(\omega^p)}$ is Lipschitz. By Lemmas \ref{lemma_holder} and
\ref{lemma_push_lipschitz}, the later sum defines a H{\"o}lder continuous function on
$\Cc_{k-p+1}$. It follows that the last identity holds everywhere on
$\Cc_{k-p+1}$. So, $T_p$ has H{\"o}lder continuous super-potentials.

Let $T$ denote the first Green current of
$f$. So, $T$ is the limit of
$d^{-n} (f^n)^*(\omega)$  in the Hartogs' sense.
By Theorem \ref{th_wedge_hartogs}, $d^{-pn}(f^n)^*(\omega^p)$
converge to $T^p$. Hence, $T_p=T^p$. 
\hfill $\square$ 

\par

\

 Here is one of our main applications of super-potentials.

\begin{theorem} \label{th_endo_generic}
There is a Zariski dense open set $\Hc_d^*(\P^k)$ in
  $\Hc_d(\P^k)$ such that if $f$ is in
  $\Hc_d^*(\P^k)$, then $d^{-pn}(f^n)^*(S)\rightarrow
T^p$ uniformly on $S\in\Cc_p$. 
In particular, for $f$ in $\Hc^*_d(\P^k)$, $T^p$ is the unique current in 
$\Cc_p$ which is $f^*$-invariant.
\end{theorem}

The open set  $\Hc_d^*(\P^k)$ is given by the following lemma.

\begin{lemma} \label{lemma_f_generic}
There is a Zariski dense open set  $\Hc_d^*(\P^k)$ in
  $\Hc_d(\P^k)$ and an integer $N\geq 1$ such that if
  $f$ is in $\Hc_d^*(\P^k)$ and if $\delta$ denotes
  the maximal multiplicity of $f^N$ at a point in $\P^k$, then
$(20k^2\delta)^{8k} < d^N$. 
\end{lemma}
\proof
Fix an $N$ large enough. Observe that the set $\Hc_d^*(\P^k)$ of $f$ satisfying the
previous inequality is a Zariski open set in $\Hc_d(\P^k)$. We only
have to construct such a map $f$ in order to obtain the density of   $\Hc_d^*(\P^k)$.
Choose a
rational map
$h:\P^1\rightarrow \P^1$ of degree $d$ whose critical points are
simple and have disjoint infinite orbits.  
Observe that the multiplicity of $h^N$ at every point is at most equal to 2.
We construct the map $f$ using an idea of Ueda. 
Let $\sigma_k$ denote the group of permutations of
$\{1,\ldots,k\}$. It acts in a canonical way on $\P^1\times
\cdots\times \P^1$, $k$ times. Using the symetric functions on
$(x_1,\ldots,x_k)\in \P^1\times\cdots\times \P^1$, one shows that
$\P^1\times\cdots\times\P^1$ divided by $\sigma_k$ is isomorphic to $\P^k$. 
Let $\pi:\P^1\times
\cdots\times \P^1\rightarrow \P^k$ denote the canonical map.
If $\widehat f$ is the endomorphism of $ \P^1\times \cdots\times
\P^1$, $k$ times,  defined by
$\widehat f(x_1,\ldots,x_k):=(h(x_1),\ldots,h(x_k))$, then there is a
holomorphic map $f:\P^k\rightarrow\P^k$ of algebraic degree $d$ such that
$f\circ\pi=\pi\circ\widehat f$. We also have $f^N\circ \pi=\pi\circ
\widehat f^N$. 
Consider a point $x$ in $\P^k$ and a
point $\widehat x$ in $\pi^{-1}(x)$. 
The multiplicity of $\widehat f^N$ at $\widehat x$ is at most equal to
$2^k$.
It follows that the multiplicity of $f^N$ at $x$ is at most equal to
$2^kk!$ since  $\pi$ has degree $k!$. Therefore, $f$ satisfies the
desired inequality if $N$ is large enough. 
\endproof

Replacing $f$ by $f^N$, one can assume that $f$ satisfies the lemma for
$N=1$.  Let $\delta$ be the maximal
multiplicity of $f$ at a point in $\P^k$.
We introduce some notations. We call {\it dynamical
  super-potential of $S$} the function $\Vc_S$ defined by
$$\Vc_S:=\Uc_S-\Uc_{T^p}-c_S\quad \mbox{where}\quad
c_S:=\Uc_S(T^{k-p+1})-\Uc_{T^p}(T^{k-p+1})$$
where $\Uc_S$ and $\Uc_{T^p}$ are the super-potentials of mean 0 of $S$
and  $T^p$. 
We also call {\it dynamical Green
quasi-potential of $S$} the form
$$V_S:=U_S-U_{T^p}-(m_S-m_{T^p}+c_S)\omega^{p-1}$$
where $U_S$, $U_{T^p}$ are the Green quasi-potentials of $S$, $T^p$ and
$m_S$, $m_{T^p}$ their means.

\begin{lemma}  \label{lemma_dyn_sp_iterate}
We have $\Vc_S(T^{k-p+1})=0$, $\Vc_S(R)=\langle
  V_S,R\rangle$ for $R$ smooth in $\Cc_{k-p+1}$, and 
$\Vc_{L(S)}=d^{-1}\Vc_S\circ\Lambda$ on $\Cc_{k-p+1}$. Moreover, 
$\Uc_S-\Vc_S$ is bounded by a
constant independent of $S$. 
\end{lemma}
\proof It is clear that  $\Vc_S(T^{k-p+1})=0$. 
Since $T^{k-p+1}$ has bounded super-potentials, $c_S$ is bounded by a
constant independent of $S$.
Hence, since $\Uc_{T^p}$ is bounded, $\Uc_S-\Vc_S$ is bounded
by a constant independent of $S$. For $R$ smooth, we have
$$\langle V_S,R\rangle  =  \big(\langle U_S,R\rangle -m_S\big) -\big(\langle U_{T^p},
R\rangle -m_{T^p}\big)-c_S
 = \Uc_S(R)-\Uc_{T^p}(R)-c_S=\Vc_S(R).$$
It remains to prove that $\Vc_{L(S)}=d^{-1}\Vc_S\circ\Lambda$.
Since $\Lambda(T^{k-p+1})=T^{k-p+1}$, we have
$\Vc_{L(S)}=d^{-1}\Vc_S\circ\Lambda=0$ at $T^{k-p+1}$. Hence, we only
have to show that $\Vc_{L(S)}-d^{-1}\Vc_S\circ\Lambda$ is constant.
By Proposition \ref{prop_pull_general}, we have
$$\Uc_{L(S)}=d^{-1}\Uc_S\circ\Lambda +\Uc_{L(\omega^p)}+\const$$
and since $L(T^p)=T^p$, this implies
$$\Uc_{T^p}=d^{-1}\Uc_{T^p}\circ\Lambda  +\Uc_{L(\omega^p)}+\const.$$
It follows that 
$$\Vc_{L(S)} = d^{-1}\Uc_S\circ\Lambda-d^{-1}\Uc_{T^p}\circ\Lambda+\const.$$
So,  $\Vc_{L(S)}-d^{-1}\Vc_S\circ\Lambda$ is constant.
\endproof

\begin{lemma} \label{lemma_estimate_W}
Let  $W_\epsilon$ be the
$\epsilon$-neighbourhood of the set $P$ of critical values of $f$ and $W^c_\epsilon$ be
the complement of $W_\epsilon$ with $0<\epsilon\ll 1$. There is a constant $c>0$
independent of $\epsilon$ such that for $R$ smooth in $\Cc_{k-p+1}$
and for $0<\epsilon'\ll\epsilon$, we have
$$\|\Lambda(R)_{\epsilon'}-\Lambda(R)\|_{\infty,W^c_\epsilon}\leq
c\|R\|_{\Cc^1}\epsilon^{-5k}\epsilon',$$
where $\Lambda(R)_{\epsilon'}$ is the $\epsilon'$-regularization of
$\Lambda(R)$, see Remark \ref{rk_theta_reg} for the terminology.
\end{lemma}
\proof
Let  $B_\epsilon$ be the ball of radius $\epsilon$ centred at a given point $a$ of $W^c_\epsilon$.
Since $B_\epsilon$ does not intersect $P$, $f$ admits $d$
inverse branches on $B_\epsilon$. More precisely, there are $d$ injective holomorphic maps
$g_i:B_\epsilon\rightarrow\P^k$ such that $f\circ g_i=\id$ on
$B_\epsilon$.
Observe that since $f$ is finite, when the diameter of a ball $B$
tends to 0, the connected components of $f^{-1}(B)$ tend to single
points. So, $g_i(B_\epsilon)$ have small size. Using Cauchy's integral, it is easy to
check that all the derivatives order $n$ of $g_i$ on $B_{\epsilon/2}$ are
 $ \lesssim\epsilon^{-n}$.  On
  $B_{\epsilon}$, we have
$$\Lambda(R)=d^{-p+1}\sum g_i^*(R).$$
For fixed local real coordinates $(x_1,\ldots,x_{2k})$, $R$ is a
combination with smooth coefficients of $dx_{i_1}\wedge \ldots\wedge
dx_{i_{2k-2p+2}}$. Hence, 
the estimate on the derivatives of $g_i$ implies that 
$$\|g_i^*(R)\|_{\Cc^1(B_{\epsilon/2})}\lesssim
\|R\|_{\Cc^1}\epsilon^{-2k+2p-3}\lesssim \|R\|_{\Cc^1}\epsilon^{-5k}.$$
It follows that
$$\|\Lambda(R)\|_{\Cc^1(W^c_{\epsilon/2})}\lesssim
\|R\|_{\Cc^1}\epsilon^{-5k}.$$
Let $\tau$ be an automorphism of $\P^k$ close enough to the
identity. Lemma \ref{lemma_deformation_c1} 
implies that
$$\|\tau_*(\Lambda(R))-\Lambda(R)\|_{\infty, W^c_\epsilon}\lesssim
\|R\|_{\Cc^1}\epsilon^{-5k}\dist(\tau,\id).$$
We then deduce the desired estimate from the definition of $\Lambda(R)_{\epsilon'}$.
\endproof

\begin{lemma} \label{lemma_delta_holder}
The quasi-potentials of $f_*(\omega)$ are $\delta^{-1}$-H{\"o}lder
continuous.
\end{lemma}
\proof
Let $B$ be a small ball in $\P^k$. The inverse image $f^{-1}(B)$ of
$B$ is a union of small open sets. Hence, there is a smooth psh
function $u$ on $f^{-1}(B)$ such that $\omega=\ddc u$ there. Define
the function $v$ on $B$ by
$$v(z):=\sum_{w\in f^{-1}(z)} u(w)$$
where the points in $f^{-1}(z)$ are repeated according to their multiplicity. It is
clear that $v$ is continuous and $\ddc v=f_*(\omega)$. We only have to
show that $v$ is $\delta^{-1}$-H{\"o}lder continuous. Recall that the
multiplicity of $f$ at every point is $\leq \delta$. By Lojasiewicz's
inequality \cite[Lemma 4.3]{DinhSibony9}, we 
can write, for $z,z'$ in $B$, 
$$f^{-1}(z)=\{w_1,\ldots ,w_{d^k}\}\quad \mbox{and}\quad 
f^{-1}(z')=\{w_1',\ldots ,w_{d^k}'\}$$
so that $\dist_\FS(w_i,w_i')\lesssim \dist_\FS(z,z')^{\delta^{-1}}$.
Hence,
$$|v(z)-v(z')|\leq d^k\|u\|_{\Cc^1}\max\dist_\FS(w_i,w_i')\lesssim
\dist_\FS(z,z')^{\delta^{-1}}.$$ 
This implies the lemma.
\endproof

\begin{lemma} \label{lemma_souci_pc}
Let $P$ denote the set of critical values of $f$ as above. If $R$ is smooth,
then $\Vc_S(\Lambda(R))=\langle V_S,\Lambda(R)\rangle_{\P^k\setminus P}.$
\end{lemma}
\proof
Observe that $\Lambda(R)$ is smooth outside $P$.
We will show that $\Uc_S(\Lambda(R))=\langle
U_S,\Lambda(R)\rangle_{\P^k\setminus P}-m_S$. This and
the same identity for $T^p$ imply the result. Since $R\leq
c\omega^{k-p+1}$ for a constant $c>0$, we have 
$$\Lambda(R)\leq cd^{1-p} f_*(\omega^{k-p+1})\leq cd^{1-p} [f_*(\omega)]^{k-p+1}.$$ 
Lemma \ref{lemma_delta_holder} and Proposition 
\ref{prop_levine_estimate} imply that 
$\langle U_{S_\theta},\Lambda(R)\rangle_{\P^k\setminus P}$ converge to 
$\langle U_S,\Lambda(R)\rangle_{\P^k\setminus P}$ when
$\theta\rightarrow 0$. So, it is enough to consider the case where $S$
is smooth. In this case, $U_S$ is smooth. Since $\Lambda(R)$ has no
mass on $P$, we have
$$\langle U_S,\Lambda(R)\rangle_{\P^k\setminus P}-m_S = 
\langle U_S,\Lambda(R)\rangle-m_S=\Uc_S(\Lambda(R)).$$
This completes the proof.
\endproof

\begin{proposition} \label{prop_cap_bound}
For every smooth form $R$ in $\Cc_{k-p+1}$, $d^{-4n/5}\Vc_S(\Lambda^n(R))$ converge to $0$
uniformly on $S$. In particular, we have 
$|\log\capacity(\Lambda^n(R))|=o(d^{4n/5})$. 
\end{proposition}

Fix an integer $n$ large enough and define $\epsilon:=d^{-n}$. In what
follows, the symbols $\lesssim$ and $\gtrsim$ mean inequalities up to
multiplicative constants which are independents of $n$ and $i$.
Observe that we can assume $S$ smooth.
Define
$\epsilon_i:=\epsilon^{(20k^2\delta)^{6ki}}$ for $0\leq i\leq n$. The main point here is
that $\epsilon_i/\epsilon_{i-1}$ has to be small.
Define also by induction $R_0:=R$ and 
$R_i:=\Lambda (R_{i-1})_{\epsilon_i}$ the $\epsilon_i$-regularization
of $\Lambda (R_{i-1})$, see Remark \ref{rk_theta_reg}  for the terminology.
Let $V_i$ be the Green dynamical quasi-potentials of $L^i(S)$. They are
 forms with bounded mass.

\begin{lemma} \label{lemma_estimate_V0}
We have $d^{-i}|\Vc_S(R_i)|\lesssim (-\log\epsilon) d^{-i/4}$.
\end{lemma}
\proof
By Proposition \ref{prop_regularization}, we have
$$\|R_i\|_{\infty}\lesssim
\epsilon_i^{-2k^2-4k}\lesssim \epsilon_i^{-4k^2}.$$ 
Hence, Lemma \ref{prop_estimate_sqp} applied to $K=\P^k$ implies that 
$$d^{-i}|\Vc_S(R_i)|\lesssim  d^{-i}(-\log \epsilon_i) =
d^{-i}(-\log\epsilon) (20k^2\delta)^{6ki}.$$
Lemma \ref{lemma_f_generic} implies the result. Recall that we suppose
$N=1$.
\endproof

\begin{lemma} \label{lemma_estimate_integral}
We have
$\langle V_{n-i},\Lambda(R_{i-1})-R_i\rangle_{\P^k\setminus P} \gtrsim  - \epsilon^i$.
\end{lemma}
\proof
Observe that 
$V_{n-i}':= V_{n-i}-c\omega^{p-1}$ is negative for some 
universal constant $c>0$. Since $\Lambda(R_{i-1})$ and $R_i$ have the
same mass, we also have
$$\langle V_{n-i}',\Lambda(R_{i-1})-R_i\rangle_{\P^k\setminus P} = 
\langle V_{n-i},\Lambda(R_{i-1})-R_i\rangle_{\P^k\setminus P}.$$ 
Proposition \ref{prop_regularization} implies
$$\|R_{i-1}\|_{\Cc^1}\lesssim \epsilon_{i-1}^{-2k^2-4k-1}\lesssim
\epsilon_{i-1}^{-5k^2}.$$ 
Let $W_i$ denote  the $\epsilon_i^{(10k)^{-1}}$-neighbourhood of $P$
and $W^c_i$  its complement. We obtain from Lemma
\ref{lemma_estimate_W} applied to $R:=R_{i-1}$ that 
$$\|\Lambda(R_{i-1})-R_i\|_{\infty,W_i^c}\lesssim
\|R_{i-1}\|_{\Cc^1} \big[\epsilon_i^{(10k)^{-1}}\big]^{-5k}\epsilon_i \lesssim
\epsilon_{i-1}^{-5k^2}\epsilon_i^{1/2}\lesssim \epsilon^i.$$
Since $V_{n-i}'$ has bounded mass, we deduce that
$$|\langle V_{n-i}',\Lambda(R_{i-1})-R_i\rangle_{W_i^c}|\lesssim \epsilon^i.$$
It remains to prove that 
$$\langle V_{n-i}',\Lambda
(R_{i-1})-R_i\rangle_{W_i\setminus P}\geq -\epsilon^i.$$ 
Since $V_{n-i}'$ is negative and $R_i$ is positive, 
it is enough to bound the integral  $\langle V_{n-i}',\Lambda
(R_{i-1})\rangle_{W_i\setminus P}$.
By Proposition \ref{prop_regularization}, we have 
$$R_{i-1}\lesssim \|R_{i-1}\|_\infty\omega^{k-p+1}\lesssim
\epsilon_{i-1}^{-4k^2} \omega^{k-p+1}.$$ 
It follows that
$$\Lambda(R_{i-1})\lesssim \epsilon_{i-1}^{-4k^2}
f_*(\omega^{k-p+1})\lesssim  \epsilon_{i-1}^{-4k^2}
[f_*(\omega)]^{k-p+1}.$$
Lemma \ref{lemma_delta_holder} and Proposition
\ref{prop_levine_estimate} then imply that
$$|\langle V_{n-i}',\Lambda (R_{i-1})\rangle_{W_i\setminus P}|\lesssim 
\epsilon_{i-1}^{-4k^2}\epsilon_i^{(10k)^{-1}(20k^2\delta)^{-k}\delta^{-k}}
\lesssim \epsilon_{i-1}^{-(20k^2\delta)^{2k}}\epsilon_i^{(20k^2\delta)^{-3k}}
\lesssim \epsilon^i.$$
This completes the proof.
\endproof

\noindent
{\bf End of the proof of Proposition \ref{prop_cap_bound}.} 
Since $\Vc_S$ is bounded from above by a constant independent of $S$,
we only have to bound $\Vc_S(\Lambda^n (R))$ from below.  By
Lemmas \ref{lemma_dyn_sp_iterate}
and \ref{lemma_souci_pc}, we have since $R_0=R$ and $R_i$ are smooth
\begin{eqnarray*}
d^{-n}\Vc_S(\Lambda^n(R)) & = & d^{-1}\Vc_{L^{n-1}(S)}(\Lambda (R_0))\\
&=&  d^{-1} \langle
V_{n-1}, \Lambda (R_0)-R_1\rangle_{\P^k\setminus P} +d^{-1} \langle V_{n-1}, R_1\rangle \\
&=&  d^{-1} \langle
V_{n-1}, \Lambda (R_0)-R_1\rangle_{\P^k\setminus P} +d^{-1} \Vc_{L^{n-1}(S)}(R_1) \\
& = &d^{-1} \langle
V_{n-1}, \Lambda (R_0)-R_1\rangle_{\P^k\setminus P} + d^{-2} \Vc_{L^{n-2}(S)}(\Lambda (R_1)). 
\end{eqnarray*}
By induction, we  obtain
\begin{eqnarray*}
d^{-n}\Vc_S(\Lambda^n(R))  & = & d^{-1} \langle
V_{n-1}, \Lambda (R_0)-R_1\rangle_{\P^k\setminus P} +\cdots \\
& & + d^{-n}\langle V_0,\Lambda
(R_{n-1})-R_n\rangle_{\P^k\setminus P} + d^{-n}\Vc_S(R_n).
\end{eqnarray*}
It follows from Lemmas \ref{lemma_estimate_V0} and
\ref{lemma_estimate_integral} that
$$d^{-n}\Vc_S(\Lambda^n(R))\gtrsim -d^{-1}\epsilon-\cdots - 
d^{-n}\epsilon^n - d^{-n/4}(-\log\epsilon)\gtrsim -\epsilon -d^{-n/4}
(-\log\epsilon).$$
Since $\epsilon=d^{-n}$, we get the result. 
\hfill $\square$ 

\bigskip

\noindent
{\bf End of the proof of Theorem \ref{th_endo_generic}.} 
Consider a current $S$ in $\Cc_p$ and a smooth form $R$ in $\Cc_{k-p+1}$.
We want to prove that $L^n(S)$ converge to $T^p$ uniformly on $S$. 
By Propositions \ref{prop_compactness_sqp} and \ref{prop_unique_sp}, it is enough to show that
$\Vc_{L^n(S)}(R)$ converge to 0 uniformly on $S$.
By Lemma \ref{lemma_dyn_sp_iterate}, we have
$$\Vc_{L^n(S)}(R)=d^{-n}\Vc_{S} (\Lambda^n (R)).$$
Proposition \ref{prop_cap_bound} implies the result.
\hfill $\square$ 

\begin{proposition} 
Assume that $f$ is in $\Hc_d^*(\P^k)$. For any $\alpha>0$,
there are constants $c>0$ and $\lambda>1$ such that
if $S$ is in $\Cc_p$ and  $\Phi$ is a test $(k-p,k-p)$-form of class $\Cc^\alpha$, then
$$|\langle d^{-pn}(f^n)^*(S)-T^p,\Phi\rangle|\leq c
\lambda^{-n}\|\Phi\|_{\Cc^\alpha}.$$
In particular, if $\varphi$ is a $\Cc^\alpha$ function such that
$\langle T^k,\varphi\rangle=0$, then 
$$\|d^{-kn}(f^n)_*(\varphi)\|_\infty\leq c\lambda^{-n}\|\varphi\|_{\Cc^\alpha}.$$ 
\end{proposition}
\proof
We prove the fisrt assertion.
Using theory of interpolation as in
Lemma \ref{lemma_compare_dist}, we only have to prove the case
$\alpha=3$. Assume that $\Phi$ has a bounded $\Cc^3$-norm. Multiplying
$\Phi$ by a constant allows to assume that
$\ddc\Phi=R^+-R^-$ where $R^\pm$ are $\Cc^1$ forms
in $\Cc_{k-p+1}$ with bounded $\Cc^1$-norm. A straighforward
computation as above gives
$$\langle
d^{-pn}(f^n)^*(S)-T^p,\Phi\rangle=d^{-n}\Vc_S(\Lambda^n(R^+))-d^{-n}\Vc_S(\Lambda^n(R^-)).$$
The estimates we obtained above give
$$d^{-n}\Vc_S(\Lambda^n(R^\pm))\gtrsim -nd^{-n/4}.$$
On the other hand, since $\Vc_S$ is bounded from above uniformly on
$S$, we have 
$$d^{-n}\Vc_S(\Lambda^n(R^\pm))\lesssim d^{-n}.$$
So, it is enough to take a $\lambda$ smaller than $d^{1/4}$. 

For the second assertion, if $\delta_a$ is the Dirac mass at $a$ then
$$\langle d^{-kn} (f^n)^*(\delta_a),\varphi\rangle = \langle \delta_a,
d^{-kn} (f^n)_*(\varphi)\rangle = d^{-kn}
(f^n)_*(\varphi)(a).$$
Since $\langle T^k,\varphi\rangle=0$, we deduce from the first
assertion that 
$$|d^{-kn} (f^n)_*(\varphi)(a)|\leq c
\lambda^{-n}\|\varphi\|_{\Cc^\alpha}.$$ 
This completes the proof.
\endproof

Note that for $\alpha\leq 2$, we can take $\lambda$ any constant
smaller than $d^{\alpha/2}$ if we replace $\Hc_d^*(\P^k)$ by a
suitable Zariski open set depending on $\lambda$. In dimension 1, Drasin-Okuyama proved in
\cite{DrasinOkuyama} that the second assertion holds for every $f$ if
$a$ is a point on the Julia set, i.e. on the support of the equilibrium
measure.


\subsection{Equidistribution problem for automorphisms} \label{section_eq_aut}

In this paragraph, we consider the class of regular polynomial
automorphisms introduced by the second author in \cite{Sibony}. Let $f$
be a polynomial automorphism of $\C^k$. We extend $f$ to a birational
map on $\P^k$ that we still denote by $f$. Let $I_+$ and $I_-$ be the
indeterminacy sets of $f$ and $f^{-1}$ respectively. 
With the notations  of Paragraph \ref{section_pull_back}, we have $I=I_+$ and
$I'=I_-$. 
They are analytic
subsets of codimension $\geq 2$ in $\P^k$. The map $f$ is said to be {\it regular} if
$I_+\cap I_-=\varnothing$. We summarize here some properties of $f$,
which are deduced from the above assumption \cite{Sibony}.

The indeterminacy sets $I_\pm$ are irreducible and there is an integer
$p$ such that $\dim I_+=k-p-1$ and $\dim I_-=p-1$. They are  contained in the hyperplane
at infinity $L_\infty$. We also have 
$f(L_\infty\setminus I_+)=f(I_-)=I_-$ and $f^{-1}(L_\infty\setminus I_-)=f^{-1}(I_+)=I_+$. 
If $d_\pm$ denote the algebraic degrees of $f^\pm$, then $d_+^p=d_-^{k-p}$.
Denote by $\K_+$
(resp. $\K_-$) 
the set of points $z$ in $\C^k$ such that the forward orbit
$(f^n(z))_{n\geq 0}$ (resp. the backward orbit $(f^{-n}(z))_{n\geq 0}$) is
bounded in $\C^k$. They are closed subsets in $\C^k$ and $\overline
\K_\pm=\K_\pm\cup I_\pm$. Moreover, $I_-$ is attracting for $f$ and
$\P^k\setminus\overline\K_+$ is the attracting basin; 
$I_+$ is attracting for $f^{-1}$ and $\P^k\setminus\overline\K_-$ is
the attracting basin. 

The positive closed $(1,1)$-currents $d_\pm^{- n} (f^{\pm
n})^*(\omega)$ converge to the Green $(1,1)$-currents $T_\pm$
associated to $f^{\pm 1}$. These currents have H{\"o}lder continuous
quasi-potentials out of $I_\pm$ and satisfy $f^*(T_+)=d_+ T_+$ and
$f_*(T_-)=d_- T_-$.  The self-intersections $T_+^p$ and $T_-^{k-p}$ are
positive closed currents of mass 1 with support in the boundaries of
$\overline \K_+$ and $\overline \K_-$ respectively. The probability measure
$\mu:=T_+^p\wedge T_-^{k-p}$ is supported in the boundary of $\K:=\K_+\cap\K_-$. 
The current $T_+^s$, $1\leq s\leq p$, is the {\it Green current of
  order $s$ of $f$} and its super-potentials are called {\it Green
  super-potentials of order $s$ of $f$}.

Let $\Cc_{k-s+1}(W)$ denote the set of currents
in $\Cc_{k-s+1}$ with compact support in an open set $W$.
We assume that $W$ is a neighbourhood of
$I_-$ such that $\overline W\cap I_+=\varnothing$. 
Since $\dim I_-=p-1$,
$\Cc_{k-s+1}(W)$ is not empty for $s\leq p$.
If $\Uc$ is a function on $\Cc_{k-s+1}(W)$, define
$$\|\Uc\|_{\infty,W}:=\sup_{R\in \Cc_{k-s+1}(W)}|\Uc(R)|.$$
In the following
result, we give a new construction of the currents $T_+^s$ and
$T_-^s$.
Note that we cannot apply the results of Paragraph
\ref{section_green_stable} here, since $\Sigma'=L_\infty$.
Indeed, we apply $f^*$ only to currents without mass on $L_\infty$.

\begin{theorem} \label{th_green_current_aut}
Let $f$ and $W$ be as above. Then,
the Green super-potentials of order $s$ of $f$, $1\leq s\leq p$, are
H{\"o}lder continuous on $\Cc_{k-s+1}(W)$.
Let $S_n$ be currents in $\Cc_s$ and $\Uc_{S_n}$
be super-potentials of $S_n$ such that
$\|\Uc_{S_n}\|_{\infty, W}=o(d_+^n)$ for an open set $W$ which contains 
$\overline\K_-$. Then, $d_+^{-sn} (f^n)^*(S_n)\rightarrow T_+^s$. 
\end{theorem}

It is shown in \cite{Sibony} that the current $f^*(\omega^s)$ is of
mass $d_+^s$ for $1\leq s\leq p$, see also Paragraph \ref{section_pull_back}. It follows that
$f_*(\omega^{k-s})$ is of mass $d_+^s$. Define $L_s:=d_+^{-s}f^*$ and $\Lambda_s:=d_+^{-s+1}f_*$.
Assume that the super-potentials of $S$ are finite on
$\Cc_{k-s+1}(W)$. Then, 
$S$ is $f^*$-admissible, because $\Lambda_s(R)$ belongs to
$\Cc_{k-s+1}(W)$ when $\supp(R)$ is close enough to $I_-$.
By Lemma \ref{lemma_pull_adm} and
Proposition \ref{prop_pull_general}, the current $f^*(S)$ is well-defined and is of mass $d_+^s$. 
Consider a
super-potential $\Uc_{L_s(\omega^s)}$ of $L_s(\omega^s)$.
Since $L_s(\omega^s)$ is smooth on $W$, it is easy to check that
$\Uc_{L_s(\omega^s)}$ is Lipschitz on $\Cc_{k-s+1}(W)$.
We first prove the following result.

\begin{proposition} \label{prop_green_current_aut}
Let $S_n$ be currents in $\Cc_s$ and $\Uc_{S_n}$ be super-potentials
of $S_n$ with $\|\Uc_{S_n}\|_{\infty,W}=o(d_+^n)$. If $S$ is a
limit value of $d_+^{-sn} (f^n)^*(S_n)$, then $S$ admits
a super-potential which is equal on $\Cc_{k-s+1}(\P^k\setminus\overline \K_+)$ to $\sum_{n\geq 0} d_+^{-n}
\Uc_{L_s(\omega^s)}\circ\Lambda_s^n$.
Moreover, this equality holds on $\Cc_{k-s+1}(\P^k\setminus I_+)$ when
$W$ contains $\overline \K_-$. 
\end{proposition}
\proof
Reducing $W$ allows to assume that $f(W)\Subset W$. If $W$ contains
$\overline\K_-$, we can keep this property. Fix an open set $W_0$
relatively compact in
$\P^k\setminus\overline \K_+$ which contains $I_-$. If $W$ contains $\overline \K_-$, we
can take $W_0$ relatively compact in $\P^k\setminus I_+$. Observe that $f^{-m}(W)$ contains
$W_0$ for $m$ large enough. So, replacing $S_n$ by
$d_+^{-sm}(f^m)^*(S_{n+m})$ and $W$ by some open set of $f^{-m}(W)$ allows to assume
that $W_0\Subset W$. 

By Proposition \ref{prop_pull_general}, there is a super-potential  of $L_s(S_n)$ 
which is equal on $\Cc_{k-s+1}(W)$ to
$d_+^{-1}\Uc_{S_n}\circ\Lambda_s+\Uc_{L_s(\omega^s)}$.
We apply again this proposition to $L_s(S_n)$. There is a super-potential  of $L_s^2(S_n)$ 
which is equal on $\Cc_{k-s+1}(W)$ to
$$d_+^{-2}\Uc_{S_n}\circ\Lambda_s^2+\Uc_{L_s(\omega^s)}+d_+^{-1}\Uc_{L_s(\omega^s)}\circ\Lambda_s.$$
By induction, $L_s^n(S_n)$ admits a super-potential $\Uc_{L_s^n(S_n)}$
which is equal to 
$$d_+^{-n}\Uc_{S_n}\circ\Lambda_s^n+\Uc_{L_s(\omega^s)}+
d_+^{-1}\Uc_{L_s(\omega^s)}\circ\Lambda_s+\cdots+d_+^{-n+1}
\Uc_{L_s(\omega^s)}\circ\Lambda_s^{n-1}$$
on $\Cc_{k-s+1}(W)$. By hypothesis, the first term tends to 0.
Hence, $\Uc_{L_s^n(S_n)}$ converge to $\sum_{n\geq 0} d_+^{-n}
\Uc_{L_s(\omega^s)}\circ\Lambda_s^n$ on $\Cc_{k-s+1}(W)$. This sum
converges since $\Uc_{L_s(\omega^s)}$ is Lipschitz on $\Cc_{k-s+1}(W)$.

By Proposition \ref{prop_compactness_sqp}, it remains
to show that $\Uc_{L_s^n(S_n)}$ are bounded from above
uniformly on $n$. For this purpose, it is enough to show that the
means  $\Uc_{L_s^n(S_n)}(\omega^{k-s+1})$ of  $\Uc_{L_s^n(S_n)}$  are bounded from above
uniformly on $n$.
If $R_0$ is a smooth form in
$\Cc_{k-s+1}(W_0)$.
We have
$$\Uc_{L_s^n(S_n)}(R_0)=d_+^{-n}\Uc_{S_n}(\Lambda_s^n(R_0))+\Uc_{L_s(\omega^s)}(R_0)+\cdots+d_+^{-n+1}
\Uc_{L_s(\omega^s)}(\Lambda_s^{n-1}(R_0)).$$
This sum is bounded from above.
On the other hand, $R_0$ admits a positive quasi-potential since it is
smooth. Proposition
\ref{lemma_sp_qp} implies the result.
\endproof

\noindent
{\bf End of the proof of Theorem \ref{th_green_current_aut}.}
Since $W$ contains $\overline\K_-$, by
Proposition \ref{prop_green_current_aut}, any cluster point of
$L_s^n(S_n)$ has a super-potential equal to  
$\sum d_+^{-n} \Uc_{L_s(\omega^s)}\circ \Lambda_s^n$ on
$\Cc_{k-s+1}(\P^k\setminus I_+)$.
Lemma \ref{prop_unique_sp} implies that 
there is only one cluster point for the sequence $L_s^n(S_n)$, hence
$L_s^n(S_n)$ converge to
a current $T_s$. This current does not depend on $S_n$ since it admits
a super-potential independent of $S_n$. For
$S_n=\omega^s$, we obtain that $T_s$ is the Green current of order $s$
of $f$. It admits a super-potential $\Uc_{T_s}$ equal to  $\sum_{n\geq 0} d_+^{-n}
\Uc_{L_s(\omega^s)}\circ\Lambda_s^n$ on $\Cc_{k-s+1}(\P^k\setminus I_+)$. Lemma \ref{lemma_holder}
implies that this function is H{\"o}lder continuous on  $\Cc_{k-s+1}(W)$.

Let $T_+:=T_1$. We want next to prove that $T_s=T_+^s$. For this
purpose, it is sufficient to show that $T_s$ and $T_l$ are wedgeable
 and  $T_s\wedge T_l=T_{s+l}$ when $s+l\leq p$. Since $s+l\leq p$,
 there is a smooth form $\Omega\in\Cc_{k-s-l+1}$ with compact support
 in $\P^k\setminus I_+$. Hence, $\Omega\wedge T_l$ has compact support
 in $\P^k\setminus I_+$ and
the super-potentials of $T_s$ are finite at $\Omega\wedge T_l$. 
It follows that $T_s$ and $T_l$ are wedgeable.

The computation in Proposition \ref{prop_green_current_aut} implies that $L_s^n(\omega^s)$
admits a super-potential $\Uc_{L_s^n(\omega^s)}$ 
which is equal to  $\sum_{i=0}^n d_+^{-i}
\Uc_{L_s(\omega^s)}\circ\Lambda_s^i$ on $\Cc_{k-s+1}(\P^k\setminus I_+)$.
Fix a real smooth test form $\Phi$  of bidegree $(k-s-l,k-s-l)$
with compact
support in $\P^k\setminus I_+$. As in Proposition
\ref{prop_unique_sp}, write
$\ddc\Phi=c(\Omega^+-\Omega^-)$ with $c>0$ and
$\Omega^\pm$ in $\Cc_{k-s-l+1}(\P^k\setminus
I_+)$. The sequence $\Omega^\pm\wedge L_l^n(\omega^l)$ converges to
$\Omega^\pm\wedge T_l$. Since these currents have supports in a
fixed compact subset of $\P^k\setminus I_+$, the values of $\Uc_{L_s^n(\omega^s)}$ at  
$\Omega^\pm\wedge L_l^n(\omega^l)$ converge to the value of
$\Uc_{T_s}$ at  $\Omega^\pm\wedge T_l$. The formula (\ref{eq_wedge}) implies
that $L_s^n(\omega^s)\wedge L_l^n(\omega^l)$ converge to
$T_s\wedge T_l$. 
On the other hand, $L_{s+l}^n(\omega^{s+l})$ and
$L_s^n(\omega^s)\wedge L_l^n(\omega^l)$ are smooth forms which are
equal outside $I_+$. They have no mass on $I_+$ because $\dim
I_+<k-s-l$. Hence, these currents are equal.
Therefore, letting $n\rightarrow\infty$ gives 
$T_{s+l}=T_s\wedge T_l$ and in particular $T_s=T_+^s$. 
\hfill $\square$

\begin{theorem}
The Green current $T_+^s$ is the most diffuse $f^*$-invariant current in $\Cc_s$.
In particular, it is extremal in the convex set
of $f^*$-invariant currents in $\Cc_s$.
\end{theorem}
\proof
It follows from the convergence in Theorem \ref{th_green_current_aut} 
that $T_+^s$ is $f^*$-invariant. Let $T$ be an $f^*$-invariant current
in $\Cc_s$ and $\Uc_T$ be a super-potential of $T$. Proposition \ref{prop_pull_general}
implies that $L_s(T)$ admits a super-potential $\Uc$ which is equal
to $d_+^{-1}\Uc_T\circ\Lambda_s +\Uc_{L_s(\omega^s)}$ on
$R$ smooth in $\Cc_{k-s+1}$. Since $L_s(T)=T$,
there is a constant $c$ such that $\Uc=\Uc_T+c$. Subtracting from
$\Uc_T$ an appropriate constant gives another super-potential that we
still denote by $\Uc_T$, such that
$$\Uc_T=d_+^{-1}\Uc_T\circ\Lambda_s +\Uc_{L_s(\omega^s)}$$
on $R$ in $\Cc_{k-s+1}$ which is smooth in a neighbourhood of $I_+$.
The condition on $R$ is invariant under
$\Lambda$. So, iterating the above identity gives
$$\Uc_T=d_+^{-n} \Uc_T\circ\Lambda^n + \sum_{i=0}^{n-1} d_+^{-i}
\Uc_{L_s(\omega^s)}\circ\Lambda_s^i.$$
Since $\Uc_T$ is bounded from above, letting $n\rightarrow\infty$, we
obtain
$$\Uc_T\leq \sum_{i=0}^\infty d_+^{-i}
\Uc_{L_s(\omega^s)}\circ\Lambda_s^i=\Uc_{T_+^s}.$$
This identity holds on smooth forms $R$ in $\Cc_{k-s+1}$. Hence,
$T_+^s$ is more diffuse than $T$.

Now, we prove that $T_+^s$ is extremal among $f^*$-invariant
currents. Assume that $T_+^s={1\over 2}(T+T')$ with $T$ and $T'$ in
$\Cc_s$ invariant by $f^*$. Let $\Uc_T$ be as above. Let $\Uc_{T'}$ be the
analogous super-potential of $T'$. It is the unique super-potential
which satisfies
$$\Uc_{T'}=d_+^{-1}\Uc_{T'}\circ\Lambda_s +\Uc_{L_s(\omega^s)}$$
on smooth forms in $\Cc_{k-s+1}$. 
Observe that ${1\over 2}(\Uc_T+\Uc_{T'})$ is a super-potential of
$T_+^s$ satisfying the same property. It follows that 
$${1\over 2}(\Uc_T+\Uc_{T'})=\Uc_{T_+^s}.$$
We deduce from the inequalities $\Uc_T\leq \Uc_{T_+^s}$ and
$\Uc_{T'}\leq \Uc_{T_+^s}$
that 
$\Uc_T$, $\Uc_{T'}$ are equal to $\Uc_{T_+^s}$.
Hence, $T=T'=T_+^s$.
This implies the result.
\endproof

In the case of bidegree $(p,p)$, we have the following stronger result
which is another main application of the super-potentials. It
was proved by Forn\ae ss and the second author in the case of
dimension $2$ \cite{FornaessSibony1}. 

\begin{theorem} \label{th_aut_unique}
The current $T_+^p$ is the unique
positive closed current of bidegree $(p,p)$ of
mass $1$ supported in $\overline\K_+$.
 The current $T_-^{k-p}$ is the unique
positive closed current of bidegree  $(k-p,k-p)$ of
mass $1$ supported in $\overline \K_-$.
\end{theorem}

In what follows, we only consider currents $S$ in $\Cc_p$ with support in
$\overline \K_+$. By Proposition \ref{prop_estimate_sqp}, their super-potentials of mean 0 are bounded on
$\Cc_{k-p+1}(W)$ uniformly on $S$ when $W\Subset \P^k\setminus
\overline \K_+$. In particular, they are bounded at
the current  $R_\infty:=(\deg I_-)^{-1}[I_-]$.
We call {\it the dynamical super-potential of $S$} the function $\Vc_S$ defined by
$$\Vc_S:=\Uc_S-\Uc_{T_+^p}-c_S\quad \mbox{where}\quad
c_S:=\Uc_S(R_\infty)-\Uc_{T_+^p}(R_\infty)$$
where $\Uc_S$, $\Uc_{T_+^p}$ are the super-potentials of mean 0 of $S$
and $T_+^p$.
We also call {\it the dynamical Green
quasi-potential of $S$} the form
$$V_S:=U_S-U_{T_+^p}-(m_S-m_{T_+^p}+c_S)\omega^{p-1}$$
where $U_S$, $U_{T_+^p}$ are the Green quasi-potentials of $S$, $T_+^p$ and
$m_S$, $m_{T_+^p}$ are their means. Denote for simplicity $L:=L_p$ and $\Lambda:=\Lambda_p$.

\begin{lemma}  \label{lemma_dyn_sp_iterate_aut}
Let $W\Subset \P^k\setminus I_+$ be an open set. Then,
 $\Vc_S(R_\infty)=0$, $\Vc_S(R)=\langle
  V_S,R\rangle$ for $R$ smooth in $\Cc_{k-p+1}(W)$, and 
$\Vc_{L(S)}=d_+^{-1}\Vc_S\circ\Lambda$ on $\Cc_{k-p+1}(W)$. Moreover, 
$\Uc_S-\Vc_S$ is bounded on  $\Cc_{k-p+1}(W)$ by a
constant independent of $S$. 
\end{lemma}
\proof It is clear that  $\Vc_S(R_\infty)=0$. 
Recall that $m_S$, $m_{T_+^p}$ and $c_S$ are bounded.
Since $\Uc_{T_+^p}$ is continuous on
$\Cc_{k-p+1}(W)$,  $\Uc_S-\Vc_S$ is 
bounded on $\Cc_{k-p+1}(W)$ by a constant independent
of $S$. 
We also have for $R$ smooth in $\Cc_{k-p+1}(W)$
$$\langle V_S,R\rangle  =  \big(\langle U_S,R\rangle -m_S\big) -\big(\langle U_{T_+^p},
R\rangle -m_{T_+^p}\big)-c_S
 = \Uc_S(R)-\Uc_{T_+^p}(R)-c_S=\Vc_S(R).$$

It remains to prove that $\Vc_{L(S)}=d_+^{-1}\Vc_S\circ\Lambda$ on $\Cc_{k-p+1}(W)$.
Observe that since $I_-$ is irreducible, $\Lambda(R_\infty)=R_\infty$.
We deduce that 
$\Vc_{L(S)}=d_+^{-1}\Vc_S\circ\Lambda=0$ at $R_\infty$. Hence, we only
have to show that $\Vc_{L(S)}-d_+^{-1}\Vc_S\circ\Lambda$ is constant.
By Proposition \ref{prop_pull_general}, see also Proposition \ref{prop_green_current_aut}, we have
$$\Uc_{L(S)}=d_+^{-1}\Uc_S\circ\Lambda +\Uc_{L(\omega^p)}+\const$$
and since $L(T_+^p)=T_+^p$, this implies
$$\Uc_{T_+^p}=d_+^{-1}\Uc_{T_+^p}\circ\Lambda  +\Uc_{L(\omega^p)}+\const.$$
It follows that 
$$\Vc_{L(S)} = d_+^{-1}\Uc_S\circ\Lambda-d_+^{-1}\Uc_{T_+^p}\circ\Lambda+\const.$$
It is clear that  $\Vc_{L(S)}-d_+^{-1}\Vc_S\circ\Lambda$ is constant.
\endproof

\noindent
{\bf Proof of Theorem \ref{th_aut_unique}.} 
Consider a current $S$ in $\Cc_p(\P^k)$ with support in $\overline\K_+$.
Define $S_n:=d_+^{pn}(f^n)_*(S)$ on $\C^k$. These currents are positive
closed with support in $\overline \K_+$. Since $\overline \K_+ =\K_+\cup I_+$,
$S_n$ are defined on $\P^k\setminus I_+$. Since $\dim I_+<k-p$, $S_n$ can
be extended to positive closed currents on $\P^k$ without mass on
$I_+$ \cite{HarveyPolking}. We also denote this extension 
by $S_n$. 
Since $f^n$ is an automorphism in
$\C^k$ we have $(f^n)^*(S_n)=d_+^{pn} S$ on $\C^k$. The equality
holds in $\P^k$ because the currents have supports in
$\overline\K_+$ and hence, have no mass at
infinity. So, necessarily $S_n$ have mass 1. Let $\Vc_{S_n}$, $\Vc_S$ denote the
dynamical super-potentials of $S_n$ and  $S$. We want to prove that
$S=T_+^p$. According to Proposition \ref{prop_unique_sp}, it is enough to show that
$\Vc_S=0$ on $\Cc_{k-p+1}(W)$ for any $W$ disjoint from $I_+$.

We have $L^n(S_n)=S$, hence Lemma \ref{lemma_dyn_sp_iterate_aut}
implies that $\Vc_S=d_+^{-n}\Vc_{S_n}\circ\Lambda^n$. Since
$\Vc_{S_n}$ is bounded from above on $\Cc_{k-p+1}(W)$ by a constant independent
of $n$, the last identity implies that $\Vc_S\leq 0$ on  $\Cc_{k-p+1}(W)$.
If $\Vc_S\not=0$ on $\Cc_{k-p+1}(W)$, there is a smooth form $R$ in
$\Cc_{k-p+1}(W)$
such that $\Vc_S(R)<0$. It follows that
$\Vc_{S_n}(\Lambda^n(R))\lesssim - d_+^n$. 
Let $W''$ be a neighbouhood of $\overline \K_+$, disjoint from $I_-$, such that
$f^{-1}(W'')\subset W''$. Hence, $\|Df^{-1}\|$ is bounded there by some
constant $M$ on $W''$. It follows that  $\|\Lambda^n(R)\|_{\infty, W''}\lesssim M^{3kn}$.
The inequality $\Vc_{S_n}(\Lambda^n(R))\lesssim - d_+^n$
contradicts
Proposition \ref{prop_estimate_sqp} which gives
$|\Vc_{S_n}(\Lambda^n(R))|\lesssim 1+\log M^{3kn}$. So, $\Vc_S=0$ on
$\Cc_{k-p+1}(W)$ and this completes the proof.
\hfill $\square$

\

The following result holds for currents of integration on generic
varieties of dimension $k-p$  in $\P^k$.

\begin{corollary}
Let $S$ be a current in $\Cc_p$ such that $\supp(S)\cap
I_-=\varnothing$. Then, $d_+^{-pn} (f^n)^*(S)$ converge to $T_+^p$.  
\end{corollary}
\proof
Let $W$ be a neighbourhood of $I_-$ such that $f(W)\Subset W$ and
$W\cap\supp(S)=\varnothing$. Hence, $f^{-n}(W)\subset f^{-n-1}(W)$ and
 $d_+^{-pn} (f^n)^*(S)$ has support in $\P^k\setminus f^{-n}(W)$. It follows that
 the limit values of $d_+^{-pn} (f^n)^*(S)$ are supported in
the complement of  $\cup_{n\geq 0} f^{-n}(W)$ which is contained in
 $\overline \K_+$. By Theorem \ref{th_aut_unique}, the only limit
 value is $T_+^p$.
\endproof

\begin{remark} \rm
In \cite{deThelin}, de Th{\'e}lin proved that the measure $\mu$ is
hyperbolic. It admits $k-p$ strictly negative and $p$ strictly
positive Lyapounov exponents. Pesin's theory implies that if a point
$a$ is  generic with respect to $\mu$, then it admits a stable manifold
of dimension $k-p$ and an unstable manifold of dimension $p$.
If $p=k-1$ and if $\tau:\C\rightarrow \overline \K_+$ is an entire
curve, using the Ahlfors' construction \cite{Ahlfors}, we obtain positive
closed $(k-1,k-1)$-currents with support in $\overline {\tau(\C)}$.
Indeed, Ahlfors' inequality implies the existence of
$(r_n)\rightarrow\infty$ such that the currents of integration on
$\tau(\Delta_{r_n})$, properly normalized, converge to a positive
closed current of mass 1. 
Theorem \ref{th_aut_unique} implies that this current is equal to
$T_+^{k-1}$. Hence $\overline {\tau(\C)}$ contains the support of $T_+^{k-1}$.
This result holds for generic stable manifolds of $\mu$. 
\end{remark}

\begin{remark} \rm
For $1\leq s\leq p$, if $S$ is a current in $\Cc_s$ with
super-potentials bounded on $\Cc_{k-s+1}(W)$ for some small
neighbourhood $W$ of $I_-$, then we can prove in the
same way that $d_+^{-sn} (f^n)^*(S)$ converge to $T_+^s$. 
The proof follows the same lines as the one in Theorem \ref{th_aut_unique}. We should choose
$W''$ large enough, in particular, we have $W''\cup W=\P^k$. 
In order to
apply Proposition \ref{prop_estimate_sqp}, we write $R$ as a
combination of a current in $\Cc_{k-p+1}(W)$ and  a smooth
form with bounded $\Cc^0$-norm. 
\end{remark}


T.-C. Dinh, UPMC Univ Paris 06, UMR 7586, Institut de
Math{\'e}matiques de Jussieu, F-75005 Paris, France. {\tt
  dinh@math.jussieu.fr}, {\tt http://www.math.jussieu.fr/$\sim$dinh} 

\

\noindent
N. Sibony,
Universit{\'e} Paris-Sud, Math{\'e}matique - B{\^a}timent 425, 91405
Orsay, France. {\tt nessim.sibony@math.u-psud.fr}

\end{document}